\documentstyle[amstex,psfig]{amsart}

 \newtheorem{thm}{Theorem}[section]
 
 \newtheorem{lem}[thm]{Lemma}
 
 \theoremstyle{definition}
 
 \theoremstyle{remark}
 
 \numberwithin{equation}{section}
 \newcommand{\set}[1]{\left\{#1\right\}}

\begin{document}

\title[Lower Order Terms]
{Lower order terms in the full moment conjecture for the Riemann zeta function}

\author[Conrey, Farmer, Keating, Rubinstein and Snaith]
{J.\ B.\ Conrey$^{1,2}$, D.\ W.\ Farmer$^1$, J.\ P.\ Keating$^2$,\\ M.\ O.\ Rubinstein$^3$ and N.\ C.\ Snaith$^{2}$\\ \\
1.\ American Institute of Mathematics\\ 360 Portage Avenue\\Palo Alto, CA 94306\\USA\\ \\
2.\ School of Mathematics \\ University of Bristol\\ Bristol BS8 1TW\\ UK \\ \\
3.\ Pure Mathematics \\ University of Waterloo\\ 200 University Ave W\\Waterloo, ON, Canada\\N2L 3G1\\}

%\author{Brian Corny, David Former, Jon Beating, Michael Ruination, Nina Snitch}

\subjclass{Primary L-functions, Random Matrix Theory}


%\date{2006}

%%% ----------------------------------------------------------------------

\begin{abstract}
We describe an algorithm for obtaining explicit expressions for lower terms
for the conjectured full asymptotics of the moments of the Riemann zeta function,
and give two distinct methods for obtaining numerical values of these coefficients.
We also provide some numerical evidence in favour of the conjecture.
\end{abstract}

%%% ----------------------------------------------------------------------
\maketitle
%%% ----------------------------------------------------------------------

\section{Introduction}

In \cite{CFKRS} the authors propose conjectures for the full asymptotics
of the moments of $L$-functions. A sample conjecture states, for
integer $k$, that the $2k$th moment of $|\zeta|$ on the half line
can be estimated using a polynomial $P_k$ of degree $k^2$, with
the polynomial given implicity as a $2k$-fold residue (see
(\ref{eq:P_k(x)}) below).

The leading term in the conjecture agrees with the Keating-Snaith 
conjecture [KS] for the leading asymptotics of the moments of $\zeta$.
Besides that, all the terms of the polynomial obtained agrees
with known theorems for $k=1,2$ \cite{I} \cite{H-B}, and the leading term agrees with
conjectures made through earlier and distinct number theoretical methods for $k=3,4$
\cite{CG} \cite{CGo}.

The method used in \cite{CFKRS} to conjecture the full asymptotics
relies on number theoretic heuristics based on the approximate
functional equation.  The conjecture is supported by the fact that the
formula coincides with an analogous expression in random matrix theory
\cite{CFKRS2} for the moments of characteristic polynomials from the
unitary group, the main difference being that the moments of $\zeta$
have extra arithmetic information that does not show up for random
matrices. Perhaps the most compelling support for the conjecture, though,
are numerics that confirm the conjectured moments.

For those numerical confirmations it is necessary to use all terms
arising in our heuristics.  The purpose of this paper is to investigate
the lower order degree terms which appear in the conjectured polynomials.
Specifically, we
\begin{enumerate}
    \item describe an algorithm to obtain meromorphic expressions in $k$ for the
    coefficients of the polynomial $P_k(x)$. Our main results are given in
    Theorems~\ref{thm:lot}--\ref{thm:b_k}.
    \item explain how one can numerically compute the lower order terms
    and to provide further experimental confirmation of the full moment conjecture, 
    including for non-integer values of $k$.
    Numerical values for the coefficients of these polynomials for $k=1,2,\ldots,7$ are listed
    in \cite{CFKRS} without explanation, with some numerics confirming the
    conjecture for $k=3$.
\end{enumerate}

At the end of the paper we also outline the analogous approach for moments of
quadratic Dirichlet $L$-functions
and of quadratic twists of an elliptic curve $L$-function, in both cases
evaluated at the critical point. These two cases are examples
of unitary-symplectic and orthogonal families respectively \cite{KaS} \cite{KeS2} \cite{CF}.

Before stating our results, we introduce notation and conjectures
from~\cite{CFKRS}.

\subsection{Moment conjecture for $\zeta$}
\label{subsection: moment conjectures zeta}

Let
$$
    \Delta(z_1,\ldots,z_{m})
    =\prod_{1\le i < j\le m}(z_j-z_i)
    = \left|  z_i^{j-1}  \right|_{m \times m}
$$
denote the Vandermonde determinant.

\proclaim Conjecture (see \cite{CFKRS}).

For positive integer $k$, and any $\epsilon > 0$,
\begin{equation}
\label{eq: moments zeta}
    \int_0^T |\zeta(1/2+it)|^{2k} dt= \int_0^T
    P_k\left(\log \tfrac{t}{2 \pi}\right) dt
    +O(T^{1/2+\epsilon}),
\end{equation}
with the constant in the $O$ term depending on $k$ and $\epsilon$,
where $P_k$ is the polynomial of degree $k^2$ given implicitly by the
$2k$-fold residue
\begin{equation}
\label{eq:P_k(x)}
     P_k(x)= \frac{(-1)^k}{k!^2}\frac{1}{(2\pi i)^{2k}}
    \oint\cdots \oint \frac{G(z_1,
    \ldots,z_{2k})\Delta^2(z_1,\ldots,z_{2k})} {\displaystyle
    \prod_{i=1}^{2k} z_i^{2k}} e^{\tfrac
    x2\sum_{i=1}^{k}z_i-z_{i+k}}~dz_1\ldots dz_{2k} ,
\end{equation}
with the path of integration over small circles about $z_i=0$, where
\begin{equation}
\label{eq:G}
    G(z_1,\ldots,z_{2k})= A_k(z_1,\ldots,z_{2k})
    \prod_{i=1}^k\prod_{j=1}^k\zeta(1+z_i-z_{j+k}) ,
\end{equation}
and $A_k$ is the Euler product
\begin{eqnarray}
  &&A_k(z_1,\ldots,z_{2k}) \notag \\
  &&=
  \label{eq:A_k}
  \prod_p
  \prod_{i,j=1}^k (1-p^{-1-z_i+z_{k+j}})
    %  \sum \begin{Sb} e_1+\cdots e_k= \\ f_1+\cdots f_k \\ e_i,f_i \geq 0 \end{Sb}
    %      p^{-\sum_1^k e_i(1/2+\alpha_i) + f_i(1/2-\alpha_{j+k}) }
    % diagonal sum above replaced with integral whose derivatives are more directly evaluated
    \int_0^1 \prod_{j=1}^k
    \left(1-\frac{e(\theta)}{p^{\frac12 +z_j}}\right)^{-1}
    \left(1-\frac{e(-\theta)}{p^{\frac12 -z_{k+j}}}\right)^{-1}\,d\theta
    \notag \\
    \\
    &&=
    \label{eq:A_k euler product}
    \prod_p \sum_{j=1}^k
    \prod_{i\neq j}
        \frac{\displaystyle \prod_{m=1}^k (1-p^{-1+z_{i+k}-z_m})}
        {1-p^{z_{i+k}-z_{j+k}}}.
\end{eqnarray}
Here $e(\theta) = \exp(2\pi i \theta)$.

We use both these expressions for the local factor for $A_k$. The first is used
in obtaining meromorphic expressions in $k$ for the coefficients of $P_k(x)$.

The second expression, derived in~\cite{CFKRS}[2.6],
is used to numerically compute $A_k(z_1,\ldots,z_{2k})$
for specific values of $z_1,\ldots,z_{2k}$.
The individual terms
in the sum over $j$ in~(\ref{eq:A_k euler product}) have poles
(though these poles cancel out when summed over $j$, see the paragraph 
following~\cite{CFKRS}[2.6.16])
and when we numerically
evaluate these terms individually, we take care to avoid the
poles by making sure that the $z_{j+k}$'s are distinct.

The main point of the conjecture is that we believe it gives the full
asymptotics of the moments of zeta.
While our numerical results in Section~\ref{section:data}
are consistent with a remainder of size $O(T^{1/2+\epsilon})$, there is some
debate regarding the error, especially in relation to moments of other
families of $L$-functions~\cite{CFKRS}~\cite{Z}, and
it would be worthwile to carry out more detailed testing concerning the nature of
the remainder.

The leading coefficient of $P_k(x)$ will be shown in Section~\ref{subsection:residue}
to equal
\begin{equation}
    \label{eq:g_k}
    a_k \prod_{j=0}^{k-1}\frac{j!}{(j+k)!},
\end{equation}
with 
\begin{equation}
    \label{eq:a_k}
    a_k= \prod_p \left(1-p^{-1}\right)^{k^2}
    {}_2F_1(k,k;1;1/p)
\end{equation}
and ${}_2F_1(a,b;c;t)$ the Gauss hypergeometric function
\begin{equation}
    {}_2F_1(a,b;c;t) =  \frac{\Gamma(c)}{\Gamma(a)\Gamma(b)} \sum_{n=0}^\infty
    \frac{\Gamma(a+n)\Gamma(b+n)}{\Gamma(c+n)}
     \ \frac{t^n}{n!}.
\end{equation}
This agrees with the leading term that was first conjectured by Keating and
Snaith~\cite{KeS}.

More generally, letting
$$
    Z(s) = \chi(s)^{-1/2} \zeta(s)
$$
with
$$
    \chi(s) = \pi^{s-1/2} \Gamma((1-s)/2) / \Gamma(s/2)
$$
we conjectured \cite{CFKRS} for shifted moments
\begin{equation}
\label{eq: shifted Z}
    \int\limits_{0}^{T} Z(1/2+it+\alpha_1)\cdots
    Z(1/2+it+\alpha_{2k})\,dt \sim 
    \int_{0}^{T}
    P_k\left(\alpha,\log \tfrac{t}{2 \pi}\right) ~dt,
\end{equation}
where
\begin{equation}
\label{eq:P_k(a,x)}
    P_k(\alpha,x)=
    \frac{(-1)^k}{k!^2}\frac{1}{(2\pi i)^{2k}}
    \oint \cdots \oint \frac{G(z_1,
    \ldots,z_{2k})\Delta(z_1,\ldots,z_{2k})^2} {\displaystyle
    \prod_{i=1}^{2k} \prod_{j=1}^{2k}(z_i-\alpha_j)} e^{\tfrac
    x2\sum_{i=1}^{k}z_i-z_{i+k}}~dz_1\ldots dz_{2k},
\end{equation} with the path of integration being small circles
surrounding the poles $\alpha_i$, and $-1/4 < \Re \alpha_j$.
%We believe that the difference between the l.h.s. and r.h.s.
%of~(\ref{eq: shifted Z}) is smaller in size than any of the terms on the
%right, i.e. that~(\ref{eq: shifted Z}) captures the full asymptotics
%of the shifted moments.
One recovers the moments of
$\zeta$ by setting the shifts $\alpha_i$ equal to $0$, and
observing that $Z(s)=Z(1-s)$. 

An alternative formulation of this conjecture also given in
\cite{CFKRS} involves a combinatorial sum and is established by the 
following lemma.

\begin{lem}

See \cite{CFKRS}, Section 2.5.

Suppose $F(u;v)=F(u_1,\ldots,u_k;v_1,\ldots,v_{k})$ is a function
of $2k$ variables, symmetric with respect to the first
$k$ variables and also symmetric with respect to the second set of
$k$ variables. Suppose also that $F$ is regular near
$(0,\ldots,0)$, and that $f(s)$ has a simple pole of
residue~$1$ at $s=0$ but is otherwise analytic in a neighbourhood
about $s=0$. Let
$$
    H(u_1,\ldots,u_k;v_1,\ldots v_k)= F(u_1,\ldots;\ldots,v_k)
    \prod_{i=1}^k\prod_{j=1}^k f(u_i-v_j).
$$
If for all $1\leq i,j \leq k$, $\alpha_i-\alpha_{j+k}$ is
contained in the region of analyticity of $f(s)$ then
\begin{eqnarray}
    \label{eq: combinatorial sum}
    &
    \frac{(-1)^k}{k!^2}\frac{1}{(2\pi i)^{2k}}
    \oint \cdots \oint \frac{H(z_1, \ldots,z_{2k})
    \Delta(z_1,\ldots,z_{2k})^2}{\prod_{i=1}^{2k}\prod_{j=1}^{2k}
    (z_i-\alpha_j)} \,dz_1\ldots dz_{2k}
    \notag \\
    =&\sum_{\sigma \in \Xi} H( \alpha_{\sigma(1)},\ldots,\alpha_{\sigma(2k)}),
\end{eqnarray}
where one integrates about small circles enclosing the
$\alpha_j$'s, and where $\Xi$ is the set of $\binom{2k}{k}$
permutations $\sigma\in S_{2k}$ such that $\sigma(1)<\cdots <
\sigma(k)$ and $\sigma(k+1)<\cdots < \sigma(2k)$.
\end{lem}

Equation (\ref{eq:P_k(x)}) allows us to obtain the coefficients of $P_k(x)$
by computing power series expansions and then the residue of the r.h.s,
giving meromorphic expressions in $k$ for the coefficients which
can also be evaluated to high precision numerically, even for non integer $k$. 
In practice we have been able to do so for the first ten coefficients of $P_k(x)$.
If $k \in {\mathbb Z}$, $k>3$, to obtain numerical values
for all $k^2$ coefficients of $P_k(x)$ we developed a second method using
equation~(\ref{eq: combinatorial sum}).
This involved taking small distinct shifts and high
working precision to capture
cancellation amongst the order $k^2$ poles of the r.h.s. of ~(\ref{eq: combinatorial sum}).

\subsection{Results}
\label{section: results}

Our first theorem below explictly gives the coefficients of
$P_k(x)$ in the full moment conjecture for the Riemann zeta function.
These are described in terms of the multivariate Taylor coefficients of
\begin{equation}
    \label{eq:series of}
    \frac{1}{a_k} A_k(z_1,\ldots,z_{2k}) \prod_{1 \leq i,j \leq k} (z_i-z_{j+k}) \zeta(1+z_i-z_{j+k}).
\end{equation}
We let $b_k(\alpha;\beta)$ denote the coefficient of
$z_1^{\alpha_1} \ldots z_k^{\alpha_k} z_{k+1}^{\beta_1} \ldots z_{2k}^{\beta_k}$
in the Taylor series of~(\ref{eq:series of}).
Let
\begin{equation}
    |\alpha|= \sum_1^k \alpha_i
\end{equation}
and likewise for $\beta$.
Notice that the function in~(\ref{eq:series of}) is symmetric
in $z_1,\ldots,z_{k}$ and in $z_{k+1},\ldots,z_{2k}$, and that 
$b_k(\alpha;\beta) = (-1)^{|\alpha|+|\beta|} b_k(\beta;\alpha)$. We may therefore collect
together terms in the Taylor series accordingly and express~(\ref{eq:series of}) as
\begin{equation}
    \sum_{\alpha;\beta}
    b_k(\alpha;\beta)
    \left(
    z_1^{\alpha_1} \cdots z_k^{\alpha_k}
    z_{k+1}^{\beta_1} \cdots z_{2k}^{\beta_k} \pm \text{sym}
    \right)
\end{equation}
The `sym' indicates that we group terms that have exponents of the same
form. Thus the sum over $\alpha;\beta$ follows the convention
that $\alpha \geq \beta$ lexicographically, we list the $\alpha_i$'s and $\beta_i$'s in
decreasing order, and we supress the $\alpha_i$'s and $\beta_i$'s that are 0.
For example, all the terms of degree $4$ are collected with
coefficients: $b_k(1,1,1,1;)$, $b_k(2,1,1;)$, $b_k(2,2;)$,
$b_k(3,1;)$, $b_k(4;)$, $b_k(1,1,1;1)$, $b_k(2,1;1)$, $b_k(3;1)$, $b_k(1,1;1,1)$, $b_k(2;1,1)$,
$b_k(2;2)$. The terms that go with $b_k(1,1,1,1;)$ are
$$\sum_{1\leq i_1 < i_2 <i_3 < i_4 \leq k}
z_{i_1} z_{i_2} z_{i_3} z_{i_4} +z_{k+i_1} z_{k+i_2} z_{k+i_3} z_{k+i_4}.$$

\begin{thm}
\label{thm:lot}
Let $P_k(x)$ be given by equation (\ref{eq:P_k(x)}).
Writing
\begin{equation}
    P_k(x)= c_0(k) x^{k^2} + c_1(k) x^{k^2-1}+ \ldots + c_{k^2}(k),
\end{equation}
we have
\begin{equation}
    \boxed{
    c_r(k) = a_k \prod_{l=0}^{k-1} \frac{l!}{(k+l)!}
    \sum_{|\alpha|+ |\beta| =r}
    2^{1-\delta(\alpha;\beta)}
    b_k(\alpha;\beta) N_k(\alpha;\beta),
    }
\end{equation}
where $a_k$ is given by (\ref{eq:a_k}), the function $\delta(\alpha;\beta)$
equals zero unless $\alpha=\beta$ in which case it equals 1, and
$N_k(\alpha;\beta)$ is defined by
\begin{equation}
    \label{eq:N thm}
    N_k(\alpha;\beta) =
    \frac{1}{2^{k^2-r}}
    \left( \prod_{l=0}^{k-1} \frac{l!}{(k+l)!} \right)^{-1}
    \sum_{{\text{rearrangements}} \atop {\text{$\sigma,\tau$ of $\alpha$ and $\beta$}}}
    \tilde{M}_k(\sigma(\alpha),\tau(\beta)).
\end{equation}
The function $\tilde{M}_k$ is given as a $2k \times 2k$ determinant
in equation (\ref{eq:M tilde}).
\end{thm}

By rearrangements, we mean distinct permutations.
Two permuatations $\sigma(\alpha)$ and $\mu(\alpha)$ are said to 
be distinct if $\alpha_{\sigma_i} \neq \alpha_{\mu_i}$ for some $i$. 
For example, if $\alpha_1=7, \alpha_2=5, \alpha_3=5$ then the two permutations
$\alpha_2, \alpha_1, \alpha_3$ and $\alpha_3, \alpha_1, \alpha_2$ are not distinct
and would only be counted once in~(\ref{eq:N thm}).

The reason for writing $N_k(\alpha;\beta)$ as we 
have, with the factor
$\left( \prod_{l=0}^{k-1} \frac{l!}{(k+l)!} \right)^{-1}$, is explained by the next theorem.
\begin{thm}
    \label{thm:N}
    $N_k(\alpha;\beta)$ is a polynomial in $k$ of degree $\leq 2(|\alpha|+|\beta|)$.
\end{thm}
This theorem allows us to determine $N_k(\alpha;\beta)$ explicitly by
evaluating~(\ref{eq:N thm}) at $2(|\alpha|+|\beta|)+1$ values of $k$ and
interpolating. A few example $N_k(\alpha;\beta)$'s are given in~(\ref{eq:N exs}).

Finally, the coefficients $b_k(\alpha;\beta)$ that appear in Theorem~\ref{thm:lot}
can also be explicitly determined.
\begin{thm}
    \label{thm:b_k}
    The Taylor coefficients $b_k(\alpha;\beta)$ of~(\ref{eq:series of}) can be written
    explicitly as a polynomial in: $k$, the Taylor coefficients of
    $s\zeta(1+s)$, and the Taylor coefficients of $\log(A_k(z_1,\ldots,z_{2k}))$.
    The latter Taylor coefficients can further be expressed explicitly as a sum over all
    primes $p$ of a rational function in: $k$, $p$, $\log(p)$, and finitely many Gauss
    hypergeomtric functions ${}_2F_1(k_1,k_2;m;1/p)$, where $k_1,k_2,m \in {\mathbb Z}$,
    $k_1,k_2 \geq k$ and $m \geq 1$.
\end{thm}
To illustrate what the last theorem looks like in practice, see
equations~(\ref{eq: b_k series}) and~(\ref{eq:B_k}).

This paper is structured as follows. In Section~\ref{subsection:residue}
we prove Theorems~\ref{thm:lot} and~\ref{thm:b_k}, and also give a procedure to determine
the polynomial and rational functions of Theorem~\ref{thm:b_k}.
In Section~\ref{section:poly} we prove Theorem~\ref{thm:N}.

Section~\ref{section: evaluating c_r} is devoted to numerical
evaluation of the lower order terms. Two different methods are
described. The first involves numerically computing the terms that appear in
in Theorem~\ref{thm:lot}, while the second
uses~(\ref{eq: combinatorial sum}), small shifts, and very high precision
to capture cancellation amongst the poles of the summand.
Data supporting the full moment conjecture is then presented in
Section~\ref{section:data}.

We also provide plots of the coefficients $c_r(k)$ as a function of $k$ for
$r \leq 7$ and also of the zeros of the polynomials $P_k(x)$ for several
values of $k$.

In Section~\ref{section:other families} we briefly describe
the analogous approach for quadratic Dirichlet $L$-functions and of
quadratic twists of an elliptic curve $L$-function.

\section{Lower order terms in the moments of $\zeta$}
\label{section: lower terms zeta}

\subsection{Evaluating the residue explicitly}
\label{subsection:residue}

The $2k$-fold residue in (\ref{eq:P_k(x)}) involves extracting the
coefficient of $\prod_1^{2k} z_i^{2k-1}$, i.e. a polynomial of
degree $2k(2k-1)$. The Vandermonde determinant has degree
$2{2k\choose2} = 2k(2k-1)$. However, the product
$\prod_{i=1}^k\prod_{j=1}^k\zeta(1+z_i-z_{j+k})$ has poles which
cancel $k^2$ of the Vandermonde factors. Hence, in
(\ref{eq:P_k(x)}), we need only take terms in the Taylor expansion
of $\exp\left(\frac{x}{2}\sum_{i=1}^k z_i-z_{i+k}\right)$ up to
degree $k^2$. Therefore, $P_k(x)$ is indeed a polynomial of degree
$k^2$ in $x$ and we write
\begin{equation}
\label{eq: coeffs P}
    P_k(x) = c_0(k) x^{k^2} + c_1(k) x^{k^2-1} + c_2(k) x^{k^2-2}
    + \cdots + c_{k^2}(k).
\end{equation}

One complication in developing expressions in $k$ for the
$c_r(k)$'s is that the Vandermonde determinant
in~(\ref{eq:P_k(x)}) prevents us from separating the integrals.
However, this can be overcome by introducing extra variables
and pulling out the Vandermonde as a differential operator.
We illustrate the method for the leading term $c_0(k)$ and then
generalize.

Noticing that $A_k(0,0,\ldots,0) = a_k$ (set all the variables to $0$ 
in~(\ref{eq:A_k}), apply Lemma~\ref{lemma:hyper}, and compare to~(\ref{eq:a_k})),
the leading term is given by
\begin{eqnarray}
\label{eq:c_0 a}
    c_0(k) x^{k^2} = \frac{a_k}{k!^2}
    \frac{1}{(2\pi i)^{2k}}
    \oint\cdots \oint
       &\Delta(z_1,\ldots,z_{2k})
        \Delta(z_1,\ldots,z_{k})
        \Delta(z_{k+1},\ldots,z_{2k}) \notag \\
       &\times
        \frac{
                \exp\left({\tfrac
                x2\sum_{i=1}^{k}z_i-z_{i+k}}\right)
             }
             {
                 {\displaystyle
                 \prod_{i=1}^{2k} z_i^{2k}}
             }
        dz_1\ldots dz_{2k}
\end{eqnarray}
(the $k^2$ poles of the $\zeta$ product have sign opposite from
the Vandermonde factors that they cancel, and these $k^2$ minuses
cancel the $(-1)^k$ in (\ref{eq:P_k(x)})). Comparing the degree of the Vandermonde
determinants in the numerator, with the degree of the denominator shows that only
terms of degree $k^2$ in the taylor expansion of the $\exp$ contribute to the 
residue.
Changing variables
$u_i=x z_i/2$ and then relabeling $u_i$ with $z_i$ gives
\begin{eqnarray}
\label{eq:c_0 b}
    \frac{x^{k^2}}{2^{k^2}} \frac{a_k}{k!^2}
    \frac{1}{(2\pi i)^{2k}}
    \oint\cdots \oint
       &\Delta(z_1,\ldots,z_{2k})
        \Delta(z_1,\ldots,z_{k})
        \Delta(z_{k+1},\ldots,z_{2k}) \notag \\
       &\times
        \frac{
                \exp\left({\sum_{i=1}^{k}z_i-z_{i+k}}\right)
             }
             {
                 {\displaystyle
                 \prod_{i=1}^{2k} z_i^{2k}}
             }
        dz_1\ldots dz_{2k}.
\end{eqnarray}
Introducing extra variables $x_i$, we consider
\begin{equation}
\label{eq:trick a}
    \frac{1}{(2\pi i)^{2k}}
    \oint\cdots \oint
        p(z_1,\ldots,z_{2k})
        \frac{
                \exp\left({\sum_{i=1}^{2k}x_i z_i}\right)
             }
             {
                 {\displaystyle
                 \prod_{i=1}^{2k} z_i^{2k}}
             }
        dz_1\ldots dz_{2k}.
\end{equation}
with $p(z)$ a polynomial in $z_1,\ldots,z_{2k}$. Pulling out the
polynomial $p$, (\ref{eq:trick a}) equals
\begin{equation}
\label{eq:trick b}
    p(\partial/\partial x_1,\ldots,\partial/\partial x_{2k})
    \frac{1}{(2\pi i)^{2k}}
    \oint\cdots \oint
        \frac{
                \exp\left({\sum_{i=1}^{2k}x_i z_i}\right)
             }
             {
                 {\displaystyle
                 \prod_{i=1}^{2k} z_i^{2k}}
             }
        dz_1\ldots dz_{2k},
\end{equation}
and taking the residue gives
\begin{equation}
\label{eq:trick c}
    p(\partial/\partial x_1,\ldots,\partial/\partial x_{2k})
    \prod_1^{2k}\frac{x_i^{2k-1}}{(2k-1)!}.
\end{equation}
Therefore, (\ref{eq:c_0 b}) equals
\begin{equation}
\label{eq:c_0 c}
    \frac{x^{k^2}}{2^{k^2}} \frac{a_k}{k!^2}
    q(\partial/\partial x_1,\ldots,\partial/\partial x_{2k})
    \prod_1^{2k}\frac{x_i^{2k-1}}{(2k-1)!}
\end{equation}
evaluated at $x_1=\ldots=x_k=1$, $x_{k+1}=\ldots=x_{2k}=-1$, with
\begin{equation}
\label{eq:triple vandermonde}
    q(\partial/\partial x_1,\ldots,\partial/\partial x_{2k})=
    \Delta(\partial/\partial x_1,\ldots,\partial/\partial x_{2k})
    \Delta(\partial/\partial x_1,\ldots,\partial/\partial x_{k})
    \Delta(\partial/\partial x_{k+1},\ldots,\partial/\partial x_{2k}).
\end{equation}
Two lemmas allow us to reduce this further.
\begin{lem}
\label{lemma:single vand}
$$
    \Delta(\partial/\partial x_1,\ldots,\partial/\partial x_{n})
    \prod_1^n f(x_i)
    = \left| f^{(j-1)}(x_i) \right|_{n \times n.}
$$
\end{lem}
\begin{pf}
This follows using the definition of the Vandermonde determinant
$$
    \Delta(\partial/\partial x_1,\ldots,\partial/\partial x_{n})
    = \left| \partial^{j-1}/\partial x_i^{j-1} \right|_{n \times n.}
$$
Noticing that row $i$ of the matrix only involves $x_i$, we factor
the product into the determinant.
\end{pf}

We can now consider the effect of applying the three Vandermonde's
in (\ref{eq:triple vandermonde}).

\begin{lem}
\label{lemma:triple vand}
$$
    \Delta(\partial/\partial x_1,\ldots,\partial/\partial x_{2k})
    \Delta(\partial/\partial x_1,\ldots,\partial/\partial x_{k})
    \Delta(\partial/\partial x_{k+1},\ldots,\partial/\partial x_{2k})
     \prod_1^{2k} f(x_i)
$$
evaluated at $x_1=\ldots=x_k=1$, $x_{k+1}=\ldots=x_{2k}=-1$ equals
$$
    k!^2 \left|
             \begin{array}{cccc}
             f(1) & f^{(1)}(1) & \ldots & f^{(2k-1)}(1) \\
             f^{(1)}(1) & f^{(2)}(1) & \ldots & f^{(2k)}(1) \\
             \vdots & \vdots & \ddots & \vdots \\
             f^{(k-1)}(1) & f^{(k)}(1) & \ldots & f^{(3k-2)}(1) \\
             f(-1) & f^{(1)}(-1) & \ldots & f^{(2k-1)}(-1) \\
             f^{(1)}(-1) & f^{(2)}(-1) & \ldots & f^{(2k)}(-1) \\
             \vdots & \vdots & \ddots & \vdots \\
             f^{(k-1)}(-1) & f^{(k)}(-1) & \ldots & f^{(3k-2)}(-1)
             \end{array}
         \right|_{2k \times 2k.}
$$
(the first $k$ rows of this  $2k \times 2k$ matrix involve $f$ and
its derivatives evaluated at $1$, while the last $k$ rows have the
entries evaluated at $-1$).
\end{lem}

\begin{pf}
Consider first what happens when we apply just the last two
$\Delta$'s. By Lemma~\ref{lemma:single vand},
\begin{eqnarray}
    &&\Delta(\partial/\partial x_1,\ldots,\partial/\partial x_{k})
    \Delta(\partial/\partial x_{k+1},\ldots,\partial/\partial x_{2k})
    \prod_1^{2k} f(x_i) \notag \\
    &=&
    \left|
        \begin{array}{cc}
            A_{k\times k} & 0_{k\times k} \\
            0_{k\times k} & B_{k\times k} \\
        \end{array}
    \right|
\end{eqnarray}
with
\begin{eqnarray}
\label{eq:A}
    A_{i,j} &=&  f^{(j-1)}(x_i) \\
    B_{i,j} &=&  f^{(j-1)}(x_{i+k}).
\end{eqnarray}
Expanding this determinant, we get a sum of $(k!)^2$ terms each of
which is a product of the form
\begin{equation}
\label{eq:expand the det}
    \text{sgn}(\mu) \ \text{sgn}(\nu) \prod_{i=1}^k f^{(\mu_i-1)}(x_i) f^{(\nu_i-1)}(x_{i+k})
\end{equation}
where $\mu_1,\ldots,\mu_k$ and $\nu_1,\ldots,\nu_k$ are
permutations of the numbers $1,2,\ldots,k$. Applying the third
Vandermonde $\Delta(\partial/\partial x_1,\ldots,\partial/\partial
x_{2k})$ to a typical such term gives, by Lemma~\ref{lemma:single vand},
\begin{equation}
\label{eq:typical term}
    \text{sgn}(\mu) \ \text{sgn}(\nu)
    \left|
             \begin{array}{cccc}
             f^{(\mu_1-1)}(x_1) & f^{(\mu_1)}(x_1) & \ldots & f^{(\mu_1+2k-2)}(x_1) \\
             f^{(\mu_2-1)}(x_2) & f^{(\mu_2)}(x_2) & \ldots & f^{(\mu_2+2k-2)}(x_2) \\
             \vdots & \vdots & \ddots & \vdots \\
             f^{(\mu_k-1)}(x_k) & f^{(\mu_k)}(x_k) & \ldots & f^{(\mu_k+2k-2)}(x_k) \\
             f^{(\nu_1-1)}(x_{k+1}) & f^{(\nu_1)}(x_{k+1}) & \ldots & f^{(\nu_1+2k-2)}(x_{k+1}) \\
             f^{(\nu_2-1)}(x_{k+2}) & f^{(\nu_2)}(x_{k+2}) & \ldots & f^{(\nu_2+2k-2)}(x_{k+2}) \\
             \vdots & \vdots & \ddots & \vdots \\
             f^{(\nu_k-1)}(x_{2k}) & f^{(\nu_k)}(x_{2k}) & \ldots & f^{(\nu_k+2k-2)}(x_{2k})
             \end{array}
    \right|_{2k \times 2k.}
\end{equation}
Setting $x_1=\ldots=x_k=1$, $x_{k+1}=\ldots=x_{2k}=-1$, we may
rearrange the first $k$ rows and the last $k$ rows so as to undo
the permutations $\mu$ and $\nu$. This introduces another
$\text{sgn}(\mu) \ \text{sgn}(\nu)$ in front of the the
determinant. Hence each such term contributes the same amount,
\begin{equation}
\label{eq:typical contribution}
       \left|
             \begin{array}{cccc}
             f(1) & f^{(1)}(1) & \ldots & f^{(2k-1)}(1) \\
             f^{(1)}(1) & f^{(2)}(1) & \ldots & f^{(2k)}(1) \\
             \vdots & \vdots & \ddots & \vdots \\
             f^{(k-1)}(1) & f^{(k)}(1) & \ldots & f^{(3k-2)}(1) \\
             f(-1) & f^{(1)}(-1) & \ldots & f^{(2k-1)}(-1) \\
             f^{(1)}(-1) & f^{(2)}(-1) & \ldots & f^{(2k)}(-1) \\
             \vdots & \vdots & \ddots & \vdots \\
             f^{(k-1)}(-1) & f^{(k)}(-1) & \ldots & f^{(3k-2)}(-1)
             \end{array}
         \right|_{2k \times 2k,}
\end{equation}
and summing over the $(k!)^2$ pairs $\mu,\nu$ gives us the Lemma.
\end{pf}
Applying Lemma~\ref{lemma:triple vand} to (\ref{eq:c_0 c}) yields
\begin{equation}
\label{eq:c_0 d}
    c_0(k) =
    \frac{a_k}{2^{k^2}}
           \left|
             \begin{array}{cccc}
             \Gamma(2k)^{-1} &  \Gamma(2k-1)^{-1}& \ldots & \Gamma(1)^{-1} \\
             \Gamma(2k-1)^{-1} & \Gamma(2k-2)^{-1} & \ldots & \Gamma(0)^{-1} \\
             \vdots & \vdots & \ddots & \vdots \\
             \Gamma(k+1)^{-1} & \Gamma(k)^{-1} & \ldots & \Gamma(-k+2)^{-1} \\
             -\Gamma(2k)^{-1} & \Gamma(2k-1)^{-1} & \ldots &  \Gamma(1)^{-1}\\
             \Gamma(2k-1)^{-1} & -\Gamma(2k-2)^{-1} & \ldots & -\Gamma(0)^{-1} \\
             \vdots & \vdots & \ddots & \vdots \\
             (-1)^{k}\Gamma(k+1)^{-1} & (-1)^{k+1}\Gamma(k)^{-1} & \ldots & (-1)^{3k-1}\Gamma(-k+2)^{-1}
             \end{array}
         \right|_{2k \times 2k.} \ \ \ \
\end{equation}
Here, we take $\Gamma(m)^{-1} = 0$ if $m \in
\set{0,-1,-2,-3,\ldots}$. The first $k$ rows and last $k$ rows are
identical except for the presence of a checkerboard pattern of
minus ones in the latter rows. The $i,j$ entry above equals
\begin{equation}
     \begin{cases}
     \Gamma(2k-i-j+2)^{-1}, &\text{if $1 \leq i \leq k$;} \\
     (-1)^{i+j-k-1}\Gamma(3k-i-j+2)^{-1}&\text{if $k+1 \leq i \leq 2k$;}
    \end{cases}
\end{equation}
We show later that (\ref{eq:c_0 d}) equals $a_k \prod_{l=0}^{k-1} l!/(k+l)!$.
See Lemmas \ref{lemma: poly1} and \ref{thm: poly2}, with
$e_i=f_i=0$, $c_i=k+i-1$.

Next we consider in (\ref{eq: coeffs P}) the $r$th term of our
polynomial $P_k(x)$. To evaluate $c_r(k)$ we examine the power
series expansion of the integrand in~(\ref{eq:P_k(x)}). As in our
consideration of $c_0(k)$, we first cancel the poles of
$\prod_{i=1}^k\prod_{j=1}^k\zeta(1+z_i-z_{j+k})$ against the
Vandermonde, and write the integral in~(\ref{eq:P_k(x)}) as
\begin{eqnarray}
\label{eq:huge integral}
    P_k(x)= &&\frac{1}{k!^2}\frac{1}{(2\pi i)^{2k}}
    \oint\cdots \oint
    \frac{
        \Delta(z_1,\ldots,z_{2k})
        \Delta(z_1,\ldots,z_{k})
        \Delta(z_{k+1},\ldots,z_{2k})
     }
     {
        \displaystyle \prod_{i=1}^{2k} z_i^{2k}
     }
    \notag
    \\
     \times
     &&A_k(z_1,\ldots,z_{2k})
    \prod_{1 \leq i,j \leq k}
    (z_i-z_{j+k}) \zeta(1+z_i-z_{j+k})
    \notag \\
    \times
    &&e^{\tfrac x2\sum_{i=1}^{k}z_i-z_{i+k}}
    dz_1\ldots dz_{2k}.
\end{eqnarray}

Because of the various symmetries satistfied by the factors of the
integrand, our job of determining the series expansion of
\begin{equation}
\label{eq:A_k and product}
    A_k(z_1,\ldots,z_{2k})
    \prod_{1 \leq i,j \leq k}
    (z_i-z_{j+k}) \zeta(1+z_i-z_{j+k})
\end{equation}
is not as difficult as might be supposed.

\subsubsection{Series for $A_k(z_1,\ldots,z_{2k})$}
\label{section:series A_k} First, we write
\begin{equation}
\label{eq:log A_k}
    A_k(z_1,\ldots,z_{2k}) =
    \exp\left( \log(A_k(z_1,\ldots,z_{2k})\right)
\end{equation}
to turn the Euler product defining $A_k$ in~(\ref{eq:A_k euler
product}) into a sum over primes. Obtaining the series for
$\log(A_k(z_1,\ldots,z_{2k}))$ will allow us, in conjunction with
the method in Section~\ref{section:multiplying} for multiplying
series, to recover the series for
$\exp\left(\log(A_k(z_1,\ldots,z_{2k}))\right)$

Because $A_k(z_1,\ldots,z_{2k})$ is symmetric in $z_1,\ldots,z_k$
and separately in $z_{k+1},\ldots,z_{2k}$, we distinguish two sets
of variables, and let
\begin{equation}
\label{eq:B}
    B_k(\alpha_1,\alpha_2,\ldots,\alpha_k;
    \beta_1,\beta_2,\ldots,\beta_k)
\end{equation}
denote the coefficient of a typical
\begin{equation}
\label{eq:typical taylor term}
    \frac{z_1^{\alpha_1}}{\alpha_1!}
    \cdots
    \frac{z_k^{\alpha_k}}{\alpha_k!}
    \frac{z_{k+1}^{\beta_1}}{\beta_1!}
    \cdots
    \frac{z_{2k}^{\beta_k}}{\beta_k!}
\end{equation}
in the multivariate Taylor expansion of
$\log(A_k(z_1,\ldots,z_{2k}))$. Here, we prefer
to keep factorials in the denominator rather than
absorb them into $B_k$ for convenience in describing
the procedure to obtain the coefficients through
differentiation.

By the above mentioned symmetry,
we use the convention, when writing $B_k(\alpha;\beta)$ of only
listing the non-zero $\alpha_i$'s and $\beta_i$'s, and writing
them in decreasing order. Also, because
\begin{equation}
\label{eq:A_k second symmetry}
    A_k(z_1,\ldots,z_{2k}) = A_k(-z_{k+1},\ldots,-z_{2k},-z_1,\ldots,-z_{k+1})
\end{equation}
we have
\begin{equation}
\label{eq:B symmetry}
    B_k(\alpha;\beta) = (-1)^{|\alpha|+|\beta|} B_k(\beta;\alpha)
\end{equation}
where $|\alpha| = \sum_1^k \alpha_i$. Therefore
\begin{eqnarray}
\label{eq:A_k exp series}
    &&A_k(z_1,\ldots,z_{2k}) = \notag \\
    &&a_k \exp \biggl(
        B_k(1;) \sum_1^k z_i-z_{i+k} +
        B_k(1,1;) \sum_{1\leq i<j\leq k}
            z_i z_j + z_{i+k}z_{j+k} \notag \\
     &&+B_k(1;1) \sum_{1\leq i,j\leq k}
            z_i z_{j+k}
       +B_k(2;) \sum_1^{2k} \frac{z_i^2}{2} \notag \\
     &&+B_k(1,1,1;) \sum_{1\leq i<j<l \leq k}
         z_i z_j z_l - z_{i+k}z_{j+k}z_{l+k} + \cdots \notag
    \biggr). \\
\end{eqnarray}
The $a_k$ factor comes from the value of the function at the
origin $A_k(0,\ldots,0) = a_k$.

Next, let $l=l(\alpha)$ denote the number of non-zero
$\alpha_i$'s, and $m=m(\beta)$ denote the number of non-zero
$\beta_i$'s. Since we are interested in extracting $c_r(k)$, we
only need to consider the power series expansion of $A_k$ up to
degree $r$, i.e. $l+m \leq r$. Since we are assuming in our
evaluation of $B_k(\alpha;\beta)$ that the $\alpha_i$'s and
$\beta_i$'s are in decreasing order, we focus on the first $l$
$z_i$'s ($m$ $z_{i+k}$'s respectively) and we have, together
with~(\ref{eq:A_k euler product}),
\begin{equation}
\label{eq:B as sum over primes}
    B_k(\alpha;\beta) =
    \sum_p
    \prod_{i=1}^l
    \frac{\partial^{\alpha_i}}{\partial z_i^{\alpha_i}}
    \prod_{i=1}^m
    \frac{\partial^{\beta_i}}{\partial z_{i+k}^{\beta_i}}
    \log\left( f_k(1/p;z) \right) \biggr|_{z=0}
\end{equation}
where
\begin{equation}
\label{eq:f_k}
    f_k(t;z) =
    \prod_{1\leq i,j \leq k} (1-t^{1+z_i-z_{k+j}})
    %\sum \begin{Sb} e_1+\cdots e_k= \\ f_1+\cdots f_k \end{Sb}
    %  t^{\sum_1^k e_i(1/2+z_i) + f_i(1/2-z_{i+k})}
    %diagonal sum replaced by integral expression
    \int_0^1 \prod_{j=1}^k
    \left(1-e(\theta)t^{\frac12 +z_j}\right)^{-1}
    \left(1-e(-\theta)t^{\frac12 -z_{k+j}}\right)^{-1}\,d\theta
\end{equation}
Since we are assuming $\alpha_{l+1}=\ldots=\alpha_{k}=0$,
$\beta_{m+1}=\ldots=\beta_{k}=0$, with $l+m \leq r$, we may as
well immediately set $z_{r+1}=\ldots=z_{k}=0$,
$z_{k+r+1}=\ldots=z_{2k}=0$. Therefore,
\begin{eqnarray}
\label{eq:log first piece}
    \sum_{1\leq i,j \leq k}
    &\log\left( 1-t^{1+z_i-z_{k+j}} \right)
    \equiv
    \displaystyle\sum_{1\leq i,j \leq r}
    \log\left( 1-t^{1+z_i-z_{k+j}} \right)
    \notag \\
    &+
    \displaystyle\sum_{i=1}^r (k-r) \left( \log(1-t^{1+z_i}) +
    \log(1-t^{1-z_{i+k}})\right).
\end{eqnarray}
By equivalent, we mean that both expressions have the same series
expansion in $z$ up to terms involving just $z_1,\ldots,z_r$, though
not including the constant which, on the l.h.s., equals $k^2 \log(1-t)$.
The main point in doing this reduction is to get rid of the $k$
dependence in the summands.

A symbolic differentiation package (such as Maple) can then be
used to compute
\begin{equation}
\label{eq:partial diff operator}
    \prod_{i=1}^l
    \frac{\partial^{\alpha_i}}{\partial z_i^{\alpha_i}}
    \prod_{i=1}^m
    \frac{\partial^{\beta_i}}{\partial z_{i+k}^{\beta_i}}
    \biggr|_{z=0}
\end{equation}
applied to the r.h.s. of~(\ref{eq:log first piece}) as a rational
function in $k$, $\log(t)$, and $t$. We list the terms up to degree $2$:
\begin{eqnarray}
\label{eq:series log first piece}
    &&\sum_{1\leq i,j \leq k}
    \log\left( 1-t^{1+z_i-z_{k+j}} \right)
    =
    k^2 \log(1-t) -
    \frac{kt\log(t)}{1-t} \sum_1^k z_i-z_{i+k} \notag \\
    %+0 \sum_{1\leq i<j\leq k} z_i z_j + z_{i+k}z_{j+k} \notag \\
    &&+\frac{t \log(t)^2}{(1-t)^2} \sum_{1\leq i,j\leq k} z_i z_{j+k}
    -\frac{kt\log(t)^2}{(1-t)^2} \sum_1^{2k} \frac{z_i^2}{2} + \ldots
\end{eqnarray}
The coefficient above of the $\sum_{1\leq i<j\leq k} z_i z_j + z_{i+k}z_{j+k}$ term equals
zero.

Next, applying~(\ref{eq:partial diff operator}) to
\begin{equation}
\label{eq:log second piece}
    \log\left(
         %\sum \begin{Sb} e_1+\cdots e_k= \\ f_1+\cdots f_k \end{Sb}
         %  t^{\sum_1^k e_i(1/2+z_i) + f_i(1/2-z_{i+k})}
         %diagonal sum replaced by integral expression
         \int_0^1 \prod_{j=1}^k
         \left(1-e(\theta)t^{\frac12 +z_j}\right)^{-1}
         \left(1-e(-\theta)t^{\frac12 -z_{k+j}}\right)^{-1}\,d\theta
    \right)
\end{equation}
we end up, by the chain rule, with a rational expression involving partial derivatives
of the form
\begin{eqnarray}
\label{eq:partial diff int}
    \prod_{i=1}
    \frac{\partial^{c_i}}{\partial z_i^{c_i}}
    \prod_{i=1}
    \frac{\partial^{d_i}}{\partial z_{i+k}^{d_i}}
    \int_0^1 \prod_{j=1}^k
    \left(1-e(\theta)t^{\frac12 +z_j}\right)^{-1}
    \left(1-e(-\theta)t^{\frac12 -z_{k+j}}\right)^{-1}\,d\theta
    \biggr|_{z=0}.\notag \\
\end{eqnarray}
Now,
\begin{equation}
    \left(1-e(\theta)t^{\frac12 +z}\right)^{-1}
    =
    \sum_{m=0}^\infty (e(\theta)t^{\frac12+z})^m
\end{equation}
so 
\begin{equation}
    \label{eq:diff one factor}
    \frac{\partial^{c}}{\partial z^{c}} 
    \left(1-e(\theta)t^{\frac12 +z}\right)^{-1}
    \biggr|_{z=0}
    =
    \log(t)^c \sum_{m=0}^\infty (e(\theta)t^{\frac12})^m m^c.
\end{equation}
The sum above is of the form
\begin{equation}
     \sum_{m=0}^\infty w^m m^c
\end{equation}
which can be evaluate by applying $(w d/dw)^c$ to the geometric
series $1/(1-w) = \sum_{m=0}^\infty w^m$. This can be expressed either
in terms of Stirling numbers of the second kind or, alternatively in terms of 
Eulerian numbers~\cite{St}:
\begin{equation}
     \sum_{m=0}^\infty w^m m^c
     = \sum_{l=0}^c l! S(l,i) w^l (1-w)^{-l-1}
    = (1-w)^{-c-1} \sum_{l=0}^{c-1} E(c,l) w^{l+1}
\end{equation}
(the latter sum is taken to be 1 if c=0). We prefer to use the latter.
Thus, (\ref{eq:diff one factor}) equals
\begin{equation}
    \log(t)^c (1-e(\theta)t^{\frac12})^{-c-1}\sum_{l=0}^{c-1} E(c,l) e((l+1)\theta)t^{(l+1)/2}.
\end{equation}
Likewise,
\begin{eqnarray}
    &&\frac{\partial^{d}}{\partial z^{d}} 
    \left(1-e(-\theta)t^{\frac12-z}\right)^{-1}
    \biggr|_{z=0}
    = \notag \\
    &&(-\log(t))^d 
    (1-e(-\theta)t^{\frac12})^{-d-1}\sum_{l=0}^{d-1} E(d,l) e(-(l+1)\theta)t^{(l+1)/2}.
\end{eqnarray}

Applying this to (\ref{eq:partial diff int}) and expanding out, we need
to evaluate integrals of the form
\begin{equation}
    \int_0^1
    \left(1-e(\theta)t^{1/2}\right)^{-k-\sum c_i}
    \left(1-e(-\theta)t^{1/2}\right)^{-k-\sum d_i}
    e(C\theta)
    \,d\theta
\end{equation}
where $C \in {\mathbb Z}$, $-\sum d_i \leq  C \leq \sum c_i$.
The integral above can be expressed in terms of Gauss' hypergeometric series.
\begin{lem}
\label{lemma:hyper}
Let $A,B,C \in {\mathbb Z}$, $A,B \geq 1$, $0 \leq t < 1$. 
If $C \geq 0$ then
\begin{eqnarray}
    \int_0^1
    \left(1-e(\theta)t^{1/2}\right)^{-A}
    \left(1-e(-\theta)t^{1/2}\right)^{-B}
    e(C\theta)
    \,d\theta \notag \\
    = t^{C/2} { B+C-1 \choose C} {}_2F_1(A,B+C; C+1; t).
\end{eqnarray}
If $C<0$, then 
\begin{eqnarray}
    \int_0^1
    \left(1-e(\theta)t^{1/2}\right)^{-A}
    \left(1-e(-\theta)t^{1/2}\right)^{-B}
    e(C\theta)
    \,d\theta \notag \\
    = t^{|C|/2} { A+|C|-1 \choose |C|} {}_2F_1(B,A+|C|; |C|+1; t).
\end{eqnarray}
\end{lem}
\begin{pf}
Assume $C\geq 0$.
We can expand $\left(1-e(\theta)t^{1/2}\right)^{-A}$ and $\left(1-e(-\theta)t^{1/2}\right)^{-B}$
using the binomial series:
\begin{eqnarray}
     \left(1-e(\theta)t^{1/2}\right)^{-A}
     &=& 1 + A e(\theta) t^{1/2}
     + \frac{A(A+1)}{2!} (e(\theta) t^{1/2})^2
     + \frac{A(A+1)(A+2)}{3!} (e(\theta) t^{1/2})^3 + \ldots \notag \\
     \left(1-e(-\theta)t^{1/2}\right)^{-B}
     &=& 1 + B e(-\theta) t^{1/2}
     + \frac{B(B+1)}{2!} (e(-\theta) t^{1/2})^2
     + \frac{B(B+1)(B+2)}{3!} (e(-\theta) t^{1/2})^3 + \ldots \notag
\end{eqnarray}
Multiply these series together. The integral will pull out the coefficient of
$e(-C \theta)$, which equals
\begin{eqnarray}
    &&t^{C/2}
    \left(
        \frac{B\ldots(B+C-1)}{C!}+
        \frac{B\ldots(B+C)}{(C+1)!}At+
        \frac{B\ldots(B+C+1)}{(C+2)!}\frac{A(A+1)}{2!}t^2+ \ldots
    \right) \notag \\
    &&= t^{C/2} { B+C-1 \choose C} {}_2F_1(A,B+C; C+1; t).
\end{eqnarray}

The second formula in the lemma can be obtained by conjugating
the first and interchanging the role of $A$ and $B$.

\end{pf}

Using this lemma, we can write out the Taylor series
of~(\ref{eq:log second piece})
\begin{eqnarray}
    &&\log\left(
         \int_0^1 \prod_{j=1}^k
         \left(1-e(\theta)t^{\frac12 +z_j}\right)^{-1}
         \left(1-e(-\theta)t^{\frac12 -z_{k+j}}\right)^{-1}\,d\theta
    \right)
    = \notag \\
    &&\log({}_2F_1(k,k;1;t)) +
    \frac{t\log(t) k {}_2F_1(k+1,k+1;2;t)}{{}_2F_1(k,k;1;t)}
    \sum_1^k z_i-z_{i+k} + \notag \\
    &&+
    \left(
         \frac{-(t \log(t))^2 k^2{}_2F_1(k+1,k+1;2;t)^2}{{}_2F_1(k,k;1;t)^2} +
         \frac{(t \log(t))^2{k+1\choose 2}{}_2F_1(k+2,k+2;3;t)}{{}_2F_1(k,k;1;t)}
    \right) \notag \\
    &&\times \sum_{1\leq i<j\leq k} z_i z_j + z_{i+k}z_{j+k} \notag \\
    &&+
    \biggl(
         \frac{(t \log(t))^2 k^2{}_2F_1(k+1,k+1;2;t)^2}{{}_2F_1(k,k;1;t)^2} - \notag \\
         &&\frac{t\log(t)^2{}_2F_1(k+1,k+1;1;t)}{{}_2F_1(k,k;1;t)}
    \biggr)
    \times \sum_{1\leq i,j\leq k} z_i z_{j+k} \notag \\
    &&+
    \biggl(
         \frac{-(t \log(t))^2 k^2{}_2F_1(k+1,k+1;2;t)^2}{{}_2F_1(k,k;1;t)^2} +
         \frac{(t \log(t))^2{k+1\choose 2}{}_2F_1(k+2,k+2;3;t)}{{}_2F_1(k,k;1;t)} \notag \\
         &&\frac{t\log(t)^2 k{}_2F_1(k+2,k+1;2;t)}{{}_2F_1(k,k;1;t)}
    \biggr)
    \sum_1^{2k} \frac{z_i^2}{2} + \ldots
\end{eqnarray}

Combining the above with (\ref{eq:series log first piece}) we have that the first few coefficients
$B_k$ in (\ref{eq:A_k exp series}) are given by
\begin{eqnarray}
    \label{eq:B_k}
    B_k(1;) &=& \sum_p
    \frac{k\log(p)}{p-1} -
    \frac{\log(p) k\,{}_2F_1(k+1,k+1;2;1/p)}{p\,{}_2F_1(k,k;1;1/p)}
    \notag \\
    B_k(1,1;) &=& -\sum_p
    \left(
         \frac{\log(p)^2 k^2{}_2F_1(k+1,k+1;2;1/p)^2}{p^2\,{}_2F_1(k,k;1;1/p)^2} -
         \frac{\log(p)^2{k+1\choose 2}{}_2F_1(k+2,k+2;3;1/p)}{p^2\,{}_2F_1(k,k;1;1/p)}
    \right)
    \notag \\
    B_k(1;1) &=& \sum_p
    \frac{p\log(p)^2}{(p-1)^2} +
    \biggl(
         \frac{\log(p)^2 k^2{}_2F_1(k+1,k+1;2;1/p)^2}{p^2 {}_2F_1(k,k;1;1/p)^2} - \notag \\
         &&\frac{\log(p)^2{}_2F_1(k+1,k+1;1;1/p)}{p\,{}_2F_1(k,k;1;1/p)}
    \biggr)
    \notag \\
    B_k(2;) &=& -\sum_p
    \frac{kp\log(p)^2}{(p-1)^2} +
    \biggl(
         \frac{\log(p)^2 k^2{}_2F_1(k+1,k+1;2;1/p)^2}{p^2\,{}_2F_1(k,k;1;1/p)^2} - \notag \\
         &&\frac{\log(p)^2{k+1\choose 2}\,{}_2F_1(k+2,k+2;3;1/p)}{p^2\,{}_2F_1(k,k;1;1/p)} -
         \frac{\log(p)^2 k\,{}_2F_1(k+2,k+1;2;1/p)}{p\,{}_2F_1(k,k;1;1/p)}
    \biggr)
    \notag \\
\end{eqnarray}
% I double checked that these formulas give the same numerical results of our original computation
% for k=3 and 2. 

\subsubsection{Series for the $\zeta$ product}
\label{section:series for zeta product}
Let
\begin{equation}
\label{eq:zeta series}
    s \ \zeta(1+s) = 1 + \gamma_0 s - \gamma_1 s^2 +
    \frac{\gamma_2}{2!} s^3 -\frac{\gamma_3}{3!} s^4 + \cdots
\end{equation}
be the Laurent expansion of $s \zeta(1+s)$ about $s=0$, where the
$\gamma_n$'s generalize Euler's constant
\begin{equation}
\label{eq:gamma_n}
    \gamma_n = \lim_{m \to \infty}
    \sum_{k=1}^m \frac{\log(k)^n}{k}
    - \frac{\log(m)^{n+1}}{n+1}.
\end{equation}
As with the series for
$A_k$, we can here exploit the symmetries satisfied by the product
\begin{equation}
\label{eq:zeta product}
    \prod_{1 \leq i,j \leq k}
    (z_i-z_{j+k}) \zeta(1+z_i-z_{j+k}).
\end{equation}
We first set $z_{r+1}=\ldots=z_{k}=0$, $z_{k+r+1}=\ldots=z_{2k}=0$
before applying~(\ref{eq:partial diff operator}), so
that~(\ref{eq:zeta product}) is equivalent, in its series
expansion up to terms involving just $z_1,\ldots,z_r$,  to
\begin{equation}
    \prod_{1 \leq i,j \leq r}
    (z_i-z_{j+k}) \zeta(1+z_i-z_{j+k})
    \prod_{i=1}^r (z_i \zeta(1+z_i))^{k-r}
    (-z_{i+k} \zeta(1-z_{i+k}))^{k-r}.
\end{equation}
Again, one may use a symbolic differentiation package to
evaluate~(\ref{eq:partial diff operator}) applied to the above,
and thus obtain the coefficients of the multivariate series
expansion
\begin{eqnarray}
\label{eq:zeta product series}
    &&\prod_{1 \leq i,j \leq k}
    (z_i-z_{j+k}) \zeta(1+z_i-z_{j+k})
    =
    1+ \gamma k \sum_1^k z_i-z_{i+k} \notag \\
    &&+\gamma^2 k^2 \left( \sum_{1\leq i<j\leq k} z_i z_j + z_{i+k}z_{j+k} \right)
    +(2\gamma_1 +\gamma^2 - \gamma^2 k^2) \sum_{1\leq i,j\leq k} z_i z_{j+k} \notag \\
    &&+(\gamma^2 k^2 -\gamma^2 k - 2\gamma_1 k) \sum_1^{2k} \frac{z_i^2}{2} + \ldots
\end{eqnarray}

\subsubsection{Multiplying series together}
\label{section:multiplying} Let us be given two multivariate series of the
form that appears, for example, in~(\ref{eq:A_k exp series})
\begin{equation}
\label{eq:multivariate series}
    \sum_{\alpha;\beta}
    C_k(\alpha;\beta)
    \left(
    z_1^{\alpha_1} \cdots z_k^{\alpha_k}
    z_{k+1}^{\beta_{1}} \cdots z_{2k}^{\beta_{k}} \pm \text{sym}
    \right).
\end{equation}
In (\ref{eq:A_k exp series}) we pulled out the constant term from the series,
but generally, the above can have a non-zero constant term $C_k(;)$.
The `sym' indicates that we group together terms with 
exponents of the same form as explained in the introduction.

%For example, all the terms of degree $4$
%are collected with coefficients: $C_k(1,1,1,1;)$, $C_k(2,1,1;)$, $C_k(2,2;)$,
%$C_k(3,1;)$, $C_k(4;)$, $C_k(1,1,1;1)$, $C_k(2,1;1)$, $C_k(3;1)$, $C_k(1,1;1,1)$, $C_k(2;1,1)$,
%$C_k(2;2)$. The terms that go with $C_k(1,1,1,1;)$ are
%$\sum_{1\leq i_1 < i_2 <i_3 < i_4 \leq k}
%z_{i_1} z_{i_2} z_{i_3} z_{i_4} +z_{k+i_1} z_{k+i_2} z_{k+i_3} z_{k+i_4}$.
%Thus the sum over $\alpha;\beta$ follows the convention that $\alpha \geq \beta$ lexicographically,
%we list the $\alpha_i$'s and $\beta_i$'s in decreasing order, and we supress the $\alpha_i$'s
%and $\beta_i$'s that are 0.

We can easily obtain the coefficients
of their product by examining, for a given $\alpha;\beta$ the
various pairs $\gamma_1;\lambda_1$, $\gamma_2;\lambda_2$ with
$\gamma_1+\gamma_2=\alpha$, and $\lambda_1+\lambda_2=\beta$. Some
care is needed in carrying this out. While in collecting terms 
by $C_k(\alpha;\beta)$ we use the conventions in the above paragraph,
$\gamma$ and $\lambda$ need not
satisfy $\gamma \geq \lambda$. For example, say with $k=3$, a
term of the form $z_1^3 z_2 z_5^2$ can arise through
multiplication in 24 ways as $(z_1^{\gamma_{1,1}}
z_2^{\gamma_{1,2}} z_3^{\gamma_{1,3}} z_4^{\gamma_{2,1}}
z_5^{\gamma_{2,2}} z_6^{\gamma_{2,3}}) (z_1^{\lambda_{1,1}}
z_2^{\lambda_{1,2}} z_3^{\lambda_{1,3}} z_4^{\lambda_{2,1}}
z_5^{\lambda_{2,2}} z_6^{\lambda_{2,3}}))$ with
$\gamma_{1,1}+\lambda_{1,1}=3$,$\gamma_{1,2}+\lambda_{1,2}=1$,$\gamma_{2,2}+\lambda_{2,2}=2$,
and all the others equal to zero. For each of these 24 ways, one
needs to look up the corresponding coefficients of both series by
sorting the $\gamma$ and $\lambda$ and possibly swapping,
using (\ref{eq:B symmetry}), so that $\gamma \geq \lambda$.

In this manner, we are able, given the series for
$A_k(z_1,\ldots,z_{2k})$ and the series in~(\ref{eq:zeta product
series}), to obtain the series for the second line
in~(\ref{eq:huge integral})
\begin{eqnarray}
\label{eq:series for integrand}
    &&A_k(z_1,\ldots,z_{2k})
    \prod_{1 \leq i,j \leq k}
    (z_i-z_{j+k}) \zeta(1+z_i-z_{j+k})
    \notag \\
    && = a_k \sum_{\alpha;\beta}
    b_k(\alpha;\beta)
    \left(
    z_1^{\alpha_1} \cdots z_k^{\alpha_k}
    z_{k+1}^{\beta_1} \cdots z_{2k}^{\beta_k} \pm \text{sym}
    \right)
\end{eqnarray}
(to determine the multivariate series for $A_k$ from that of
$\log(A_k(z_1,\ldots,z_{2k}))$, one uses the Taylor series for
$\exp$, applied to $\log(A_k)$, and the above multiplication
algorithm).

We list the first few terms:
\begin{eqnarray}
    \label{eq: b_k series}
    &&A_k(z_1,\ldots,z_{2k})
    \prod_{1 \leq i,j \leq k}
    (z_i-z_{j+k}) \zeta(1+z_i-z_{j+k})
    \notag \\
    =
    &&a_k \biggl(
        1+ (\gamma k + B_k(1;)) \sum_1^k z_i-z_{i+k} \notag \\
        &&+(\gamma^2 k^2 + B_k(1,1;) + B_k(1;)^2 + 2 \gamma k B_k(1;))
         \sum_{1\leq i<j\leq k} z_i z_j + z_{i+k}z_{j+k} \notag \\
        &&+(2\gamma_1 +\gamma^2 - \gamma^2 k^2  +B_k(1;1) -B_k(1;)^2 -2\gamma k B_k(1;)) \sum_{1\leq i,j\leq k} z_i z_{j+k} \notag \\
        &&+\frac{1}{2}(\gamma^2 k^2 -\gamma^2 k - 2\gamma_1 k +B_k(2;) +B_k(1;)^2 + 2\gamma k B_k(1;)  ) \sum_1^{2k} z_i^2 + \ldots
    \biggr)\notag \\
\end{eqnarray}

\subsection{Determining $c_r(k)$}
\label{subsubsection: c_r}

Extracting the terms of degree $r$ from the Taylor expansion
(\ref{eq:series for integrand}) and substituting into~(\ref{eq:huge integral}) we have
\begin{eqnarray}
\label{eq:c_r integral}
    c_r(k)x^{k^2-r} = &&\frac{a_k}{k!^2}\frac{1}{(2\pi i)^{2k}}
    \oint\cdots \oint
    \frac{
        \Delta(z_1,\ldots,z_{2k})
        \Delta(z_1,\ldots,z_{k})
        \Delta(z_{k+1},\ldots,z_{2k})
     }
     {
        \displaystyle \prod_{i=1}^{2k} z_i^{2k}
     }
    \notag
    \\
    \times
    &&\sum_{|\alpha|+|\beta|=r}
    b_k(\alpha;\beta)
    \left( z_1^{\alpha_1} \cdots z_k^{\alpha_k} 
           z_{k+1}^{\beta_1} \cdots z_{2k}^{\beta_k} 
           \pm \text{sym} 
    \right) \notag \\
    \times
    &&e^{\tfrac x2\sum_{i=1}^{k}z_i-z_{i+k}}
    dz_1\ldots dz_{2k},
\end{eqnarray}
Now, each term in the second line of the integrand with the exponents of the same form
integrates the same since the integrand is a symmetric function of
$z_1,\ldots, z_k$ and of $z_{k+1},\ldots, z_{2k}$, and also because the contribution
from $z_1^{\alpha_1} \cdots z_k^{\alpha_k}  z_{k+1}^{\beta_1} \cdots z_{2k}^{\beta_k}$
is the same as for $ z_{1}^{\beta_1} \cdots z_{k}^{\beta_k} z_{k+1}^{\alpha_1} \cdots z_{2k}^{\alpha_k}$,
as can be seen by changing variables $u_i=-z_i$ and using
$b_k(\alpha;\beta) = (-1)^{|\alpha|+|\beta|} b_k(\beta;\alpha)$.

Let $d_k(\alpha;\beta)$ denote the number of terms of a given form of exponent.
For example, $d_k(1;)=2k$ since there are $2k$ terms in $\sum_{i=1}^k z_i-z_{i+k}$.
Once can write down a formula for $d_k(\alpha;\beta)$ as, up to a factor of 2,
a multinomial coefficient in terms
of the multiplicities of the values assumed by the $\alpha_i$'s and $\beta_i$'s.
Let $m_\alpha(j)$ denote the number of occurences of $j$ amongst $\alpha_1,\ldots,\alpha_k$,
and likewise for $\beta$.

%XXXXX maybe get rid of J_1 and J_2 to keep things cleaner.
Let $J_1$ denote the largest value amongst the $\alpha_i$'s and $J_2$ the largest value
amongst the $\beta_i$'s. Let
\begin{equation}
   \label{eq:delta alpha beta}
   \delta(\alpha;\beta) =
   \begin{cases}
       1 \quad \text{if $\alpha = \beta$} \\
       0 \quad \text{if $\alpha \neq \beta$}
   \end{cases}.
\end{equation}
We have introduced $\delta(\alpha;\beta)$ to
take into account that the $\pm \text{sym}$ in (\ref{eq:c_r integral}) 
collects together the terms corresponding to $\alpha;\beta$ and to $\beta;\alpha$.
We have
\begin{eqnarray}
    \label{eq:d_k}
    &&d_k(\alpha;\beta) = \notag \\
    &&2^{1-\delta(\alpha;\beta)}
    {k \choose m_\alpha(0)}
    {k -m_\alpha(0) \choose m_\alpha(1)}
    {k -m_\alpha(0) - m_\alpha(1) \choose m_\alpha(2)}
    \cdots
    {k -m_\alpha(0) \ldots - m_\alpha(J_1-1)  \choose m_\alpha(J_1)} \notag \\
    &&\times
    {k \choose m_\beta(0)}
    {k -m_\beta(0) \choose m_\beta(1)}
    {k -m_\beta(0) - m_\beta(1) \choose m_\beta(2)}
    \cdots
    {k -m_\beta(0) \ldots - m_\beta(J_2-1)  \choose m_\beta(J_2)} \notag \\
    &&=
    2^{1-\delta(\alpha;\beta)}
    (k!)^2
    \prod_{j=0}^{J_1} \frac{1}{m_\alpha(j)!}
    \prod_{j=0}^{J_2} \frac{1}{m_\beta(j)!},
\end{eqnarray}
since there are ${k \choose m_\alpha(0)}$ ways to choose which $z_i$'s, $1 \leq i \leq k$, 
have exponent $0$, then ${k -m_\alpha(0) \choose m_\alpha(1)}$ ways to decide which of the remaining
$z_i$'s have exponent $1$, etc, and likewise for $\beta$. In the simplification to obtain the
second line we used $\sum_{j=0}^{J_1} m_\alpha(j) = k$, and similarly for $\beta$.

Therefore, counting the number of terms that are collected for a given
$b_k(\alpha;\beta)$ we get
\begin{eqnarray}
\label{eq:c_r integral 2}
    &&c_r(k)x^{k^2-r} = \frac{a_k}{k!^2}\frac{1}{(2\pi i)^{2k}}
    \sum_{|\alpha|+ |\beta| =r}
    b_k(\alpha;\beta) d_k(\alpha;\beta)
    \notag \\
    &&\times
    \oint\cdots \oint
    \frac{
        \Delta(z_1,\ldots,z_{2k})
        \Delta(z_1,\ldots,z_{k})
        \Delta(z_{k+1},\ldots,z_{2k})
     }
     {
        \displaystyle \prod_{i=1}^{2k} z_i^{2k}
     }
    %\notag
    %\\
    %\times
    z_1^{\alpha_1} \cdots z_k^{\alpha_k} z_{k+1}^{\beta_1} \cdots z_{2k}^{\beta_k} \notag \\
    &&\times
    e^{\tfrac x2\sum_{i=1}^{k}z_i-z_{i+k}}
    dz_1\ldots dz_{2k},
\end{eqnarray}

Pulling out
\begin{eqnarray}
\label{eq:vands and series}
    \Delta(z_1,\ldots,z_{2k})
    \Delta(z_1,\ldots,z_{k})
    \Delta(z_{k+1},\ldots,z_{2k})
    \sum_{|\alpha|+ |\beta| =r}
    b_k(\alpha;\beta) d_k(\alpha;\beta)
    z_1^{\alpha_1} \cdots z_k^{\alpha_k}
    z_{k+1}^{\beta_1} \cdots z_{2k}^{\beta_k} \notag \\
\end{eqnarray}
from the integral, we have, as in our consideration of $c_0(k)$,
\begin{eqnarray}
\label{eq:c_r a}
    c_r(k) =
    \frac{a_k}{2^{k^2-r} k!^2}
    q_2(\partial/\partial x_1,\ldots,\partial/\partial x_{2k})
    q(\partial/\partial x_1,\ldots,\partial/\partial x_{2k})
    \prod_1^{2k}\frac{x_i^{2k-1}}{(2k-1)!}
    \notag \\
\end{eqnarray}
evaluated at $x_1 = \ldots = x_k=1, x_{k+1} = \ldots = x_{2k}=-1$,
where $q$ is given by (\ref{eq:triple vandermonde}) and
\begin{equation}
\label{eq:q2}
    q_2(z_1,\ldots,z_{2k}) =
    \sum_{|\alpha|+ |\beta| =r}
    b_k(\alpha;\beta) d_k(\alpha;\beta)
    z_1^{\alpha_1} \cdots z_k^{\alpha_k}
    z_{k+1}^{\beta_1} \cdots z_{2k}^{\beta_k} .
\end{equation}
As in the proof of Lemma \ref{lemma:triple vand}, (\ref{eq:c_r a}) equals
\begin{equation}
\label{eq:c_r b}
    \frac{a_k}{2^{k^2-r} k!^2}
    q_2(\partial/\partial x_1,\ldots,\partial/\partial x_{2k})
    \sum_{\mu,\nu}
    g(\mu,\nu)
\end{equation}
where $g(\mu,\nu)$ equals (\ref{eq:typical term}) (with
$f(x)=x^{2k-1}/(2k-1)!$)  and is comprised of a sign and a
determinant. The sum is over all $k!^2$ pairs of permutations
$\mu,\nu$.

Applying $q_2(\partial/\partial
x_1,\ldots,\partial/\partial x_{2k})$ to these determinants we get
\begin{equation}
\label{eq:c_r c}
    \frac{a_k}{2^{k^2-r}}
    \sum_{|\alpha|+ |\beta| =r}
    \sum_{\mu,\nu}
    \frac{b_k(\alpha,\beta) d_k(\alpha;\beta)}{k!^2}
    \text{sgn}(\mu) \ \text{sgn}(\nu)
    M_k(\mu,\nu,\alpha,\beta)
\end{equation}
with
\begin{eqnarray}
\label{eq:M}
    &&M_k(\mu,\nu,\alpha,\beta) = \notag \\
    &&\left|
             \begin{array}{cccc}
             f^{(\mu_1-1+\alpha_1)}(x_1) & f^{(\mu_1+\alpha_1)}(x_1) & \ldots & f^{(\mu_1+2k-2+\alpha_1)}(x_1) \\
             f^{(\mu_2-1+\alpha_2)}(x_2) & f^{(\mu_2+\alpha_2)}(x_2) & \ldots & f^{(\mu_2+2k-2+\alpha_2)}(x_2) \\
             \vdots & \vdots & \ddots & \vdots \\
             f^{(\mu_k-1+\alpha_k)}(x_k) & f^{(\mu_k+\alpha_k)}(x_k) & \ldots & f^{(\mu_k+2k-2+\alpha_k)}(x_k) \\
             f^{(\nu_1-1+\beta_1)}(x_{k+1}) & f^{(\nu_1+\beta_1)}(x_{k+1}) & \ldots & f^{(\nu_1+2k-2+\beta_1)}(x_{k+1}) \\
             f^{(\nu_2-1+\beta_2)}(x_{k+2}) & f^{(\nu_2+\beta_2)}(x_{k+2}) & \ldots & f^{(\nu_2+2k-2+\beta_2)}(x_{k+2}) \\
             \vdots & \vdots & \ddots & \vdots \\
             f^{(\nu_k-1+\beta_k)}(x_{2k}) & f^{(\nu_k+\beta_k)}(x_{2k}) & \ldots & f^{(\nu_k+2k-2+\beta_k)}(x_{2k})
             \end{array}
    \right|_{2k \times 2k.}
\end{eqnarray}
Setting $x_1=\ldots=x_k=1$, $x_{k+1}=\ldots=x_{2k}=-1$, and
rearranging rows (to undo the $\mu$ and $\nu$) we get
\begin{equation}
\label{eq:c_r d}
    c_r(k) = \frac{a_k}{2^{k^2-r}}
    \sum_{|\alpha|+ |\beta| =r}
    \sum_{\sigma,\tau}
    \frac{b_k(\alpha;\beta) d_k(\alpha;\beta)}{k!^2}
    \tilde{M}_k(\sigma(\alpha),\tau(\beta))
\end{equation}
with
\begin{eqnarray}
\label{eq:M tilde}
        &&\tilde{M}_k(\sigma(\alpha),\tau(\beta)) = \notag \\
           &&(-1)^{\sum \beta_i}
           \left|
             \begin{array}{cccc}
             \Gamma(2k-\alpha_{\sigma_1})^{-1} &  \Gamma(2k-1-\alpha_{\sigma_1})^{-1}& \ldots & \Gamma(1-\alpha_{\sigma_1})^{-1} \\
             \Gamma(2k-1-\alpha_{\sigma_2})^{-1} & \Gamma(2k-2-\alpha_{\sigma_2})^{-1} & \ldots & \Gamma(-\alpha_{\sigma_2})^{-1} \\
             \vdots & \vdots & \ddots & \vdots \\
             \Gamma(k+1-\alpha_{\sigma_k})^{-1} & \Gamma(k-\alpha_{\sigma_k})^{-1} & \ldots & \Gamma(2-k-\alpha_{\sigma_k})^{-1} \\
             -\Gamma(2k-\beta_{\tau_1})^{-1} & \Gamma(2k-1-\beta_{\tau_1})^{-1} & \ldots &  \Gamma(1-\beta_{\tau_1})^{-1}\\
             \Gamma(2k-1-\beta_{\tau_2})^{-1} & -\Gamma(2k-2-\beta_{\tau_2})^{-1} & \ldots & -\Gamma(-\beta_{\tau_2})^{-1} \\
             \vdots & \vdots & \ddots & \vdots \\
             (-1)^k\Gamma(k+1-\beta_{\tau_k})^{-1} & (-1)^{k+1}\Gamma(k-\beta_{\tau_k})^{-1} & \ldots & (-1)^{3k-1}\Gamma(2-k-\beta_{\tau_k})^{-1}
             \end{array}
         \right|_{2k \times 2k.} \notag \\
\end{eqnarray}
The extra factor of $(-1)^{\sum \beta_i}$ comes from the extra powers of $-1$ that are pulled out of
the bottom $k$ rows of the matrix. Notice that in order for $\tilde{M}_k$ to be non-zero,
the $\alpha_{\sigma_i}+i-1$'s must be a distinct subset of $\set{0,1,2,...,2k-1}$,
and similarly for the $\beta_{\tau_i}+i-1$'s. The implication of this latter point is
discussed further in Lemma~\ref{lemma:det 0} below.

For any $\alpha,\beta$, many of the $k!^2$ pairs of permutations
$\sigma,\tau$ will give the same determinant because of
multiplicity amongst the $\alpha_i$'s and $\beta_i$'s. In fact,
since $r$ is fixed, most of the $\alpha_i$'s and $\beta_i$'s will
equal zero. As before, let $m_\alpha(j)$ denote the number of
occurrences of $j$ in $\alpha_1,\ldots,\alpha_k$, and similarly for
$\beta$.

Then~(\ref{eq:c_r d}) becomes
\begin{equation}
\label{eq:c_r e}
    c_r(k) = \frac{a_k}{2^{k^2-r}}
    \sum_{|\alpha|+ |\beta| =r}
    \frac{b_k(\alpha;\beta) d_k(\alpha;\beta)}{k!^2}
    \prod_{j=0}^{J_1} m_\alpha(j)!
    \prod_{j=0}^{J_2} m_\beta(j)!
    \sum_{{\text{rearrangements}} \atop {\text{$\sigma,\tau$ of $\alpha$ and $\beta$}}}
    \tilde{M}_k(\sigma(\alpha),\tau(\beta)).
\end{equation}
By rearrangements, we mean distinct permutations as explained following~(\ref{eq:N thm}).
Recall that $J_1$ and $J_2$ denote the largest
value amongst the $\alpha_i$'s and $\beta_i$'s respectively.

Notice that this simplifies since, by (\ref{eq:d_k}),
\begin{equation}
    \frac{d_k(\alpha;\beta)}{k!^2}
    \prod_{j=0}^{J_1} m_\alpha(j)!
    \prod_{j=0}^{J_2} m_\beta(j)! = 2^{1-\delta(\alpha;\beta)}
\end{equation}
is constant, where $\delta(\alpha;\beta)$ is given by (\ref{eq:delta alpha beta}).
%yahoo, so simple!

Therefore
\begin{equation}
\label{eq:c_r f}
    c_r(k) = \frac{a_k}{2^{k^2-r}}
    \sum_{|\alpha|+ |\beta| =r}
    2^{1-\delta(\alpha;\beta)}
    b_k(\alpha;\beta)
    \sum_{{\text{rearrangements}} \atop {\text{$\sigma,\tau$ of $\alpha$ and $\beta$}}}
    \tilde{M}_k(\sigma(\alpha),\tau(\beta)).
\end{equation}

Expression (\ref{eq:c_r e}) can be pared down further by realizing
that `all of the action' takes place in rows $k-|\alpha|+1,\ldots,k$ and
$2k-|\beta|+1,\ldots,2k$. By this we mean that we need only focus on the
rearrangements that have $\alpha_{\sigma_1}= \ldots =
\alpha_{\sigma_{k-|\alpha|}} = 0$, and $\beta_{\sigma_1}= \ldots =
\beta_{\sigma_{k-|\beta|}} = 0$, because otherwise the determinant will
equal zero.

\begin{lem}
\label{lemma:det 0}
The determinant $\tilde{M}_k(\sigma(\alpha),\tau(\beta))$ in
(\ref{eq:M tilde}) equals zero, unless $\alpha_{\sigma_1} = \ldots
= \alpha_{\sigma_{k-|\alpha|}}= \beta_{\tau_1} = \ldots =
\beta_{\tau_{k-|\beta|}}=0$ (in which case it might, or might not, equal
zero).

\end{lem}
\begin{pf}
Assume that $\tilde{M}_k(\sigma(\alpha),\tau(\beta)) \neq 0$.
Let $\delta_1 := \alpha_{\sigma_i} \geq 1$ be equal to the first
non zero $\alpha_{\sigma}$. But this forces $\delta_2 :=
\alpha_{\sigma_{i+\delta_1}}$ to also be $\geq 1$, otherwise rows
$i$ and $i+\delta_1$ would coincide and the determinant would be
zero. But then $\delta_3 := \alpha_{\sigma_{i+\delta_1+\delta_2}}$
must also be $\geq 1$ otherwise rows $i+\delta_1$ and
$i+\delta_1+\delta_2$ would coincide. Continue in this fashion
until reaching beyond the $k$th row, $i+\delta_1+\delta_2+ \cdots
+ \delta_j > k$. But, $\delta_1+ \cdots \delta_j \leq \sum_{m=1}^k
\alpha_m = |\alpha|$, so $i+|\alpha|>k$ i.e. $i>k-|\alpha|$. 
Thus, $\alpha_{\sigma_1} =
\ldots = \alpha_{\sigma_{k-|\alpha|}}=0$. Similarly, $\beta_{\tau_1} =
\ldots = \beta_{\tau_{k-|\beta|}}=0$, the only difference in the proof
being that the rows would coincide up to a factor of $\pm 1$.

\end{pf}

The above lemma greatly improves the speed with which we can
evaluate (\ref{eq:c_r e}) since all but $O_r(1)$ of the terms can be
discarded. Consider
\begin{equation}
    \frac{1}{2^{k^2-r}}
    \sum_{{\text{rearrangements}} \atop {\text{$\sigma,\tau$ of $\alpha$ and $\beta$}}}
    \tilde{M}_k(\sigma(\alpha),\tau(\beta)).
\end{equation}
We will prove in Section~\ref{section:poly} using the theory of factorial
Schur functions that the above is equal to $\prod_{l=0}^{k-1} l!/(k+l)!$ times
a polynomial in $k$ of degree $\leq 2(|\alpha| + |\beta|)$,
or else is the $0$ polynomial. 

Hence, if we let $N_k(\alpha;\beta)$ denote the polynomial
%warning- this definition of N_k is slightly different than what we used
%when we intially worked out the lower terms. The difference is that previously
%we kept d_k(alpha;beta) as a polynomial in k, and moved the 1/k!^2 times the
%m_alpha(j)! m_beta(j)! 's into N_k. The current approach is better since
%we have managed to eliminate the d_k. The fact that both cases presumably lead to
%polynomials says something like: N_k should be divisible by 
% k*(k-1)...(m_alpha(0)+1) * k*(k-1)...(m_beta(0)+1)
\begin{equation}
    \label{eq:N}
    N_k(\alpha;\beta) = 
    \frac{1}{2^{k^2-r}}
    \left( \prod_{l=0}^{k-1} \frac{l!}{(k+l)!} \right)^{-1}
    \sum_{{\text{rearrangements}} \atop {\text{$\sigma,\tau$ of $\alpha$ and $\beta$}}}
    \tilde{M}_k(\sigma(\alpha),\tau(\beta))
\end{equation}
we have thus arrived at:
\begin{equation}
    \label{eq:c_r g}
    %\boxed{
    c_r(k) = a_k \prod_{l=0}^{k-1} \frac{l!}{(k+l)!}
    \sum_{|\alpha|+ |\beta| =r}
    2^{1-\delta(\alpha;\beta)}
    b_k(\alpha;\beta) N_k(\alpha;\beta).
    %}
\end{equation}
The factor $a_k \prod_{l=0}^{k-1} l!/(k+l)!$ is equal to the leading coefficient $c_0(k)$.
The $b_k(\alpha,\beta)$'s are the Taylor coefficients of
$\frac{1}{a_k} A_k(z_1,\ldots,z_{2k}) \prod_{1 \leq i,j \leq k} (z_i-z_{j+k}) \zeta(1+z_i-z_{j+k})$,
and the first few are listed in (\ref{eq: b_k series}). The function $\delta(\alpha;\beta)$
equals zero unless $\alpha=\beta$ in which case it equals 1.
We have thus managed to express
$c_r(k)$ as equal to $c_0(k)$ times a polynomial in $k$ with
coefficients linear in the $b_k(\alpha,\beta)$'s.

Knowing that $N_k(\alpha;\beta)$ is a polynomial of degree $\leq 2(|\alpha| + |\beta|)$
allows us to determine it for a given $\alpha;\beta$ by evaluating~(\ref{eq:N})
at $2(|\alpha| + |\beta|)+1$ different values of $k$ and writing the unique polynomial
of degree $\leq 2(|\alpha| + |\beta|)$ that interpolates those values. Since the arithmetic
just involves rational numbers it can be performed exactly. When evaluating the r.h.s.
of~(\ref{eq:N}) one should make sure to exploit Lemma~\ref{lemma:det 0} so as to 
only evaluate $O_r(1)$ of the rearrangements.

In this way, one can find, for example,
\begin{eqnarray}
    \label{eq:N exs}
    N_k(1;) &=& k^2 \notag \\ 
    N_k(2;) &=& 0 \notag \\  %there are two rearrangements that contribute: alpha_{sigma_k}=2 and
                    %alpha_{sigma_{k-1}}=2. the first has poly equal to k^4/2-k^2/2 and the other
                    %has poly equal to -(k^4/2-k^2/2). Together they give 0.
    N_k(1,1;) &=& k^2(k-1)(k+1)/2 \notag \\
    N_k(1;1) &=& -k^2(k-1)(k+1).
\end{eqnarray}
This allows us to write down formulae for $c_1(k)$ and $c_2(k)$:
\begin{eqnarray}
    \label{eq:c1}
    c_1(k) &=&  \left(a_k \prod_{l=0}^{k-1} \frac{l!}{(k+l)!} \right)  2 k^2 b_k(1;) \notag \\
           &=&  \left(a_k \prod_{l=0}^{k-1} \frac{l!}{(k+l)!} \right)  2 k^2 (\gamma k + B_k(1;))
\end{eqnarray}
and, after simplifying,
\begin{eqnarray}
    \label{eq:c2}
    c_2(k) &=& \left(a_k \prod_{l=0}^{k-1} \frac{l!}{(k+l)!} \right)
           k^2 (k-1)(k+1) (b_k(1,1;) - b_k(1;1)) \notag \\
           &=& \left(a_k \prod_{l=0}^{k-1} \frac{l!}{(k+l)!} \right)
           k^2 (k-1)(k+1) \notag \\
           &&\times (
               %(\gamma^2 k^2 + B_k(1,1;) + B_k(1;)^2 + 2 \gamma k B_k(1;)) -
               %(2\gamma_1 +\gamma^2 - \gamma^2 k^2  +B_k(1;1) -B_k(1;)^2 -2\gamma k B_k(1;))
               %note the extra factor of 1/2 in the first N_k is canceled by the extra power of
               %2^(1-delta).
               %simplified:
               2(B_k(1;)+\gamma k)^2 - \gamma^2 -2\gamma_1+B_k(1,1;)-B_k(1;1)
           ).
           \notag \\
\end{eqnarray}
The $B_k$'s above are given in (\ref{eq:B_k}). In practice,
we were thus able to explicitly determine the first nine lower order
terms, $r \leq 9$. 

\section{Proof that $N_k(\alpha;\beta)$ is a polynomial.}
\label{section:poly}

Throughout this section we use the following notation.
Let $e_1,e_2+1,\ldots,e_k+k-1$ be distinct integers and
$f_1,f_2+1,\ldots,f_k+k-1$ be distinct integers. If
$f_1,f_2+1,\ldots,f_k+k-1$ is a subset of $0,1,\ldots,2k-1$, let
$c_1,\ldots,c_k$ be the complementary subset. Later we will
introduce some extra assumptions on the $e_i$'s and $f_i$'s,
namely that most of them are equal to zero.

The following lemma expresses the kind of $2k \times 2k$
determinant that appears in the formula for $N_k$ as a
a $k \times k$ determinant involving binomial coefficients.

\begin{lem}
\label{lemma: poly1} 
Let $e_i$ and $f_i$ be given as above. Then
\begin{eqnarray}
\label{eq:M tilde 2}
           &\left|
             \begin{array}{cccc}
             \Gamma(2k-e_{1})^{-1} &  \Gamma(2k-1-e_{1})^{-1}& \ldots & \Gamma(1-e_{1})^{-1} \\
             \Gamma(2k-1-e_{2})^{-1} & \Gamma(2k-2-e_{2})^{-1} & \ldots & \Gamma(-e_{2})^{-1} \\
             \vdots & \vdots & \ddots & \vdots \\
             \Gamma(k+1-e_{k})^{-1} & \Gamma(k-e_{k})^{-1} & \ldots & \Gamma(2-k-e_{k})^{-1} \\
             -\Gamma(2k-f_{1})^{-1} & \Gamma(2k-1-f_{1})^{-1} & \ldots &  \Gamma(1-f_{1})^{-1}\\
             \Gamma(2k-1-f_{2})^{-1} & -\Gamma(2k-2-f_{2})^{-1} & \ldots & -\Gamma(-f_{2})^{-1} \\
             \vdots & \vdots & \ddots & \vdots \\
            (-1)^k\Gamma(k+1-f_{k})^{-1} & (-1)^{k+1}\Gamma(k-f_{k})^{-1} & \ldots & (-1)^{3k-1}\Gamma(2-k-f_{k})^{-1}
             \end{array}
         \right|_{2k \times 2k.} \notag \\
         =
         &
         \text{sgn}(f)
         \left(
             \prod_{l=0}^{k-1}
             \frac{(e_{l+1}+l)!(f_{l+1}+l)!}{l!(k+l)!} 2^{c_{l+1}-e_{l+1}-l}
         \right)
         \left|
               c_{j} \choose e_i+i-1
         \right|_{k \times k},
\end{eqnarray}
where $\text{sgn}(f)$ is $(-1)$ raised to the number of transpositions needed to get
$f_1+1, f_2+2, \ldots, f_k +k$ sorted into increasing numerical order. 

Also note, if $e_i$'s and $f_i$'s are not distinct as required, then two of the rows on
the l.h.s. will coincide up to sign, and hence the determinant will equal zero.
\end{lem}

\begin{pf}
  Introducing a $(2k-j)!$ in column $j$,
  $1/(e_i+i-1)!$ in row $i$ and $(-1)^{i-1}/(f_i+i-1)!$ in row $k+i$,
  $i=1,\ldots,k$,
  the l.h.s. of (\ref{eq:M tilde 2}) equals
\begin{equation}
     \label{eq:M step1}
    (-1)^{\frac{(k-1)k}{2}} \prod_{l=0}^{k-1}
         \frac{(e_{l+1}+l)!(f_{l+1}+l)!}{l!(k+l)!}
    \left|
             \begin{array}{cccc}
             2k-1 \choose e_{1} &  2k-2 \choose e_{1}& \ldots & 0 \choose e_{1} \\
             2k-1 \choose e_{2}+1 &  2k-2 \choose e_{2}+1 & \ldots & 0 \choose e_{2}+1 \\
             \vdots & \vdots & \ddots & \vdots \\
               2k-1 \choose e_{k}+k-1 &  2k-2 \choose e_{k}+k-1& \ldots & 0 \choose e_{k}+k-1 \\
             -  {2k-1 \choose f_{1}} & + {2k-2 \choose f_{1}} & \ldots & +  {0 \choose f_{1}}\\
             -  {2k-1 \choose f_{2}+1} & + {2k-2 \choose f_{2}+1} & \ldots & +  {0 \choose f_{2}+1}\\
             \vdots & \vdots & \ddots & \vdots \\
             -  {2k-1 \choose f_{k}+k-1} & + {2k-2 \choose f_{k}+k-1} & \ldots & +  {0 \choose
             f_{k}+k-1}
             \end{array}
         \right|_{2k \times 2k.} \ \ \ \ \
\end{equation}
   Multiplying the matrix in (\ref{eq:M step1}) on the right by
   a unit, i.e. by
   \begin{equation}
       \label{eq:right multiply}
       \left|
           j-1 \choose 2k-i
       \right|_{2k\times 2k}
       \ \ = \ \ \ \ (-1)^k
   \end{equation}
   the $ij$ entry is, for $1 \leq i \leq k$, $1 \leq j \leq 2k$
   \begin{equation}
       \label{eq: upper entries}
       \sum_{m=0}^{2k-1} {m \choose e_i+i-1}{j-1 \choose m} =
       2^{j-e_i-i} {j-1 \choose e_i+i-1},
   \end{equation}
   and  the $(i+k)j$ entry is, for $1 \leq i \leq k$, $1 \leq j \leq 2k$
   \begin{equation}
       \label{eq: lower entries}
       \sum_{m=0}^{2k-1} {m \choose f_i+i-1}{j-1 \choose m}(-1)^m =
       \begin{cases}
           (-1)^{f_i+i-1} &\text{if $j = f_i+i$} \\
           0 &\text{otherwise}.
       \end{cases}
   \end{equation}
   These binomial identities can be proven by noticing that ${m \choose A}{B \choose m}
   = {B \choose A}{B-A \choose m-A}$ and using the binomial
   theorem.

   Expanding the determinant of this new matrix along the last $k$
   rows, and pulling out powers of $2$'s from the first $k$ rows
   gives the lemma. The $\text{sgn}(f)$ that appears in the lemma
   can be obtained as follows. The last $k$ rows of the matrix given by
   (\ref{eq: upper entries}) and (\ref{eq: lower entries}) consist of $\pm 1$'s
   with a $(-1)^{f_i+i-1}$ appearing in row $k+i$ and column $f_i+i$, $1 \leq i \leq k$.
   We can swap the columns and rows of this matrix so that the bottom left $k\times k$
   submatrix becomes diagonal, and the lower right submatrix becomes, $0_{k \times k}$,
   the zero submatrix.

   There are many ways to do so, but to end up with the
   determinant in the lemma, one should make sure that we do not rearrange the relative
   ordering of the columns corresponding to $c_1,\ldots,c_k$. If the quantities
   $f_i+i$ appear in increasing numerical value, one can simply swap column
   $f_i+i$ with its neighbouring columns on the left, one at a time, until it
   sits in the $i$th column. This introduces a $(-1)^{\sum_{i=1}^k f_i}$
   into the determinant.

   However, if the $f_i+i$'s appear out of order, in order to preserve
   the ordering of the columns corresponding to $c_1,\ldots,c_k$, one should
   first swap rows so as to put the lower $k\times 2k$ submatrix into reduced row
   echelon form. For example if one has $f_{i_1}+i_1 > f_{i_2}+i_2$, but $i_1<i_2$ then one
   should swap rows $k+i_1$ and $k+i_2$. This has the effect of placing the $(-1)^{f_{i_1}+i_1-1}$
   in entry $(i_2, f_{i_1}+i_1) = (i_2, (f_{i_1}+i_1-i_2)+ i_2)$, and the
   $(-1)^{f_{i_2}+i_2-1}$ in entry  $(i_1, f_{i_2}+i_2) = (i_1, (f_{i_2}+i_2-i_1)+ i_1)$.
   The horizontal displacement then needed to get these entries into the $i_2$nd
   and $i_1$st columns is therefore unchanged and equal to $f_{i_1}+f_{i_2}$. Relabeling and repeating
   if necessary, one sees that the contribution to the determinant from the row and and column swaps that
   get the lower left $k\times k$ submatrix into diagonal form is
   \begin{equation}
      \label{eq:minus 1}
      (-1)^{\text{sgn}(f)+\sum f_i},
   \end{equation}
   where $\text{sgn}(f)$ accounts for the number of transpositions needed to get 
   $f_1+1,\ldots, f_k+k$ into increasing numerical order.

   One now easily evaluates the determinant by expanding along the lower diagonal matrix.
   This submatrix begins at entry $k+1,1$ and this contributes a $(-1)^k$ to the determinant.
   One also needs to multiply the diagonal entries themselves, and this contributes
   a $(-1)^{\sum f_i+i-1}$.

   Collecting the powers of $-1$ that appear in (\ref{eq:M step1}) and (\ref{eq:right multiply})
   and multiplying by (\ref{eq:minus 1}) and by the two factors in the previous paragraph we
   obtain the sign that appears in the lemma.
\end{pf}

Now (\ref{eq:M tilde 2}) equals
\begin{equation}
    \label{eq:factorial schur 1}
    (-1)^{k(k-1)/2}
    \text{sgn}(f)
    \prod_{l=0}^{k-1}
    \frac{(f_{l+1}+l)!}{l!(k+l)!} 2^{c_{l+1}-e_{l+1}-l}
    \left|
        (c_i)_{\lambda_j+k-j}
    \right|_{k \times k}
\end{equation}
with $(z)_\mu = z(z-1)\ldots(z-\mu+1)$ the {\it descending}
factorial, and $\lambda_j = e_{k-j+1}$.
The determinant above is essentially a factorial Schur function.

We have introduced the $\lambda_j$'s and taken the tranpose
so as to conveniently apply theorems of MacDonald and Chen-Louck
(see \cite{CL}[Theorems 3.2,3.3]) concerning the factorial
Schur function. The extra $(-1)^{k(k-1)/2}$ above comes from 
swapping the $i,j$ entry with the $i,k-j+1$ entry.

In our application to $N_k(\alpha;\beta)$, we found, in 
Lemma~\ref{lemma:det 0}, that most terms in~(\ref{eq:N}) can
be discarded, and only terms with
$\alpha_{\sigma_1} = \ldots
= \alpha_{\sigma_{k-|\alpha|}}= \beta_{\tau_1} = \ldots =
\beta_{\tau_{k-|\beta|}}=0$ contribute. 

So assume that, for some $s \geq 0$, $e_1 = \ldots = e_{k-s} = 0$, and
hence $\lambda_{s+1} = \ldots = \lambda_k =0$. Next, write
$c_j = k+j-1 -\epsilon_j$, where $\epsilon_j$ is a non-negative integer.
Assume that, for some $t \geq 0$, $\epsilon_j=0$ for all $j>t$, i.e.
that $c_j=k+j-1$, if $j>t$.

To apply their theorems, one must first assume that the $\lambda_j$'s are
decreasing $\lambda_1 \geq \lambda_2 \ldots$. We can assume this
condition by rearranging the first $s$ columns of the matrix in
(\ref{eq:factorial schur 1}) if necessary. This will change the sign of the
determinant by a power of $(-1)$ that depends only on $\lambda_1,\ldots,\lambda_s$.
However, since we are actually
permuting the $\lambda_j + k-j$'s rather than the $\lambda_j$'s,
some care is needed.

Say $\lambda_j < \lambda_{j+1}$. We are assuming that the
$\lambda_j + k-j$'s are distinct (we've assumed the $e_j+j-1$'s to
be distinct), thus $\lambda_j < \lambda_{j+1}$ actually implies
that
\begin{equation}
    \lambda_j < \lambda_{j+1} - 1
\end{equation}
since otherwise one would have two neighboring $\lambda_j + k-j$'s
that were equal.

Swapping columns $j$ and $j+1$, the subscript for the $j$th column
is then $\lambda_{j+1} + k-j-1 = (\lambda_{j+1}-1)+k-j$, and for
the $(j+1)$st column is then $\lambda_{j} + k-j =
(\lambda_{j}+1)+k-j-1$. Therefore we have replaced
$(\lambda_1,\ldots,\lambda_j,\lambda_{j+1},\ldots)$ with
$(\lambda_1,\ldots,\lambda_{j+1}-1,\lambda_{j}+1,\ldots)$ in which
$\lambda_{j+1}-1 \geq \lambda_{j}+1$. Also notice that swapping
the two columns only permutes the subscripts which therefore
remain distinct.

Continuing in this fashion, we end up with $\tilde{\lambda}_1 \geq
\ldots \geq \tilde{\lambda}_s$, with the $\tilde{\lambda}_i$
obtained from the $\lambda_i$'s by the above swapping procedure.
Notice that
\begin{equation}
    \label{eq:sum lambda}
    \sum \tilde{\lambda}_i = \sum \lambda_i
\end{equation}
since each swap adds one and subtracts one from the $\lambda_i$'s.

\begin{thm}
\label{thm: poly2} Assume that $c_j=k+j-1-\epsilon_j$, with
$\epsilon_j=0$ if $j > t$, that $\tilde{\lambda}_1 \geq \tilde{\lambda}_2 \ldots
\geq \tilde{\lambda}_s$ and that $\tilde{\lambda}_{s+1}=\ldots=\tilde{\lambda}_k=0$. Then
\begin{equation}
  \left|
        (c_i)_{\tilde{\lambda}_j+k-j}
  \right|_{k \times k} = (-1)^{k(k-1)/2} \Delta(c) 
  \text{$\times$(polynomial in $k$ of degree $\leq 2\sum \tilde{\lambda_i}$)}
\end{equation}
\end{thm}

\begin{pf}
By \cite{CL}[Theorems 3.2,3.3,and page 4150] one has
\begin{equation}
  \left|
        (c_i)_{\tilde{\lambda}_j+k-j}
  \right|_{k \times k}
  =
  (-1)^{k(k-1)/2} \Delta(c)
  \left|
        w_{\tilde{\lambda}_i-i+j}(c+j-1)
  \right|_{k \times k}
\end{equation}
with $c+j-1=(c_1+j-1,\ldots,c_k+j-1)$, $\Delta(c) =\prod_{1\le i <
j\le m}(c_j-c_i)$, and
\begin{equation}
    w_m(z) = \sum_{i_1 \leq \ldots \leq i_m \leq k}
    y_{i_1} (y_{i_2}-1)\ldots(y_{i_m}-m+1)
\end{equation}
with $y_i=z_i-i+1$ and the conventions that $w_0(z)=1$, and
$w_m(z)=0$ if $m<0$. The main point is that our $k \times k$
determinant has been replaced by an $s \times s$ determinant.
Furthermore,
\begin{equation}
    w_m(c+j-1) = \sum_{i_1 \leq \ldots \leq i_m \leq k}
    (k+j-1-\epsilon_{i_1})(k+j-2-\epsilon_{i_2})\ldots(k+j-m-\epsilon_{i_m}).
\end{equation}
The terms in this sum with all $i$'s $> t$ contribute
$$
    {k-t \choose m} (k+j-1)\ldots(k+j-m)
$$
since $\epsilon_i=0$ if $i>t$. This is a polynomial of degree $2m$
in $k$.

The terms with $i_1 \leq t$ and the other $i$'s $> t$ contribute
$$
    {k-t \choose m-1} (k+j-2)\ldots(k+j-m)
    \sum_{i_1=1}^t (k+j-1-\epsilon_{i_1}),
$$
and the terms with $i_1 \leq i_2 \leq t$ contribute
$$
    {k-t \choose m-2} (k+j-3)\ldots(k+j-m)
    \sum_{1 \leq i_1 \leq i_2 \leq t} (k+j-1-\epsilon_{i_1})(k+j-2-\epsilon_{i_2}),
$$
both of which are polynomials in $k$ of degree $2m-1$ and $2m-2$ respectively. In
this fashion one sees that the entries of
$$
    \left|
        w_{\tilde{\lambda}_i-i+j}(c+j-1)
    \right|_{s \times s}
$$
are polynomials in $k$. Furthermore, expanding this determinant we get a sum of
products of entries, one from each row and column. However, $w_m(z)=0$ if
$m<0$, so 
only the terms with all $\tilde{\lambda}_i-i+j \geq 0$ contribute,
and the degree of such a term is then $2 \sum \tilde{\lambda}_i$.
\end{pf}

Applying Theorem~\ref{thm: poly2} we find that (\ref{eq:factorial schur 1}) equals
\begin{equation}
    \Delta(c) \prod_{\l=0}^{k-1} \frac{(f_{l+1} + l)!}{l!(k+l)!}
     2^{c_{l+1}-e_{l+1}-l}
     \text{$\times$(polynomial in $k$ of degree $\leq 2\sum {\lambda_i}$)}.
\end{equation}
Here we have absorbed the factor of $\pm 1$ into the polynomial and have also
used (\ref{eq:sum lambda}). 

Thus, given that $\alpha_{\sigma_1}, \alpha_{\sigma_2}+1,\ldots,
\alpha_{\sigma_k}+k-1$ are distinct and $\beta_{\tau_1}, \beta_{\tau_2}+1,\ldots,
\beta_{\tau_k}+k-1$ are distinct, we get that
a typical $\tilde{M}_k(\sigma(\alpha),\tau(\beta))$
appearing in (\ref{eq:N}) equals:
\begin{eqnarray}
     \label{eq:M as poly}
     \tilde{M}_k(\sigma(\alpha),\tau(\beta))
     =
     \Delta(\gamma) \prod_{\l=0}^{k-1} \frac{(\beta_{\tau_{l+1}} + l)!}{l!(k+l)!}
     2^{\gamma_{l+1}-\alpha_{\sigma_{l+1}}-l} \notag \\
     \text{$\times$(polynomial in $k$ of degree $\leq 2 |\alpha|$)}
\end{eqnarray}
where $\gamma$ is the complementary subset of the $\beta_{\tau_i}+i-1$'s.
\begin{equation}
    \set{\gamma_1,\ldots,\gamma_k}
    =
    \set{0,1,\ldots,2k-1}
    -
    \set{
         \beta_{\tau_1},\beta_{\tau_2}+1,\ldots,\beta_{\tau_k}+k-1
    }
\end{equation}
and $\Delta(\gamma) = \prod_{i<j} (\gamma_j-\gamma_i)$. 

Now, by Lemma~\ref{lemma:det 0}, most of the $\alpha_{\sigma_i}$'s and
$\beta_{\tau_i}$'s are 0,
\begin{equation}
    \alpha_{\sigma_1} = \ldots
    = \alpha_{\sigma_{k-|\alpha|}}= \beta_{\tau_1} = \ldots =
    \beta_{\tau_{k-|\beta|}}=0,
\end{equation}
and $\gamma_j=k+j-1-\epsilon_j$ where $\epsilon_j=0$ if $j>|\beta|$.
The latter can be seen by noticing that
$\beta_{\tau_i} + i - 1 \leq |\beta|+k-1$ so that $\set{k+|\beta|,\ldots,2k-1}$
is a subset of $\gamma$. Thus, starting from the end, $\epsilon_k=0$, hence
$\epsilon_{k-1}=0$, $\ldots \epsilon_{|\beta|+1}=0$.

Hence,
\begin{equation}
    \Delta(\gamma) = 
    \prod_{1 \leq i < j \leq k} (j-i+\epsilon_i-\epsilon_j) 
    =
    \prod_{1 \leq i < j \leq k} (j-i) 
    \prod_{1 \leq i < j \leq k} (1+\frac{\epsilon_i-\epsilon_j}{j-i}).
\end{equation}
The first product gives
\begin{equation}
    \prod_{l=1}^k l!
\end{equation}
and, because $\epsilon_i=\epsilon_j=0$ if $i,j>|\beta|$ the second factor equals
\begin{equation}
    \prod_{1 \leq i < j \leq |\beta|} (1+\frac{\epsilon_i-\epsilon_j}{j-i})
    \prod_{{1 \leq i \leq |\beta|} \atop{ |\beta| < j \leq k}} \frac{j-i+\epsilon_i}{j-i}
\end{equation}
However, the product over $1 \leq i < j \leq |\beta|$ is a rational number. Furthermore,
 most of the numerator
of the product over $1 \leq i \leq |\beta|$,$|\beta| < j \leq k$ cancels with the denominator
leaving a polynomial in $k$ of degree $|\epsilon|$.

But $|\beta| = |\epsilon|$ which is easily verified as follows. The union of the $\beta_i+i-1$'s 
and $\gamma_i$'s give $\set{0,1,\dots,2k-1}$, so
\begin{equation}
    \sum_{i=1}^k \beta_i+i-1 +\gamma_i = \sum_{i=0}^{2k-1} i.
\end{equation}
Substituting $\gamma_i = k+i-1-\epsilon_i$ and simplifying gives $\sum \beta_i = \sum \epsilon_i$.

Collecting the above together gives 
\begin{equation}
    \Delta(\gamma) = \prod_{l=1}^k l! 
     \text{$\times$(polynomial in $k$ of degree $|\beta|$)}
\end{equation}
(we take the polyomial to be 1 if $|\beta|=0$).

Next, we determine the power of $2$ appearing in (\ref{eq:M as poly}):
\begin{equation}
     \label{eq:twos}
     \sum_{l=0}^k \gamma_{l+1}-\alpha_{\sigma_{l+1}}-l
     =
     \sum_{l=0}^k k+l-\epsilon_{l+1}-\alpha_{\sigma_{l+1}}-l
     = k^2 +|\beta| - |\alpha|.
\end{equation}
Finally,
\begin{equation}
     \prod_{\l=0}^{k-1} \frac{(\beta_{\tau_{l+1}} + l)!}{l!} =
     \text{polynomial in $k$ of degree $|\beta|$},
\end{equation}
because $\beta_{\tau_1} = \ldots = \beta_{\tau_{k-|\beta|}}=0$, and where we regard 
the l.h.s as a function of $k$ with $\beta_{\tau_{k-j}}$ fixed for $0 \leq j \leq |\beta|$.
Therefore, most of the numerator cancels with the denominator except for $|\beta|$ factors
each of which is a polynomial of degree 1 in $k$. We have therefore shown that 
\begin{equation}
    \tilde{M}_k(\sigma(\alpha),\tau(\beta)) = 
    2^{k^2} \prod_{l=0}^k \frac{l!}{(k+l)!} 
     \text{$\times$(polynomial in $k$ of degree $\leq 2 (|\alpha|+|\beta|$))}.
\end{equation}
Here we have absorbed the extra $2^{|\beta|-|\alpha|}$ from (\ref{eq:twos}) into the polynomial.
This proves that $N_k(\alpha;\beta)$ given by (\ref{eq:N}) is a polynomial in $k$
of degree $\leq  2 (|\alpha|+|\beta|)$.

\section{Numerical evaluation of $c_r(k)$}
\label{section: evaluating c_r}

Two methods were developed to numerically compute the coefficients
$c_r(k)$ of the lower order terms. The first relied on~(\ref{eq:c_r g})
and we used Maple~\cite{M} to take advantage of its symbolic capabilities.
This approach had the advantage of allowing us to obtain the coefficients to
many digits precision, and also to make sense of the conjecture for non-integer
values of $k$. This method suffered the disadvantage of being difficult to implement,
even using a high level symbolic package, and required much computational power, so that
we only determined $c_r(k)$ in this way up to $r \leq 9$. This sufficed to compute all the 
lower terms, for $k=3$, since $P_k(x)=c_0(k) x^{k^2} + c_1(k) x^{k^2-1}+ \ldots
+ c_{k^2}(k)$ is a polynomial in $x$ of degree $k^2$.

The second method was comparatively easy to implement, and allowed us
to obtain many more coefficients. However, it is limited to integer values of $k$,
and also presents more difficulties in acceleration therefore yielding lower precision.

\subsection{Method 1}
\label{section:method1}

A table of the polynomials $N_k(\alpha;\beta)$ of degree $\leq 2(|\alpha|+|\beta|)$
was prepared by evaluating~(\ref{eq:N}) at slightly more than $2(|\alpha|+|\beta|)+1$ values
of $k$ and interpolating the unique polynomial of said degree fitting those points. The
extra points were thrown in for good measure as a check against errors.
This was done for all $0 \leq |\alpha|+|\beta| \leq 9$.

Next a corresponding table of the coefficients $b_k(\alpha;\beta)$ was prepared,
expressed symbolically as a polynomial in the $\gamma_j$'s and $B_k(\alpha;\beta)$'s as
described in
Section~\ref{subsection:residue}, with $B_k(\alpha;\beta)$ given
as a sum over primes, with the summand equal to $\log(p)^{|\alpha|+|\beta|}$ times a
function rational in $p$ and in Gauss hypergeometric functions of
the form ${}_2F_1(k+A,k+B;C;1/p)$, where $A,B,C$ are non-negative integers, $C \geq 1$.
A few example $B_k(\alpha;\beta)$'s are listed in~(\ref{eq:B_k}).

We were then able to obtain, for a given $k$, numerical values of the coefficients
$c_r(k)$, for $0 \leq r \leq 9$. Because $B_k(\alpha;\beta)$
is expressed as an infinite sum over primes, we used standard methods to 
accelerate its convergence. Namely, we evaluated the first few terms, $p\leq P$,
to high precision. Then, to evaluate the tail end of the sum, $p>P$, we used Maple's 
series routine to determined the first few terms
of the series expansion in $1/p$ of the summand, writing it in the form 
\begin{equation}
    \label{eq:tail}
    (\log{p})^r \sum_{j=2}^6 \frac{d_j}{p^j}.
\end{equation}
where $r=|\alpha|+|\beta|$, and the $d_j$'s depend on the summand hence on $\alpha$ and
$\beta$. 

To evaluate a sum of the form
\begin{equation}
    \sum_{p>P} \frac{\log(p)^r}{p^j} 
\end{equation}
we first wrote it as a full sum minus the front end:
\begin{equation}
    \sum_{p} \frac{\log(p)^r}{p^j} - \sum_{p\leq P} \frac{\log(p)^r}{p^j}.
\end{equation}
The second sum was evaluated by summing the terms $p\leq P$, while first sum was
computed using Mobius inversion:
\begin{equation}
    \log\zeta(s) = \sum_{m=1}^\infty \frac{1}{m} \sum_{p} \frac{1}{p^{ms}}, \quad \Re(s)>1,
\end{equation}
so
\begin{equation}
     \sum_{p} \frac{1}{p^{s}} = 
     \sum_{m=1}^\infty
      \frac{\mu(m)}{m} \log\zeta(ms),
\end{equation}
and hence
\begin{equation}
     \label{eq:logzeta}
     \sum_{p} \frac{\log(p)^r}{p^{s}} = 
     (-1)^r
     \sum_{m=1}^\infty
     \frac{\mu(m)}{m} (\log\zeta(ms))^{(r)}.
\end{equation}
Now $(\log\zeta(ms))^{(r)}$ decreases exponentially fast in $m$, 
as can be seen by considering its Dirichlet series which is 
dominated by the first term,
and only a handful of $m$ on the r.h.s. of~(\ref{eq:logzeta}) are needed 
to evaluate $\sum_{p} \frac{\log(p)^r}{p^{s}}$ to a given precision.

The factor $a_k$ given by (\ref{eq:a_k}) can be evaluated to high precision in a similar way.

In this manner we were able to compute $c_r(k)$, $0 \leq r \leq 9$ for various $k$.
For example, for $k=3$ we obtained the coefficients to about $30$
decimal places. The actual precision can be predicited from the size of $P$,
as the overall error in using~(\ref{eq:tail}) to approximate the summand for the terms $p>P$
is $O(\log(P)^{r-1}/P^6)$. In practice, we took larger and larger values of $P$ until
the numerics stabilized to a precision that we found satisfying.

\subsection{Method 2}
\label{section:method2}

The second method we developed to compute $c_r(k)$ used the combinatorial
sum~(\ref{eq: combinatorial sum}), small shifts, and very high precision to capture
cancellation amongst the high order poles of the terms in the sum. Because this method
requires very little symbolically, this was implemented in {\tt C++} using NTL~\cite{S} 
to carry out multiprecision arithmetic. 

The basic idea is as follows. The polynomial $P_k(x)$ given by~(\ref{eq:P_k(x)})
can be regarded as a special case
of the function $P_k(\alpha,x)$ given by~(\ref{eq:P_k(a,x)}), namely with 
$\alpha_1 = \ldots = \alpha_{2k}= 0$. One can then use~(\ref{eq: combinatorial sum}) 
to evaluate
$P_k(\alpha,x)$. However, the terms in~(\ref{eq: combinatorial sum}) have poles
if the $\alpha_i$'s are not distinct.
So, we cannot simply substitute $\alpha = {\bf 0}$ and sum the terms numerically.
Instead we take the limit as $\alpha \to {\bf 0}$ with the condition that
the $\alpha_i$'s are distinct. One must also use very high precision to capture cancellation
amongst the terms which individually become very large when $\alpha$ is small.

More precisely, let
\begin{equation}
    H(z_1,\ldots,z_{2k};x) =
    \exp\left( \frac{x}{2} \sum_1^k z_j-z_{j+k}\right)
    A_k(z_1,\ldots,z_{2k})
    \prod_{i=1}^k\prod_{j=1}^k\zeta(1+z_i-z_{j+k}),
\end{equation}
and let $\epsilon_j = j \epsilon$, where $\epsilon \in {\mathbb C}$.
Then, by (\ref{eq: combinatorial sum})
\begin{equation}
    P_k(x) =
    \lim_{\epsilon \to 0}
    \sum_{\sigma \in \Xi} H( \epsilon_{\sigma(1)},\ldots,\epsilon_{\sigma(2k)};x),
\end{equation}
where $\Xi$ is the set of $\binom{2k}{k}$
permutations $\sigma\in S_{2k}$ such that $\sigma(1)<\cdots <
\sigma(k)$ and $\sigma(k+1)<\cdots < \sigma(2k)$.

Therefore, expanding $\exp$ in its Taylor series, and pulling out the coefficient
of $x^{k^2-r}$, we get
\begin{equation}
    \label{eq:via small shifts}
    c_r(k) =
    \frac{1}{2^{k^2-r} (k^2-r)!}
    \lim_{\epsilon \to 0}
    \sum_{\sigma \in \Xi} H_r( \epsilon_{\sigma(1)},\ldots,\epsilon_{\sigma(2k)}),
\end{equation}
where
\begin{equation}
    \label{eq:H}
    H_r(z_1,\ldots,z_{2k}) =
    \left(\sum_1^k z_j-z_{j+k}\right)^{k^2-r}
    A_k(z_1,\ldots,z_{2k})
    \prod_{i=1}^k\prod_{j=1}^k\zeta(1+z_i-z_{j+k}).
\end{equation}
The only complication in evaluating the above for a given $k$ and $\epsilon$ 
is that $A_k(z_1,\ldots,z_{2k})$ is expressed as an infinite product over 
primes~(\ref{eq:A_k}). To evaluate that product, we broke it up into
$p \leq P$ and $p>P$, with $P$ large. For the primes $p\leq P$, we
used~(\ref{eq:A_k euler product}) to evaluate the contribution 
from $p$, each factor only requiring finitely many arithmetic steps.

For the contribution from the larger primes, $p>P$, we used
a quadratic approximation for the local factor appearing 
in~(\ref{eq:A_k}):
\begin{equation}
    \label{eq:local factor deg 2}
    1- \sum_{{1 \leq i_1<i_2 \leq k}\atop {1 \leq j_1<j_2 \leq k}}
    p^{-2-z_{i_1}-z_{i_2}+z_{k+j_1}+z_{k+j_2}}.
\end{equation}
This approximation can be obtained by substituting $u_j=p^{-1/2-z_j}$ and 
$w_j=p^{-1/2+z_{k+j}}$ in the local factor of~(\ref{eq:A_k}),
and working out the terms up to degree four. 
Only terms of even degree appear because the intergral over $\theta$ pulls
out just the terms with the same number of $u$'s and $w$'s.
So, we expand each geometric series appearing in the integral
over $\theta$ up to degree two, multiply them out, and collect terms with 
the same number of $u$'s and $w$'s. 

The first factor $\prod (1-u_iw_j)$ appears precisely to cancel
the terms of degree two in the second factor, so we need only determine
the terms of degree four. Noticing that the local factor is symmetric
separately in the $u$'s and $w$'s, and also if we swap $u$ and $w$,
we can determine the terms of degree four 
by simply computing all representative fourth order partial derivatives
that have the same number of $u$'s as $w$'s,
evaluated at $u_j=w_j=0$, $1 \leq j \leq k$.
For instance, it is enough to immediately set $u_j=w_j=0$ if $3 \leq j \leq k$ and then
take the partial derivatives:
$\frac{\partial}{\partial u_1}\frac{\partial}{\partial u_2}
\frac{\partial}{\partial w_1}\frac{\partial}{\partial w_2}$,
$\frac{\partial}{\partial u_1}\frac{\partial}{\partial u_2}\frac{\partial^2}{\partial w_1^2}$,
$\frac{\partial^2}{\partial u_1^2}\frac{\partial^2}{\partial w_1^2}$,
evaluated at $u_1=u_2=w_1=w_2=0$. Doing so gives that the local factor, up to terms of degree four,
equals
\begin{equation}
    1- \sum_{{1 \leq i_1<i_2 \leq k}\atop {1 \leq j_1<j_2 \leq k}}
    u_{i_1}u_{i_2}w_{j_1}w_{j_2},
\end{equation}
thus giving (\ref{eq:local factor deg 2}).

Therefore, we used
\begin{eqnarray}
    &&\prod_{p>P} \prod_{{1 \leq i_1<i_2 \leq k}\atop {1 \leq j_1<j_2 \leq k}}
    (1-p^{-2-z_{i_1}-z_{i_2}+z_{k+j_1}+z_{k+j_2}})
    = \notag \\
    &&\frac{
        \prod_{{1 \leq i_1<i_2 \leq k}\atop {1 \leq j_1<j_2 \leq k}} 
        \zeta(2+z_{i_1}+z_{i_2}-z_{k+j_1}-z_{k+j_2})^{-1}
    }
    {
        \prod_{p \leq P} 
        \prod_{{1 \leq i_1<i_2 \leq k}\atop {1 \leq j_1<j_2 \leq k}}
        (1-p^{-2-z_{i_1}-z_{i_2}+z_{k+j_1}+z_{k+j_2}})
    }
\end{eqnarray}
to approximate the contribution to (\ref{eq:A_k}) for the primes $p>P$.

To compute ~(\ref{eq:via small shifts}) to $D$ decimal places, we should take
$\epsilon$ roughly of size $10^{-D}$ and then use about $(r+1) D$ digits working precision
to account for
cancellation in~(\ref{eq:via small shifts}) amongst the
order $r$ pole in $\epsilon$ of the summands. Since $r$ can be
as large as $k^2$, we used $(k^2+8) D$ digits working precision,
the $+8$ taken for extra leeway, and also chose $D$ to be slightly
larger than the desired final precision.

\section{Verifying the full moment conjecture}
\label{section:data}

In ~\cite{CFKRS} we presented numerical data supporting the conjecture
described in Section~\ref{subsection: moment conjectures zeta} for
$k=3$. Here we give some more data supporting the conjecture, for 
integer $k \leq 7$, and also for several real and complex values of $k$. 
Our data supports the conjecture,
but is not too extensive as our main effort was put towards developing
ways to evaluate the lower terms rather than to large scale verification
of the conjecture. Nonetheless, even moderate data strongly supports
the conjecture.

One experiment we carried out involved comparing the two quantitites
\begin{equation}
    \label{eq:moment zeta}
    \int_C^D |\zeta(\tfrac12 +it)|^{2k} dt
\end{equation}
and
\begin{equation}
     \label{eq:P_k integral}
     \int_C^D P_k(\log(t/2\pi)) dt,
\end{equation}
for seventeen intervals $[C,D]$ of length $50000$, and $k=3,4,5,6,$ and $7$.
We also examined the conjecture for several non-integer values of $k$, in
the latter case interpreting $P_k(x)$ as an asymptotic series rather
than as polynomial.

For both integer and non-integer $k$, we used method 1
of Section~\ref{section:method1} to compute
the first few lower terms to high precision.
For $k=3,4,5,6,7,8$ we also computed $c_r(k)$ for $r\leq k^2$ using method 2.

Tables of the coefficients $c_r(k)$, $k=3,4,5,6,7$ can be found in~\cite{CFKRS}.
As $k$ increases, it seems from numerics that the first few leading order terms
have much smaller coefficients than the later lower order terms. For example, when
$k=4$, the leading term
as listed in~\cite{CFKRS} is $c_0(4)= .24650183919342276\times 10^{-12}$, compared to
the largest value $c_{13}(4)=38.203306$. To verify the full moment conjecture,
we needed to evaluate $P_k(x)$ at $x=\log(t/(2\pi))$. Therefore, even if $t$ is
moderately large, say $10^6$, the main contribution actually comes from substantially
lower order terms. In the range we examined, the main contribution for $k=3,4,5,6,7$ came,
respectively, from the the terms $r=3,6,12,20,30$ and their immediate neighbours.

To compare to actual moment data for $\zeta$ we used Mathematica~\cite{Ma}
to numerically integrate powers of $|\zeta(1/2+it)|$ for each interval.
Due the oscillatory nature of $\zeta$,
we performed the integration between consecutive zeros
of $\zeta$ on the halfline using a table of zeros computed with the $L$-function
calculator~\cite{R}.

Tables \ref{table:345} and \ref{table:67} gives the values of~(\ref{eq:P_k integral}) 
and~(\ref{eq:moment zeta}) 
for $[C,D]=[50000n,50000(n+1)]$, $n=0,1,\ldots,16$, and $k=3,4,5,6,7$. 
The data for $k=3$ is a subset of the data given in~\cite{CFKRS}, but otherwise,
the data here is new. We see that the pairs of columns track one another nicely.

\begin{table}[h!tb]
\centerline{
\begin{tabular}{|c||c|c||c|c||c|c||}
\hline
$n$& conj. $k=3$ & data $k=3$ & conj. $k=4$ & data $k=4$& conj. $k=5$ & data $k=5$ \cr
\hline
0 & $7.23687 \times 10^{9}$ & $7.23101 \times 10^{9}$ &	$1.89527 \times 10^{12}$  & $1.88501 \times 10^{12}$ &	$6.00428 \times 10^{14}$ &   $5.91051 \times 10^{14}$ \cr
1 & $1.56965 \times 10^{10}$ & $1.57239 \times 10^{10}$ &	$5.67575 \times 10^{12}$  & $5.70833 \times 10^{12}$ &	$2.45298 \times 10^{15}$ &   $2.47886 \times 10^{15}$ \cr
2 & $2.15687 \times 10^{10}$ & $2.15368 \times 10^{10}$ &	$9.17127 \times 10^{12}$  & $9.12987 \times 10^{12}$ &	$4.68619 \times 10^{15}$ &   $4.64908 \times 10^{15}$ \cr
3 & $2.63814 \times 10^{10}$ & $2.62463 \times 10^{10}$ &	$1.24573 \times 10^{13}$  & $1.23432 \times 10^{13}$ &	$7.10198 \times 10^{15}$ &   $7.04187 \times 10^{15}$ \cr
4 & $3.05562 \times 10^{10}$ & $3.06922 \times 10^{10}$ &	$1.55847 \times 10^{13}$  & $1.5683  \times 10^{13}$ &	$9.63318 \times 10^{15}$ &   $9.6445  \times 10^{15}$ \cr
5 & $3.42903 \times 10^{10}$ & $3.44143 \times 10^{10}$ &	$1.8585 \times 10^{13}$  & $1.87265 \times 10^{13}$ &	$1.22457 \times 10^{16}$ &   $1.24349 \times 10^{16}$ \cr
6 & $3.76958 \times 10^{10}$ & $3.76835 \times 10^{10}$ &	$2.14798 \times 10^{13}$  & $2.15861 \times 10^{13}$ &	$1.4919 \times 10^{16}$ &   $1.51619 \times 10^{16}$ \cr
7 & $4.08439 \times 10^{10}$ & $4.05663 \times 10^{10}$ &	$2.42845 \times 10^{13}$  & $2.37201 \times 10^{13}$ &	$1.76398 \times 10^{16}$ &   $1.66972 \times 10^{16}$ \cr
8 & $4.37832 \times 10^{10}$ & $4.39075 \times 10^{10}$ &	$2.70108 \times 10^{13}$  & $2.724   \times 10^{13}$ &	$2.03988 \times 10^{16}$ &   $2.06017 \times 10^{16}$ \cr
9 & $4.65486 \times 10^{10}$ & $4.65312 \times 10^{10}$ &	$2.96679 \times 10^{13}$  & $2.94271 \times 10^{13}$ &	$2.3189 \times 10^{16}$ &   $2.26023 \times 10^{16}$ \cr
10 & $4.91663 \times 10^{10}$ & $4.91363 \times 10^{10}$ &	$3.22631 \times 10^{13}$ &  $3.24807 \times 10^{13}$ &	$2.60051 \times 10^{16}$ &   $2.69184 \times 10^{16}$ \cr
11 & $5.16565 \times 10^{10}$ & $5.17448 \times 10^{10}$ &	$3.48022 \times 10^{13}$ &  $3.47606 \times 10^{13}$ &	$2.88433 \times 10^{16}$ &   $2.87018 \times 10^{16}$ \cr
12 & $5.40352 \times 10^{10}$ & $5.39624 \times 10^{10}$ &	$3.72905 \times 10^{13}$ &  $3.73482 \times 10^{13}$ &	$3.17002 \times 10^{16}$ &   $3.18035 \times 10^{16}$ \cr
13 & $5.63152 \times 10^{10}$ & $5.65418 \times 10^{10}$ &	$3.97319 \times 10^{13}$ &  $4.00187 \times 10^{13}$ &	$3.45733 \times 10^{16}$ &   $3.48184 \times 10^{16}$ \cr
14 & $5.85072 \times 10^{10}$ & $5.83654 \times 10^{10}$ &	$4.21303 \times 10^{13}$ &  $4.1917  \times 10^{13}$ &	$3.74603 \times 10^{16}$ &   $3.70813 \times 10^{16}$ \cr
15 & $6.062 \times 10^{10}$ & $6.08708 \times 10^{10}$ &	$4.44887 \times 10^{13}$ &  $4.48257 \times 10^{13}$ &	$4.03594 \times 10^{16}$ & $4.08236 \times 10^{16}$ \cr
16 & $6.2661 \times 10^{10}$ & $6.27652 \times 10^{10}$ &	$4.68097 \times 10^{13}$ &  $4.69566 \times 10^{13}$ &	$4.32693 \times 10^{16}$ & $4.3287  \times 10^{16}$ \cr
%17 & $6.46366 \times 10^{10}$ & $6.42272 \times 10^{10}$ &	$4.9096 \times 10^{13}$ &  $4.86013 \times 10^{13}$ &		          &               \cr
%didn't have data for k=5, n=17 so commented out the last line
\hline
\end{tabular}
}
\caption{
This table compares the conjectured value (\ref{eq:P_k integral}) to actual data (\ref{eq:moment zeta})
for intervals $[50000n,50000(n+1)]$, $n=0,1,\ldots,16$, and $k=3,4,5$. The fit is to two or three decimal places,
consistent with the remainder stated in~(\ref{eq: moments zeta}).
}
\label{table:345}
\end{table}

\begin{table}[h!tb]
\centerline{
\begin{tabular}{|c||c|c||c|c||}
\hline
$n$& conj. $k=6$ & data $k=6$ & conj. $k=7$ & data $k=7$ \cr
\hline
0 &	$2.15456 \times 10^{17}$ &   $2.08527 \times 10^{17}$ &	$8.45652 \times 10^{19}$ &   $7.99015 \times 10^{19}$ \cr
1 &	$1.18835 \times 10^{18}$ &   $1.20686 \times 10^{18}$ &	$6.24627 \times 10^{20}$ &  $6.3773  \times 10^{20}$\cr
2 &	$2.69034 \times 10^{18}$ &   $2.66481 \times 10^{18}$ &	$1.67709 \times 10^{21}$ &   $1.66563 \times 10^{21}$ \cr
3 &	$4.56155 \times 10^{18}$ &   $4.56713 \times 10^{18}$ &	$3.18661 \times 10^{21}$ &   $3.25679 \times 10^{21}$ \cr
4 &	$6.72399 \times 10^{18}$ &   $6.61933 \times 10^{18}$ &	$5.1125  \times 10^{21}$ &   $4.87831 \times 10^{21}$ \cr
5 & $9.12928 \times 10^{18}$ &   $9.3828  \times 10^{18}$ &	$7.42365 \times 10^{21}$ &   $7.74635 \times 10^{21}$ \cr
6 &	$1.17439 \times 10^{19}$ &   $1.21474 \times 10^{19}$ &	$1.00952 \times 10^{22}$ &   $1.06992 \times 10^{22}$ \cr
7 &	$1.45431 \times 10^{19}$ &   $1.31266 \times 10^{19}$ &	$1.31065 \times 10^{22}$ &   $1.11053 \times 10^{22}$ \cr
8 &	$1.75076 \times 10^{19}$ &   $1.75386 \times 10^{19}$ &	$1.64403 \times 10^{22}$ &   $1.61306 \times 10^{22}$ \cr
9 &	$2.06221 \times 10^{19}$ &   $1.95439 \times 10^{19}$ &	$2.00815 \times 10^{22}$ &   $1.83038 \times 10^{22}$ \cr
10 &	$2.3874  \times 10^{19}$ &   $2.61353 \times 10^{19}$ &	$2.40171 \times 10^{22}$ &   $2.8627  \times 10^{22}$ \cr
11 &	$2.72527 \times 10^{19}$ &   $2.70986 \times 10^{19}$ &	$2.82354 \times 10^{22}$ &   $2.82074 \times 10^{22}$ \cr
12 &	$3.07492 \times 10^{19}$ &   $3.06639 \times 10^{19}$ &	$3.2726  \times 10^{22}$ &   $3.20372 \times 10^{22}$ \cr
13 &	$3.43557 \times 10^{19}$ &   $3.43848 \times 10^{19}$ &	$3.74797 \times 10^{22}$ &   $3.70176 \times 10^{22}$ \cr
14 &	$3.80656 \times 10^{19}$ &   $3.7414  \times 10^{19}$ &	$4.24878 \times 10^{22}$ &   $4.13975 \times 10^{22}$ \cr
15 &	$4.18729 \times 10^{19}$ &   $4.25286 \times 10^{19}$ &	$4.77427 \times 10^{22}$ &   $4.86676 \times 10^{22}$ \cr
16 &	$4.57724 \times 10^{19}$ &   $4.53193 \times 10^{19}$ &	$5.32373 \times 10^{22}$ &   $5.1628  \times 10^{22}$ \cr
17 & $4.97592 \times 10^{19}$ &   $4.98651 \times 10^{19}$ &	$5.89648 \times 10^{22}$ &   $6.0058  \times 10^{22}$ \cr
\hline
\end{tabular}
}
\caption{Conjecture v.s. data for $k=6,7$, same intervals as the previous table.
}
\label{table:67}
\end{table}

\begin{figure}[h!tb]
    \centerline{
            \psfig{figure=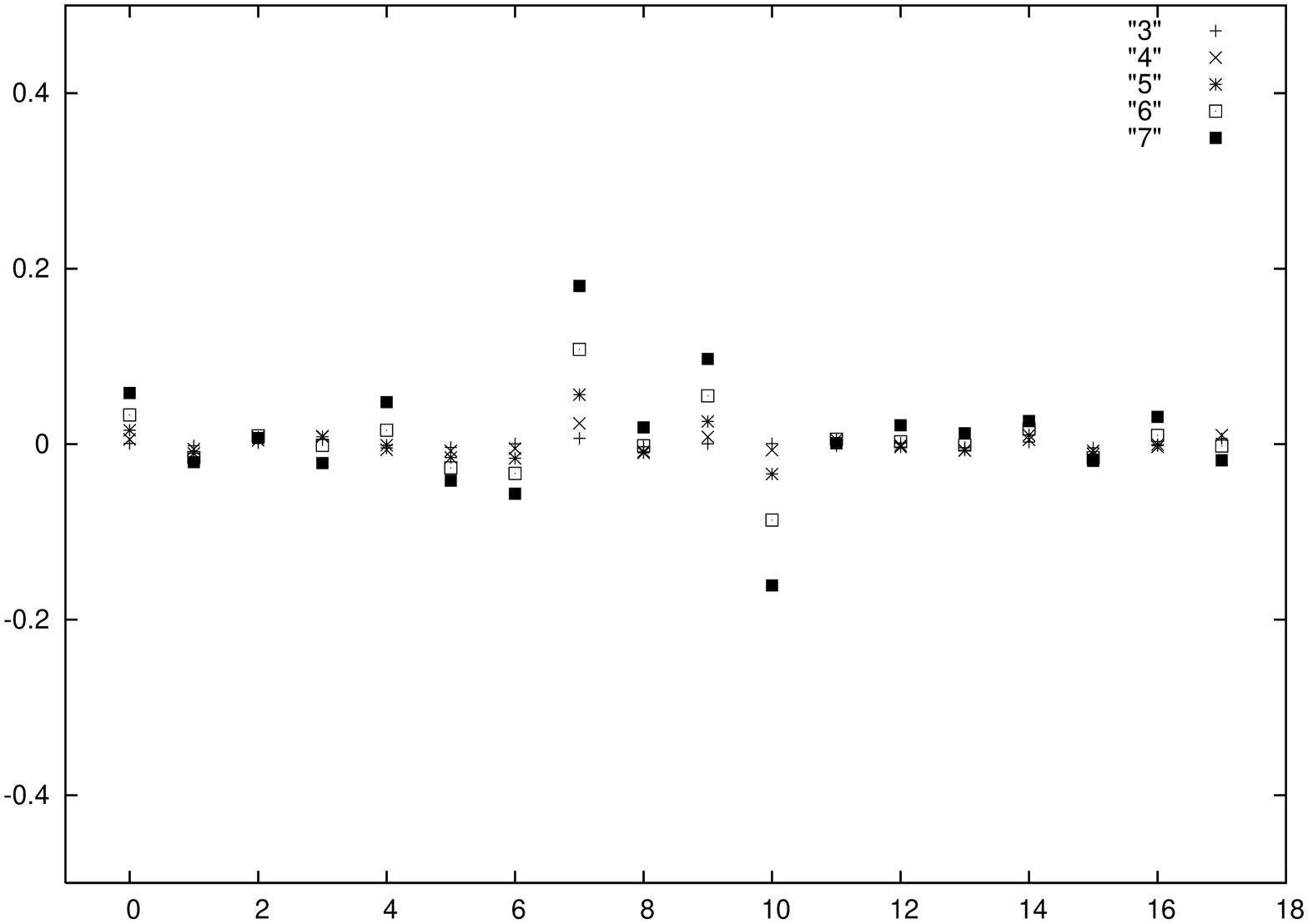,width=5in,angle=-90}
    }
    \caption
    {
         The horizontal axis is $n$, and the vertical axis depicts
         $(\text{conjecture}-\text{data})/\text{data}$
         for the values in Tables~\ref{table:345} and~\ref{table:67}. 
         For any finite interval, as $k\to\infty$, the main contribution to
         the $2k$th moment comes from the largest value of $|\zeta(1/2+it)|^{2k}$ 
         on that interval. This explains the feature that, for a fixed interval,
         the actual moment tends to progressively deviate from the conjectured value as 
         $k$ increases.
    }
    \label{fig:conj vs reality}
\end{figure}

Figure \ref{fig:conj vs reality} depicts the difference 
between~(\ref{eq:P_k integral}) and~(\ref{eq:moment zeta})
divided by~(\ref{eq:moment zeta}) for the values in Tables~\ref{table:345} 
and~\ref{table:67}.

We also present some data for non-integer $k$, specifically for $k=.5, 1.8, 3.2,$ and $.5+i$
taking just the first first few terms, $r \leq 7$, of $P_k(x)$.
For non-integer $k$ we believe, based on our numerics, that $P_k(x)$, no longer a polynomial in $x$ but
an infinite series, gives an asymptotic expansion for the $2k$th moment of $\zeta$, so that
taking more terms does not necessarily give an improvement.
We therefore compared~(\ref{eq:moment zeta}) to
\begin{equation}
     \label{eq:asympt series}
     \int_{100}^D
     \sum_{r=0}^R c_r(k) \log(t/2\pi)^{k^2-r} dt,
\end{equation}
for a few values of $D$. We present our data in Tables~\ref{table:.5}--\ref{table:.5+i},
listing for each $k$, the values of $c_0(k),\ldots,c_7(k)$, and of~(\ref{eq:moment zeta})
compared to~(\ref{eq:asympt series}) with $D=1000,10000,$ and $100000$, and $R=0,1,\ldots7$.

\begin{table}[h!tb]
\centerline{
\begin{tabular}{|c|c|c|c|c|}
\hline
$R$& $c_R(.5)$ & (\ref{eq:asympt series}), $D=1000$ 
& (\ref{eq:asympt series}), $D=10000$ & (\ref{eq:asympt series}), $D=100000$ \cr
\hline
$0$ &
$1.1299287453321533$ &
$1463.83$&
$17768.4$&
$193494.$ \cr
$1$ &
$.19628236755422853$&
$1523.55$&
$18258.1$&
$197413.$ \cr
$2$ &
$.03248602185728907$&
$1525.93$&
$18271.4$&
$197491.$ \cr
$3$ &
$-.5289095729314908$&
$1516.37$&
$18234.3$&
$197335.$ \cr
$4$ &
$3.2346669444094671$&
$1531.24$&
$18275.5$&
$197459.$ \cr
$5$ &
$-21.381296730027876$&
$1505.43$&
$18222.2$&
$197343.$ \cr
$6$ &
$166.38844209028643$&
$1559.87$&
$18310.8$&
$197488.$ \cr
$7$ &
$-1529.2695739774642$&
$1419.97$ &
$18120.1$ &
$197237.$ \cr
\hline
(\ref{eq:moment zeta})&
&
$1521.27$&
$18257.1$&
$197425.$ \cr
\hline
\end{tabular}
}
\caption{The coefficients $c_R(k)$, and conjecture v.s. data for $k=.5$ for three intervals.
The bottom row gives~(\ref{eq:moment zeta}) for the 
interval $[100,D]$, with $D=1000,10000,$ and $100000$. For each $D$, we compare this to 
the value of~(\ref{eq:asympt series}), $R=0,1,\ldots,7$.
}
\label{table:.5}
\end{table}

%k=.5+0.*I; P=5881
%c[0]=1.12992874533215334896507788143+ 0.*I;
%c[1]=.196282367554228531756982002928+ 0.*I;
%c[2]=.324860218572890685155131684573e-1+ 0.*I;
%c[3]=-.528909572931490802242248419890+ 0.*I;
%c[4]=3.23466694440946714875500352548+ 0.*I;
%c[5]=-21.3812967300278759861032468564+ 0.*I;
%c[6]=166.388442090286430079715815635+ 0.*I;
%c[7]=-1529.26957397746419926521869346+ 0.*I;
%===========================================================
%k=.5+0.*I; P=5897
%c[0]=1.12992874533215334896507788137+ 0.*I;
%c[1]=.196282367554228531756982005851+ 0.*I;
%c[2]=.324860218572890685155132528706e-1+ 0.*I;
%c[3]=-.528909572931490802242256652116+ 0.*I;
%c[4]=3.23466694440946714875530208350+ 0.*I;
%c[5]=-21.3812967300278759861161563834+ 0.*I;
%c[6]=166.388442090286430080368947572+ 0.*I;
%c[7]=-1529.26957397746419930115497087+ 0.*I;

\begin{table}[h!tb]
\centerline{
\begin{tabular}{|c|c|c|c|c|}
\hline
$R$& $c_R(3.2)$ & (\ref{eq:asympt series}), $D=1000$ 
& (\ref{eq:asympt series}), $D=10000$ & (\ref{eq:asympt series}), $D=100000$ \cr
\hline
$0$ &
$.37531596173465401\times 10^{-6}$&
$1968.83$&
$1.16353\times 10^6$&
$2.19960\times 10^8$ \cr
$1$&
$.34462154217944847\times 10^{-4}$&
$40049.5$&
$1.65169\times 10^7$&
$2.41753\times 10^9$ \cr
$2$&
$.12662390083082525\times 10^{-2}$&
$336190.$ &
$9.78885 \times 10^7$&
$1.12289 \times 10^{10}$ \cr
$3$&
$.23963666452208821 \times 10^{-1}$&
$1.52891 \times 10^6$&
$3.2097 \times 10^8$&
$2.94868 \times 10^{10}$ \cr
$4$&
$.2526426167678357$&
$4.22213 \times 10^6$&
$6.6336 \times 10^8$&
$5.06417 \times 10^{10}$ \cr
$5$&
$1.5214668466274718$&
$7.72205 \times 10^6$&
$9.6529 \times 10^8$&
$6.47041 \times 10^{10}$\cr
$6$&
$5.3060442651520751$&
$1.03793 \times 10^7$&
$1.12055 \times 10^9$&
$7.01449 \times 10^{10}$\cr
$7$&
$11.121264784324178$&
$1.16045 \times 10^7$&
$1.16894 \times 10^9$&
$7.14177 \times 10^{10}$\cr
\hline
(\ref{eq:moment zeta})&
&
$1.15305 \times 10^7$&
$1.16746 \times 10^9$&
$7.16886 \times 10^{10}$\cr
\hline
\end{tabular}
}
\caption{Conjecture v.s. data for $k=3.2$.}
\label{table:3.2}
\end{table}

%k=3.2+0.*I; P=8641
%c[0]=.375315961734654013338564668342e-6+ 0.*I;
%c[1]=.344621542179448465558500172816e-4+ 0.*I;
%c[2]=.126623900830825254032278723637e-2+ 0.*I;
%c[3]=.239636664522088214444421394139e-1+ 0.*I;
%c[4]=.252642616767835704479426328016+ 0.*I;
%c[5]=1.52146684662747180409355162430+ 0.*I;
%c[6]=5.30604426515207512371591255431+ 0.*I;
%c[7]=11.1212647843241776538677587841+ 0.*I;
%===========================================================
%k=3.2+0.*I; P=8647
%c[0]=.375315961734654013338565179668e-6+ 0.*I;
%c[1]=.344621542179448465558498566217e-4+ 0.*I;
%c[2]=.126623900830825254032280848737e-2+ 0.*I;
%c[3]=.239636664522088214444406098381e-1+ 0.*I;
%c[4]=.252642616767835704479490539218+ 0.*I;
%c[5]=1.52146684662747180409202361923+ 0.*I;
%c[6]=5.30604426515207512373217982415+ 0.*I;
%c[7]=11.1212647843241776539394200913+ 0.*I;

\begin{table}[h!tb]
\centerline{
\begin{tabular}{|c|c|c|c|c|}
\hline
$R$& $c_R(1.8)$ & (\ref{eq:asympt series}), $D=1000$ 
& (\ref{eq:asympt series}), $D=10000$ & (\ref{eq:asympt series}), $D=100000$ \cr
\hline
$0$ &
$.13885991555298723$&
$15298.6$&
$604203.$&
$1.58922\times 10^7$\cr
$1$&
$1.2590684761107478$&
$46198.5$&
$1.42746\times 10^6$&
$3.20931\times 10^7$\cr
$2$&
$2.4174835075472416$&
$59612.2$&
$1.66821\times 10^6$&
$3.56224\times 10^7$\cr
$3$&
$2.546894763686222$&
$62863.9$&
$1.70753\times 10^6$&
$3.60492\times 10^7$\cr
$4$&
$-2.21426710514627$&
$62199.9$&
$1.7021 \times 10^6$&
$3.60059\times 10^7$\cr
$5$&
$3.223904454789757$&
$62432.5$&
$1.70339\times 10^6$&
$3.60134\times 10^7$\cr
$6$&
$46.42674651960987$&
$63260.1$&
$1.70659\times 10^6$&
$3.60268\times 10^7$\cr
$7$&
$-840.1304443557953$&
$59448.3$&
$1.69608\times 10^6$&
$3.59953\times 10^7$\cr
\hline
(\ref{eq:moment zeta})&
&
$61744.5$&
$1.70134\times 10^6$&
$3.60129\times 10^7$\cr
\hline
\end{tabular}
}
\caption{Conjecture v.s. data for $k=1.8$. For this value of $k$, and the
range we examined, $R=4$ or $5$ give the best approximation.}
\label{table:1.8}
\end{table}

%k=1.8+0.*I; P=3413
%c[0]=.138859915552987225587322540551+ 0.*I;
%c[1]=1.25906847611074783031420101243+ 0.*I;
%c[2]=2.41748350754724162283139915537+ 0.*I;
%c[3]=2.54689476368622216280234473277+ 0.*I;
%c[4]=-2.21426710514626977265767967094+ 0.*I;
%c[5]=3.22390445478975732431422486753+ 0.*I;
%c[6]=46.4267465196098721224717742994+ 0.*I;
%c[7]=-840.130444355795276779436787807+ 0.*I;
%===========================================================
%k=1.8+0.*I; P=3433
%c[0]=.138859915552987225587322538434+ 0.*I;
%c[1]=1.25906847611074783031420142772+ 0.*I;
%c[2]=2.41748350754724162283136840942+ 0.*I;
%c[3]=2.54689476368622216280340743486+ 0.*I;
%c[4]=-2.21426710514626977267679871137+ 0.*I;
%c[5]=3.22390445478975732466258230657+ 0.*I;
%c[6]=46.4267465196098721125224099795+ 0.*I;
%c[7]=-840.130444355795276619171861407+ 0.*I;

%D=1000000, k=1.8, data not presented:
%6.59680704430408*10^8 (data)
%6.596476379392655*10^8 (R=7)
%6.597662893548957*10^8 (R=6)
%6.59699441919913*10^8 (etc)
%6.596505456873596*10^8
%6.600114479116199*10^8
%6.554895064803236*10^8
%6.082661833231634*10^8
%3.3556607623884577*10^8

\begin{table}[h!tb]
\centerline{
\begin{tabular}{|c|c|c|c|c|}
\hline
$R$& $c_R(.5+i)$ & (\ref{eq:asympt series}), $D=1000$ 
& (\ref{eq:asympt series}), $D=10000$ & (\ref{eq:asympt series}), $D=100000$ \cr
\hline
$0$&
$1.3117481341987813+1.211708767666727i$&
$-308.872+ 439.126i $&
$-3698.00+ 2357.78i $&
$-34129.1+ 8908.71i $\cr
$1$&
$-3.0693034820213132+2.309977688777579i$&
$-508.454 + 246.589i $&
$-4331.26+ 957.574i $&
$-35042.9 - 44.9533i $\cr
$2 $&
$23.861826126198446-5.4045694962616631i$&
$-335.035 + 646.347i $&
$-4243.1 + 2618.59i $&
$-36981.1+ 6688.15i $\cr
$3$&
$-111.54278536885322-35.79807241977336i$&
$-285.625 + 118.290i$&
$-3667.307 + 1304.32i$&
$-34054.1 + 3538.21i$\cr
$4$&
$828.16689710582718+437.514818042632i$&
$-546.679 + 1199.98i $&
$-4747.02 + 3257.57i $&
$-37543.5 + 6786.35i $\cr
$5$&
$-5808.11341189128-8339.592888954564i$&
$1514.537- 1342.42i $&
$-738.290- 0.1097i $&
$-30111.2 + 3451.52i $\cr
$6$&
$15613.29091863494+101218.4464636376i$&
$-6736.01 + 2796.92i $&
$-12358.2 + 3830.64i $&
$-45377.2 + 5624.57i $\cr
$7$&
$188541.27977634034-1175857.723687032i$&
$23708.43 - 2789.22i $&
$24043.14 + 702.364i $&
$-4932.63 + 6133.48i $\cr
\hline
(\ref{eq:moment zeta})&
&
$-340.843 + 383.859i$&
$-3946.25 + 1883.17i$&
$-35140. + 4830.47i $\cr
\hline
\end{tabular}
}
\caption{Conjecture v.s. data for $k=.5+i$. The data here is not as convincing
as for the other values of $k$, but, nonetheless, the early terms do give
a reasonable approximation, and we believe the fit would improve with more
substantial data.}
\label{table:.5+i}
\end{table}

\subsection{Zeros of $P_k(x)$}
\label{section:zeros}

While we have managed to explicitly determine the first few coefficients $c_r(k)$ of the moment
polynomials $P_k(x)$, we have not yet managed to understand certain aspects of these polynomials,
such as the uniform asymptotics of the coefficients, or uniform asymptotics of $P_k(x)$
with $x$ a function of $k$. The latter is needed, for example, to properly understand
how large $|\zeta(1/2+it)|$ can get~\cite{FGH}.

%Since the overall nature of these polynomails remains mysterious, we depict, in Figure~\ref{fig:zeros}
%the $k^2$ zeros of $P_k(x)$ for $k=2,3,4,5,6,7$ in the hope that a reader might recognize some familiar
%behaviour.
%
%\begin{figure}
%\centerline{
%    \psfig{figure=2.ps,width=3.75in,angle=270}
%    \psfig{figure=3.ps,width=3.75in,angle=270}
%}
%\centerline{
%    \psfig{figure=4.ps,width=3.75in,angle=270}
%    \psfig{figure=5.ps,width=3.75in,angle=270}
%}
%\centerline{
%    \psfig{figure=6.ps,width=3.75in,angle=270}
%    \psfig{figure=7.ps,width=3.75in,angle=270}
%}
%
%\caption[zeros]
%            {The $k^2$ zeros of the $\zeta$ moment polynomials $P_k(x)$, $2 \leq k \leq 7$}
%\label{fig:zeros}
%\end{figure}

\subsection{Plots of $c_r(k)$}
\label{section:graphs c_r}

\begin{figure}
\centerline{
    \psfig{figure=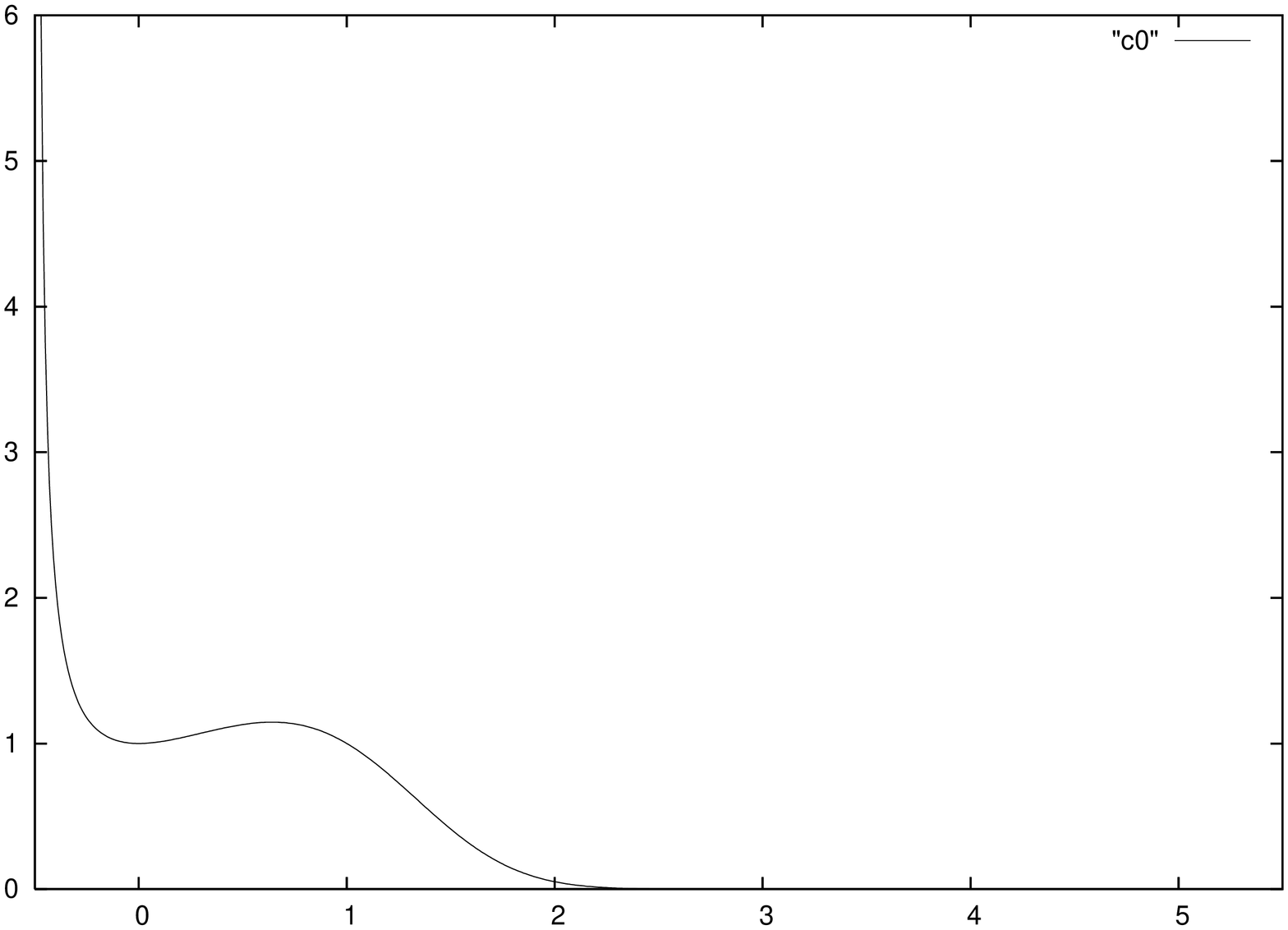,width=2.75in,angle=270}
    \psfig{figure=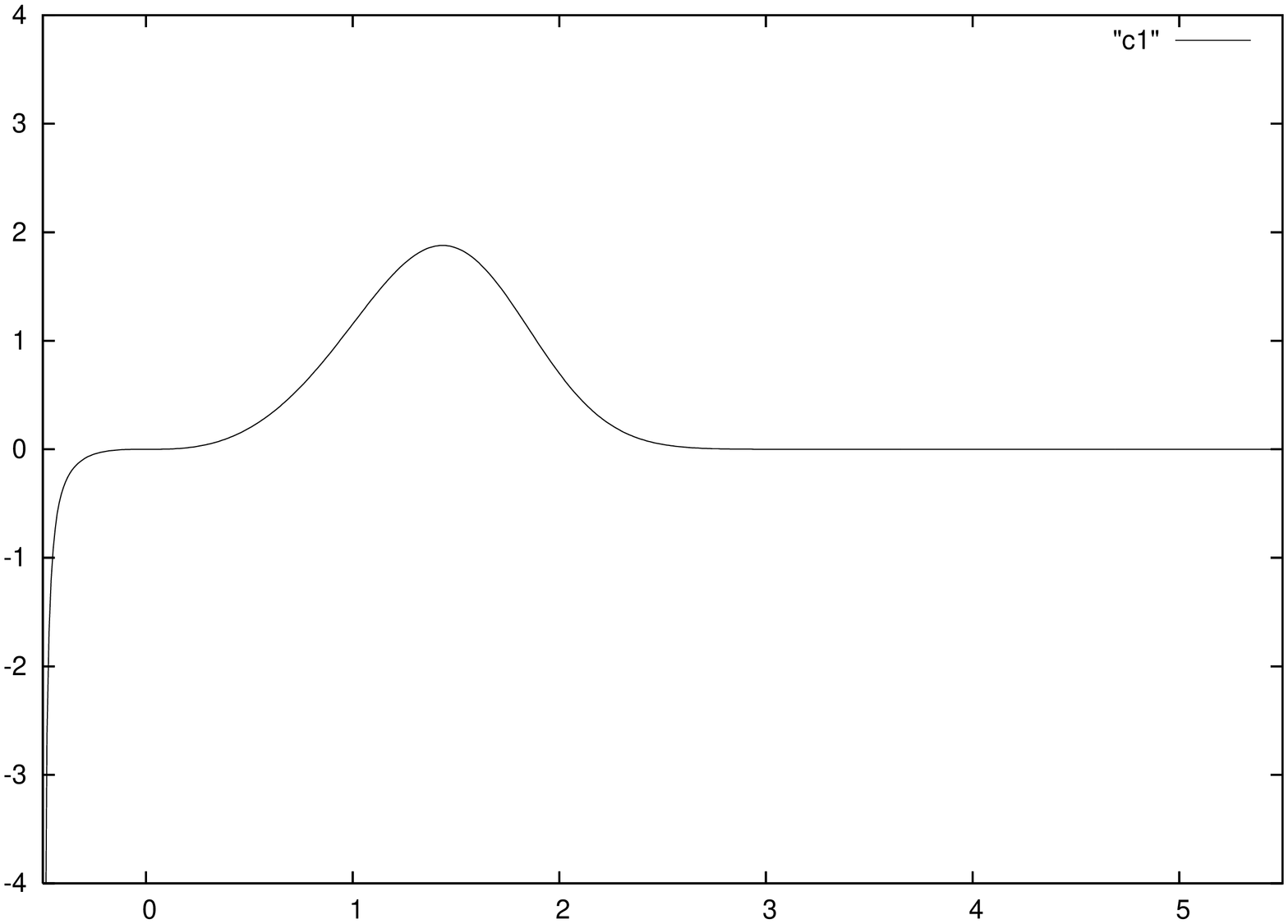,width=2.75in,angle=270}
}
\centerline{
    \psfig{figure=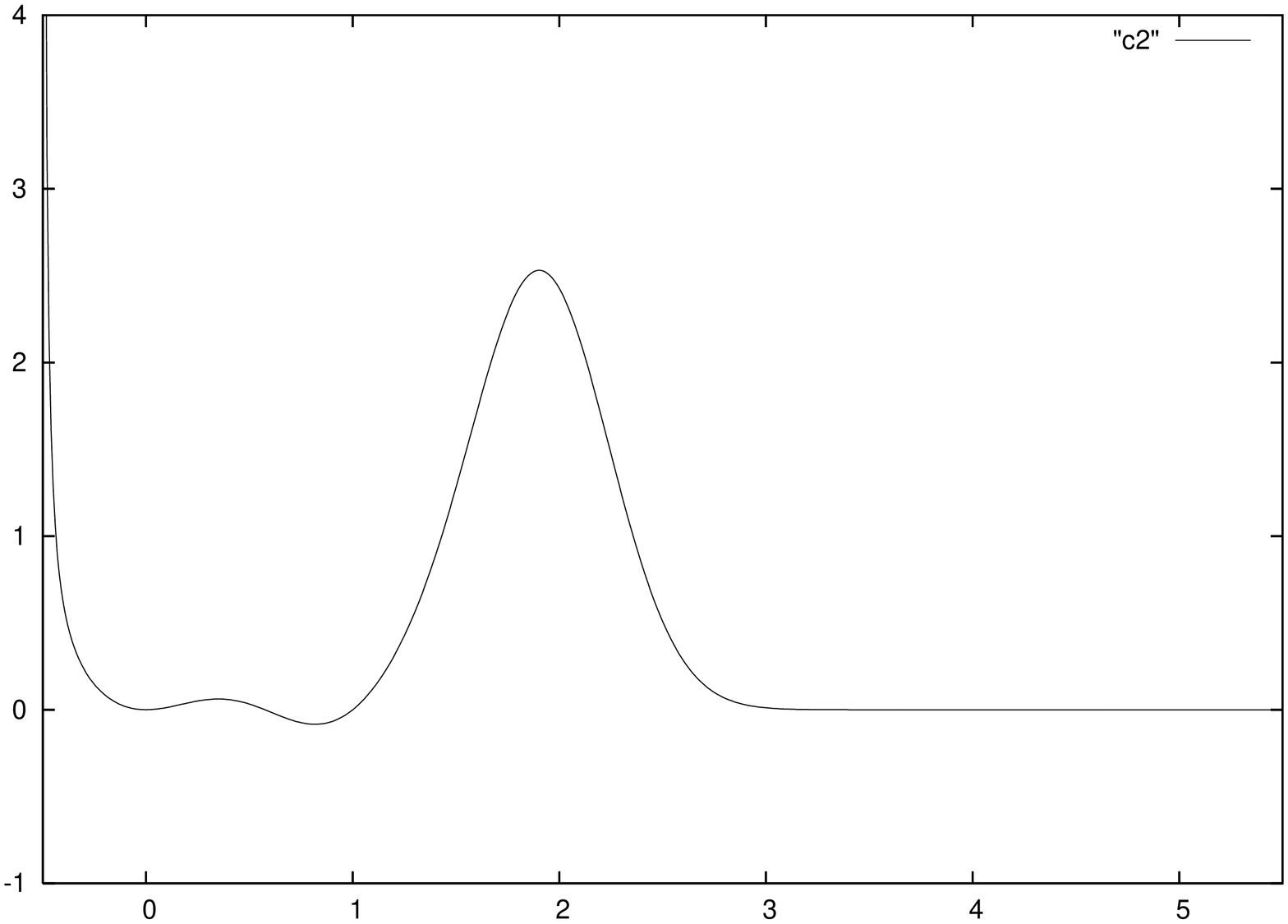,width=2.75in,angle=270}
    \psfig{figure=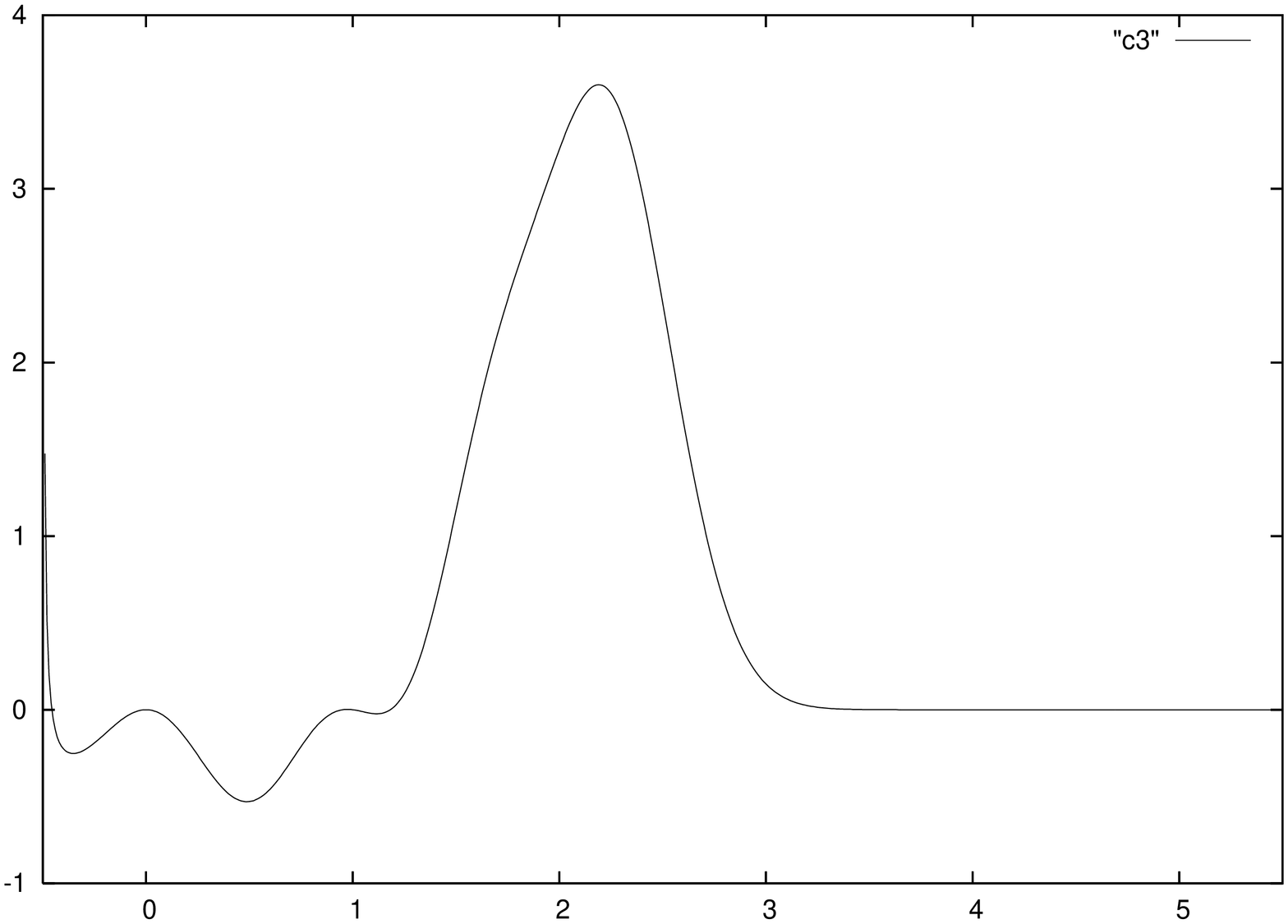,width=2.75in,angle=270}
}
\centerline{
    \psfig{figure=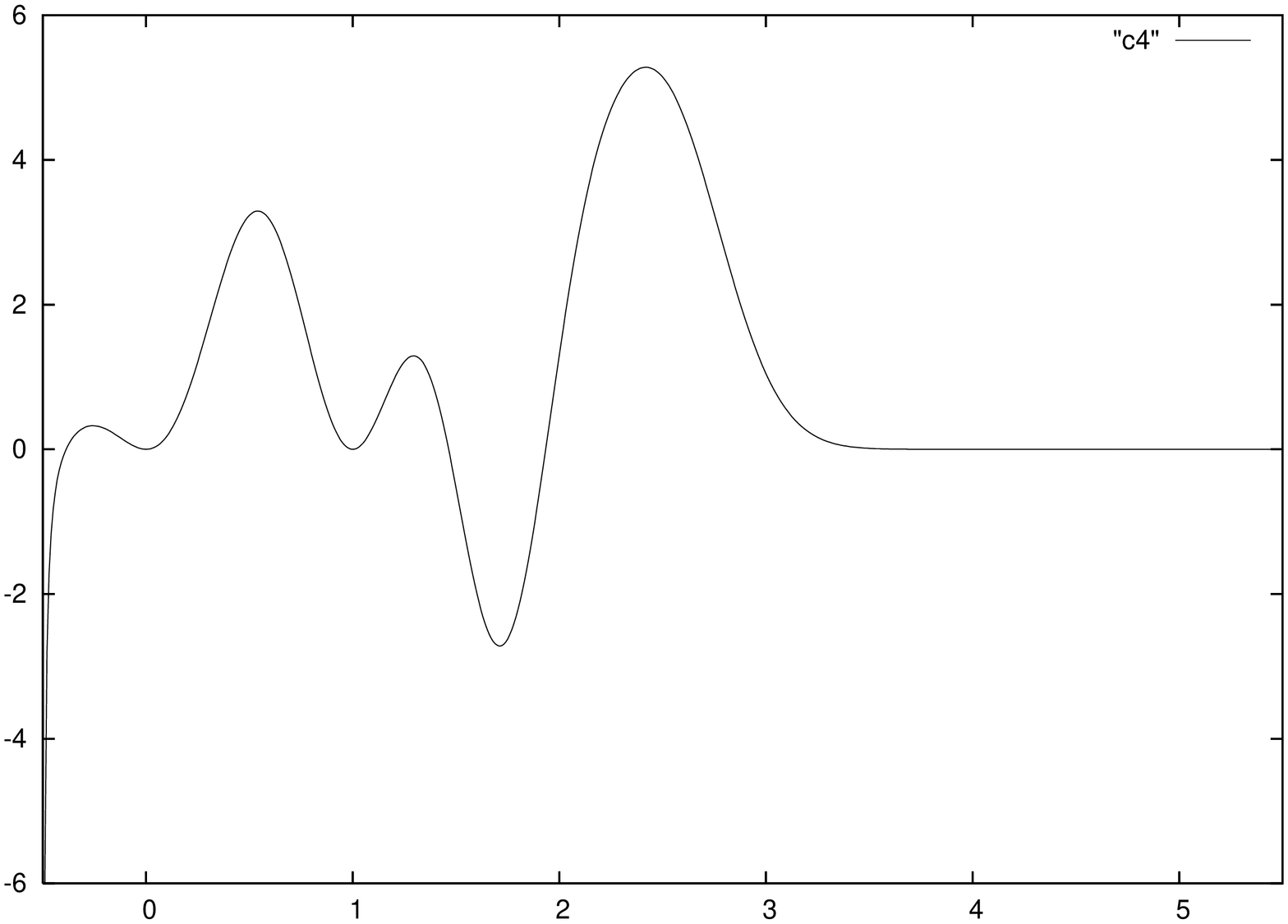,width=2.75in,angle=270}
    \psfig{figure=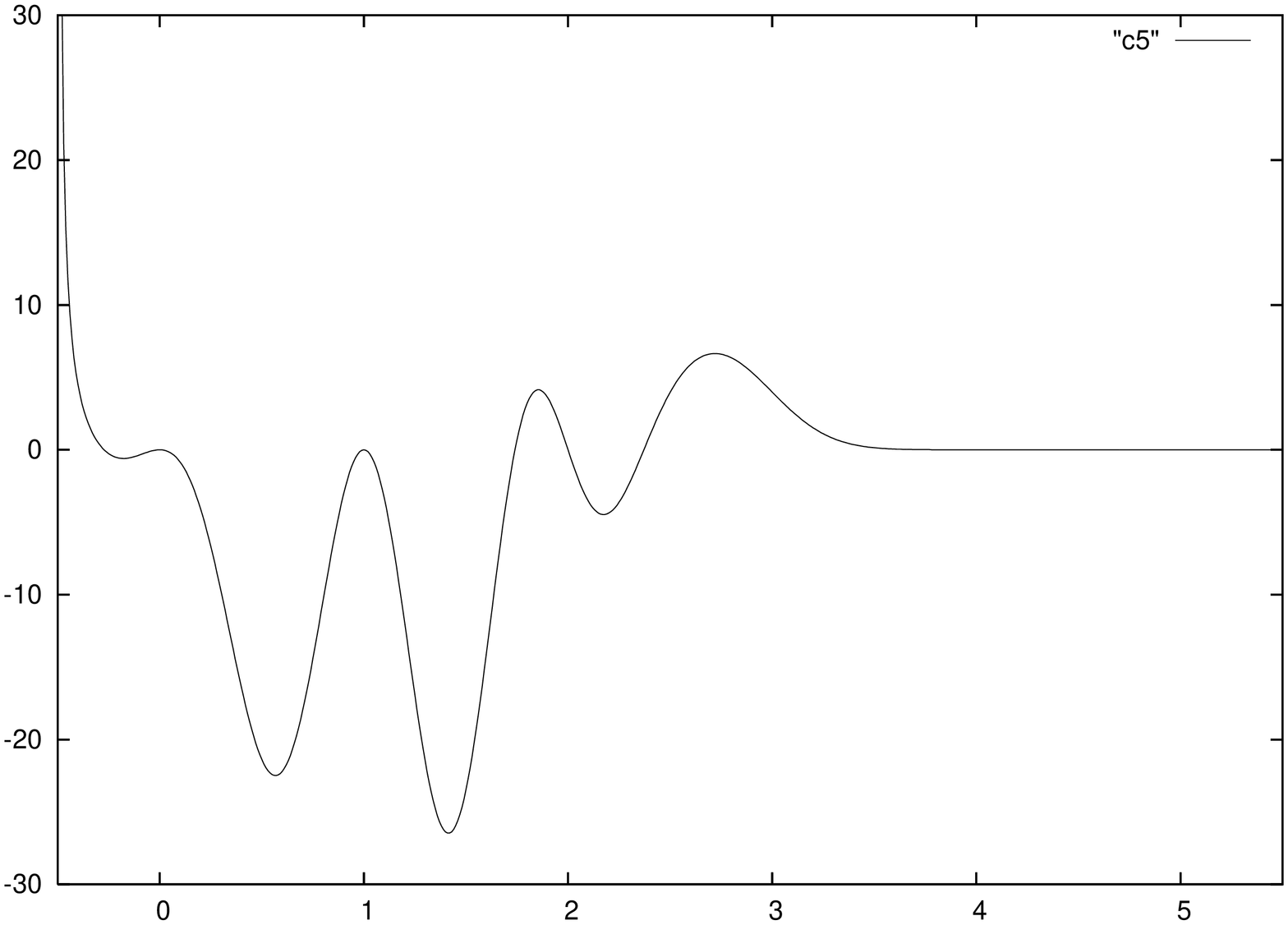,width=2.75in,angle=270}
}
\centerline{
    \psfig{figure=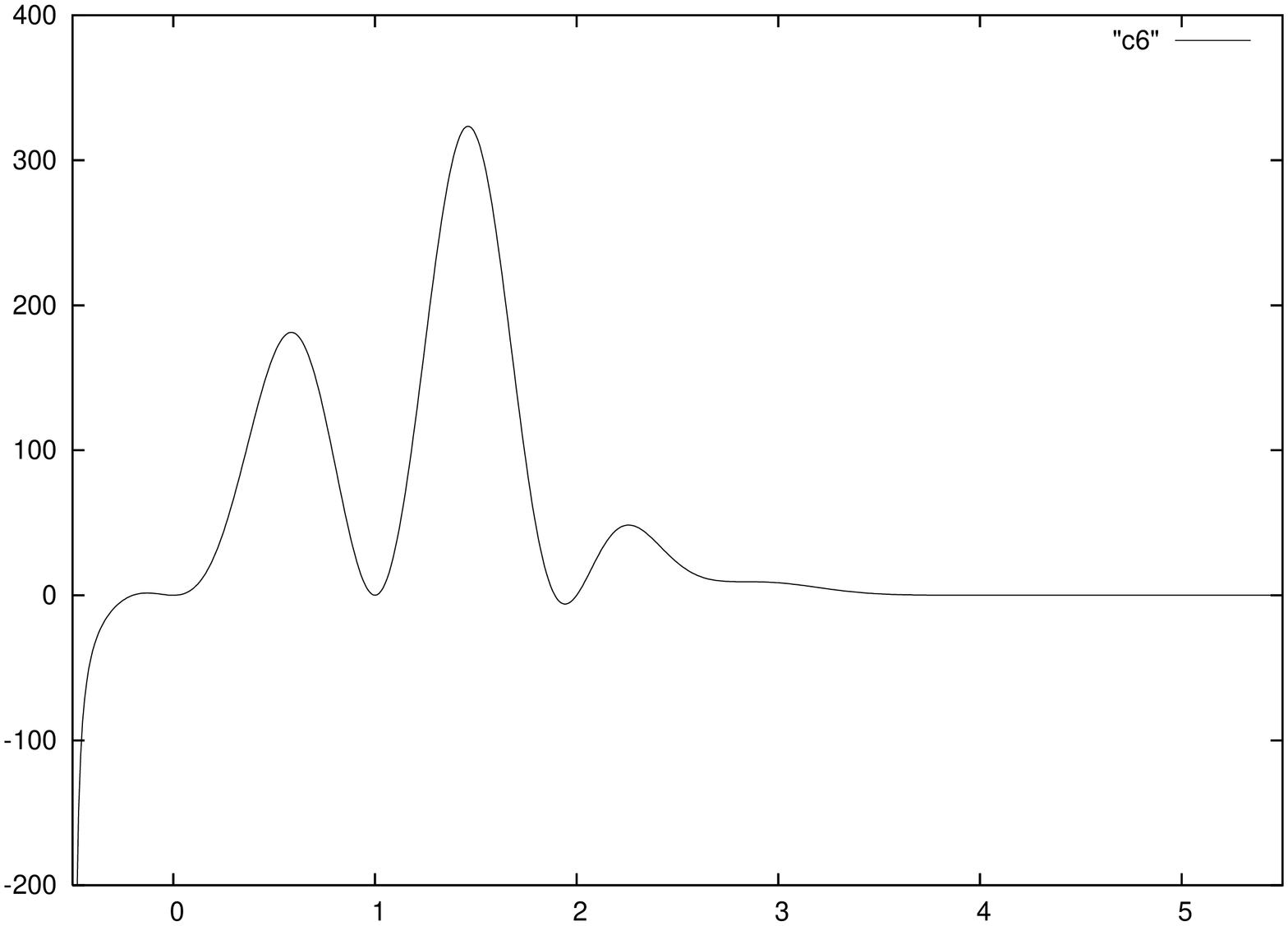,width=2.75in,angle=270}
    \psfig{figure=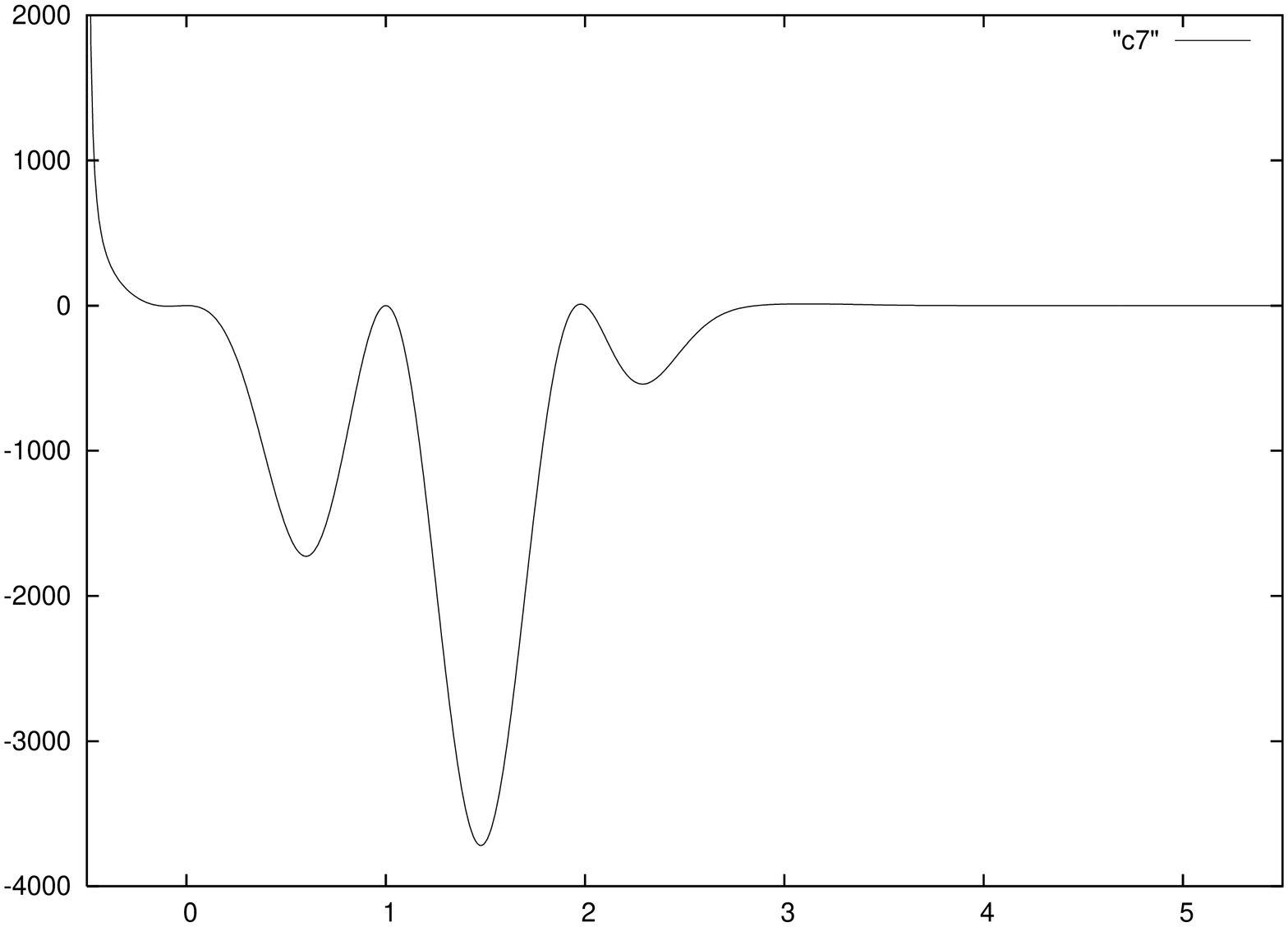,width=2.75in,angle=270}
}

\caption
            {Graphs of $c_r(k)$ with $-1/2 < k < 11/2 $, for $r=0\ldots,7$.}
\label{fig:graphs c_r}
\end{figure}

\begin{figure}
\centerline{
    \psfig{figure=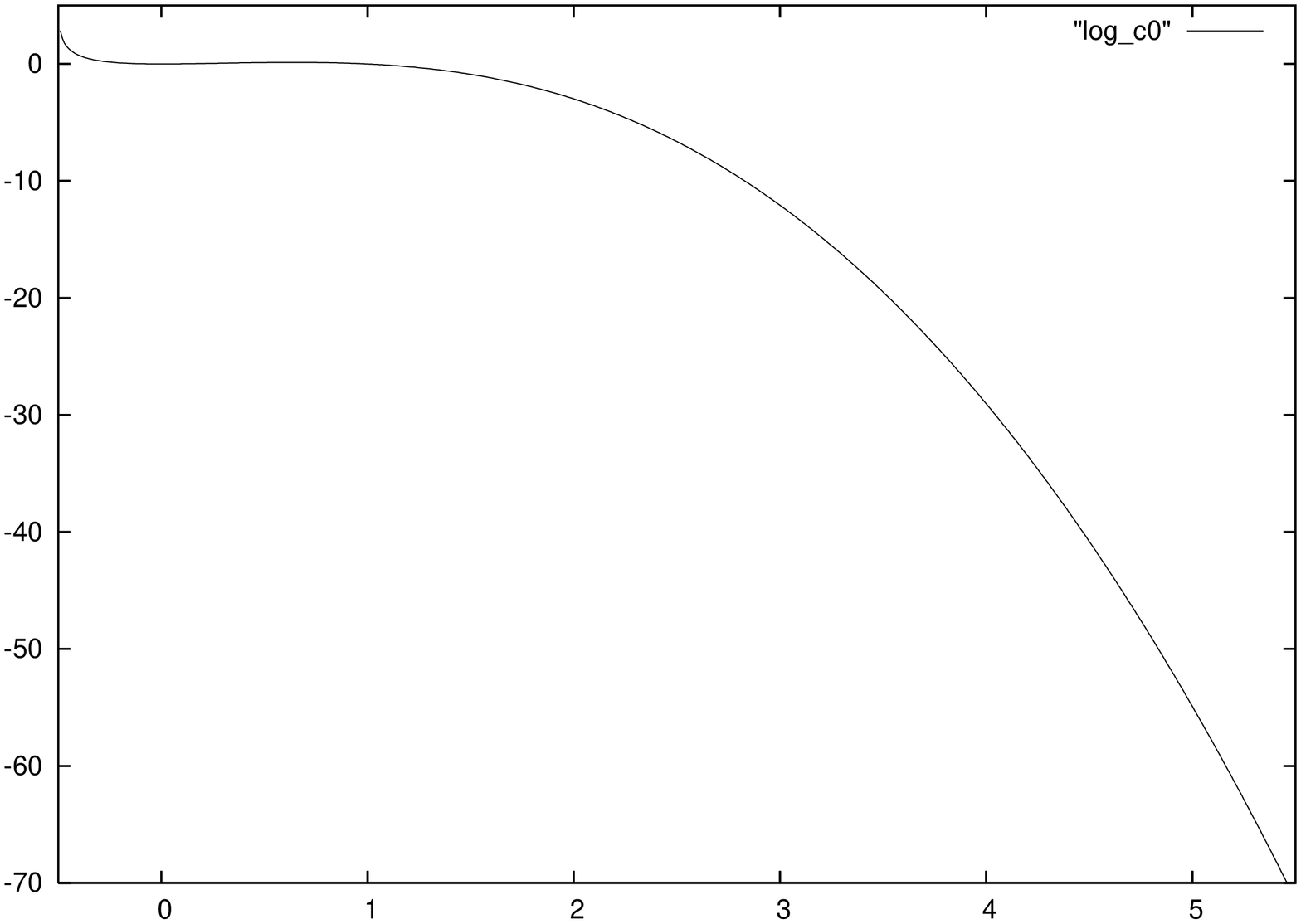,width=2.75in,angle=270}
    \psfig{figure=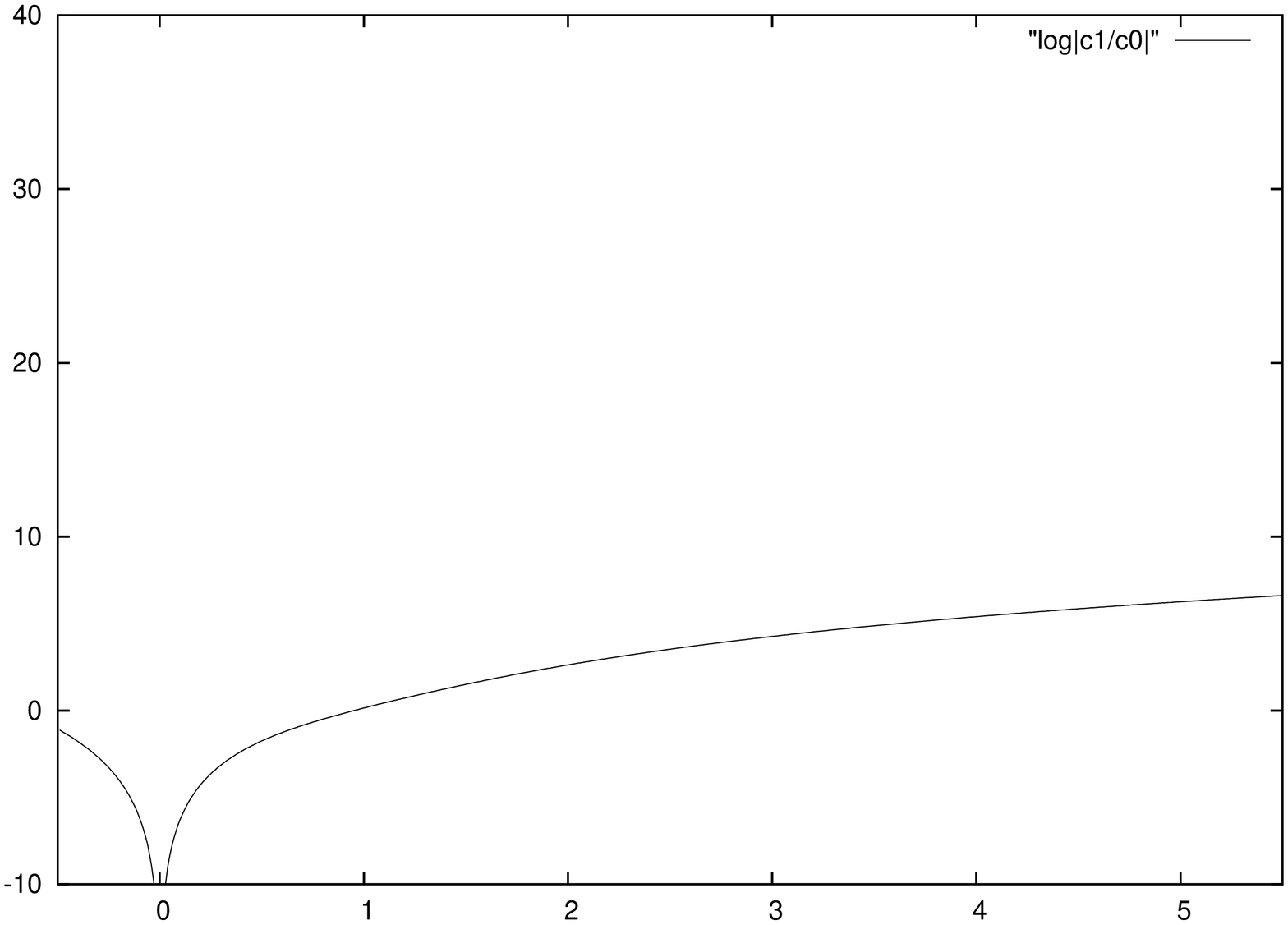,width=2.75in,angle=270}
}
\centerline{
    \psfig{figure=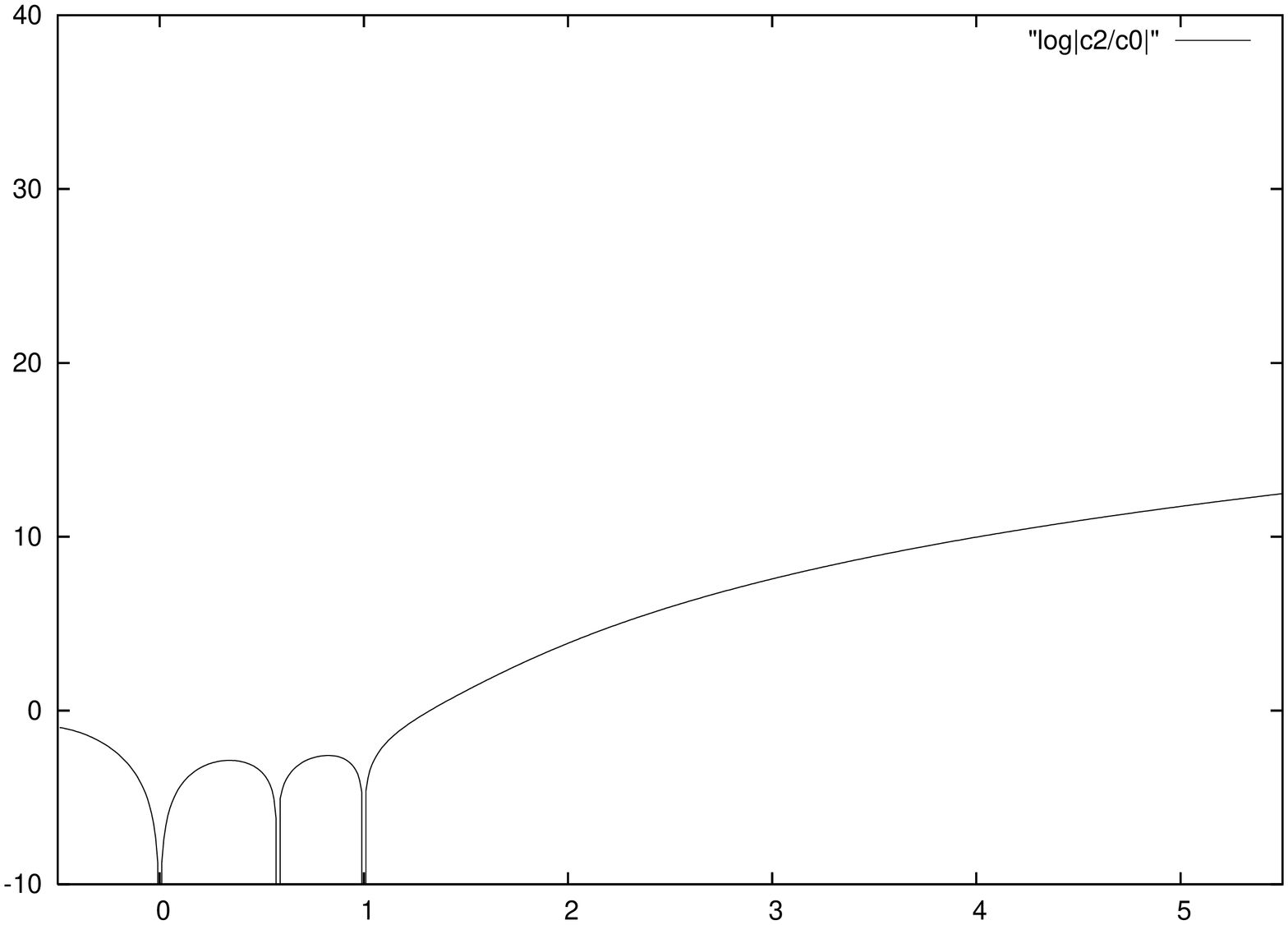,width=2.75in,angle=270}
    \psfig{figure=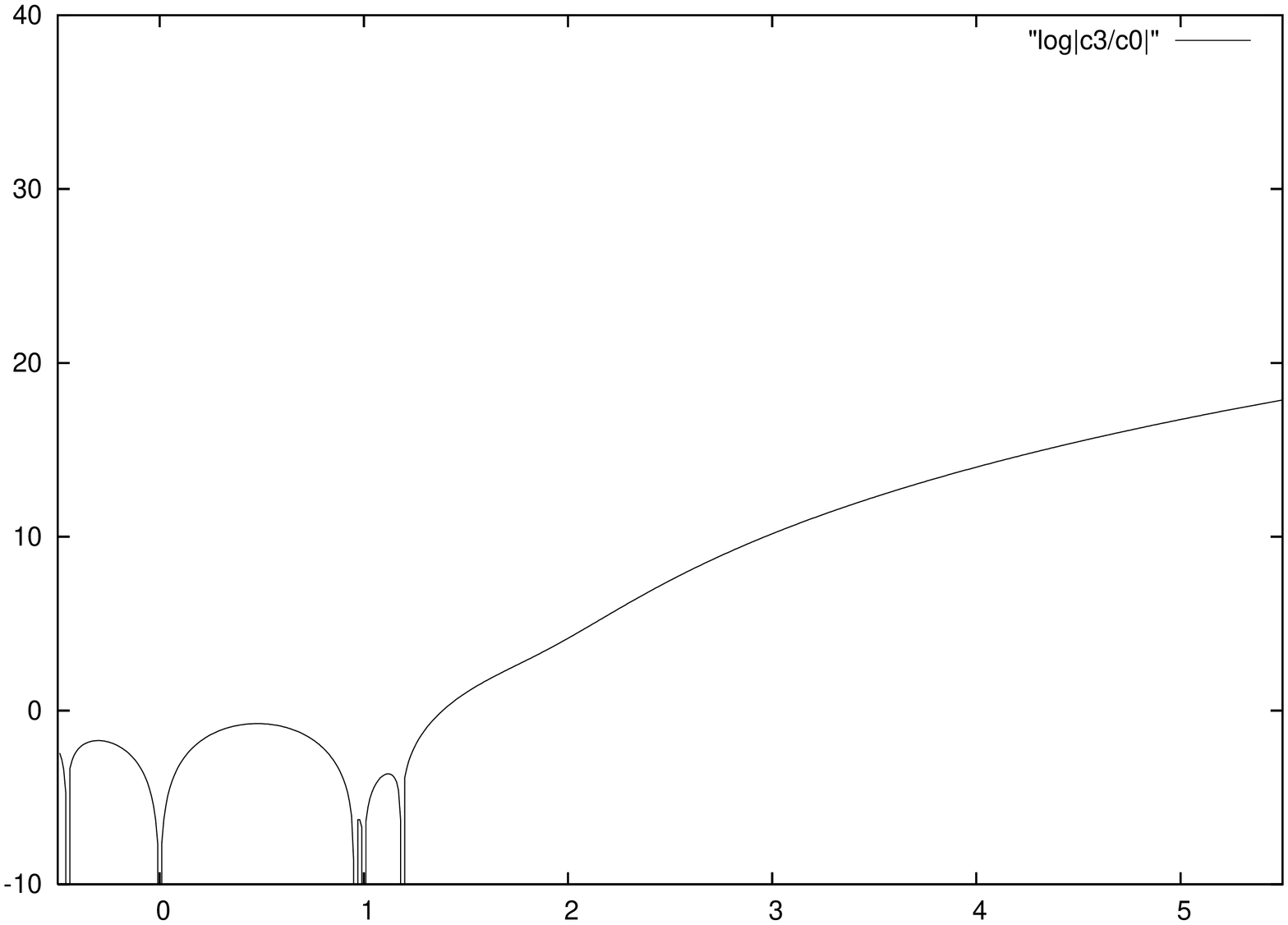,width=2.75in,angle=270}
}
\centerline{
    \psfig{figure=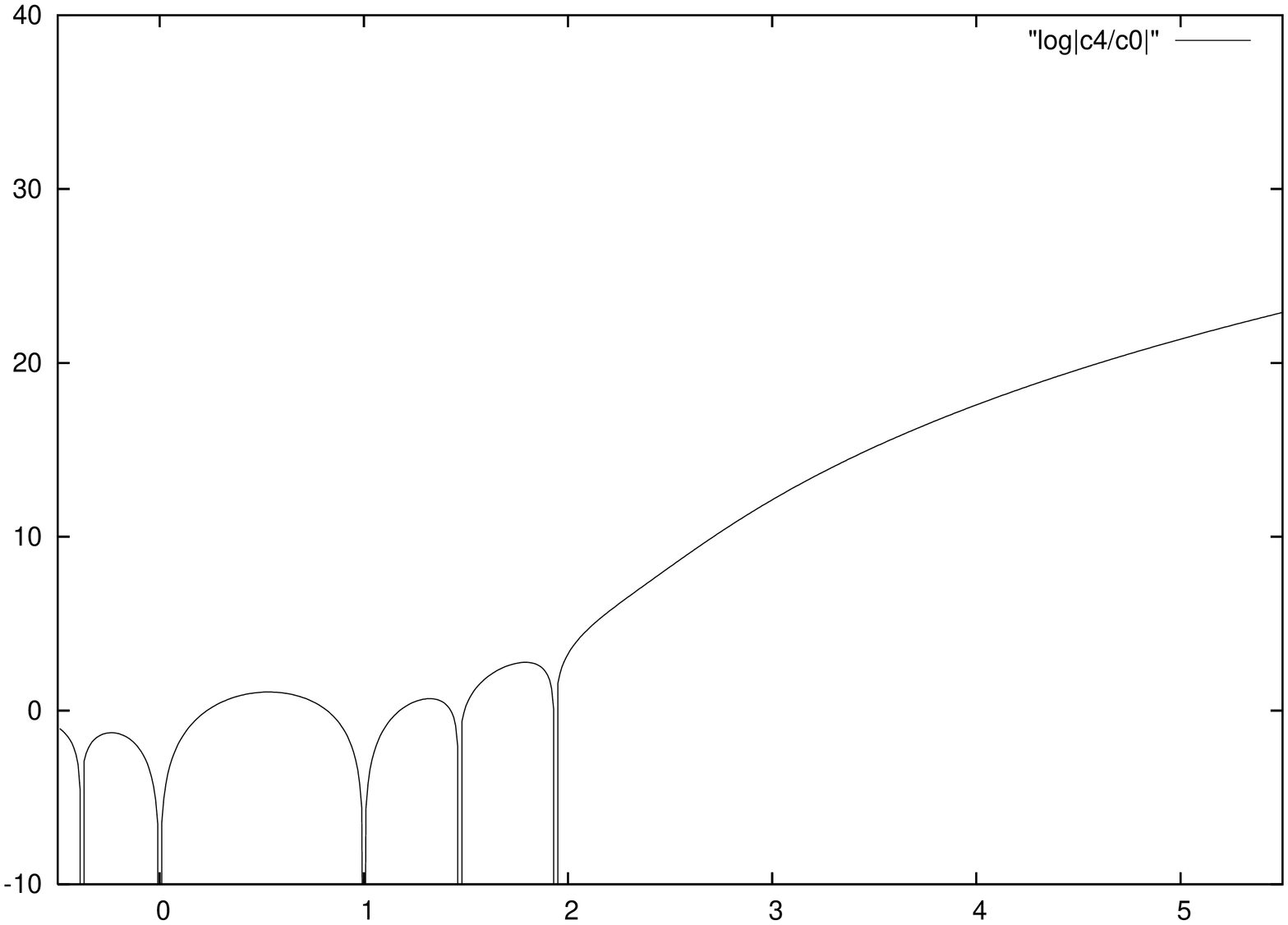,width=2.75in,angle=270}
    \psfig{figure=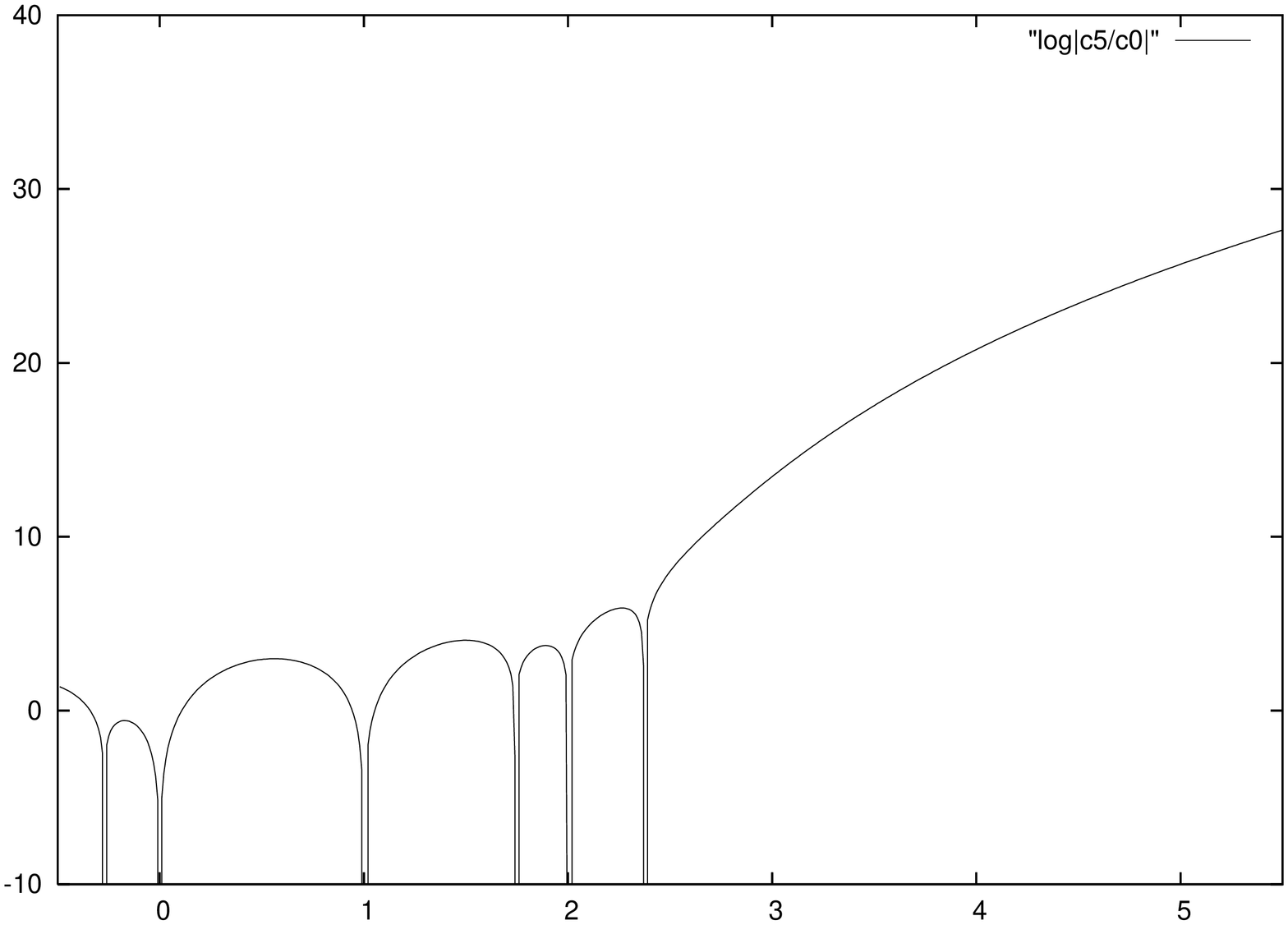,width=2.75in,angle=270}
}
\centerline{
    \psfig{figure=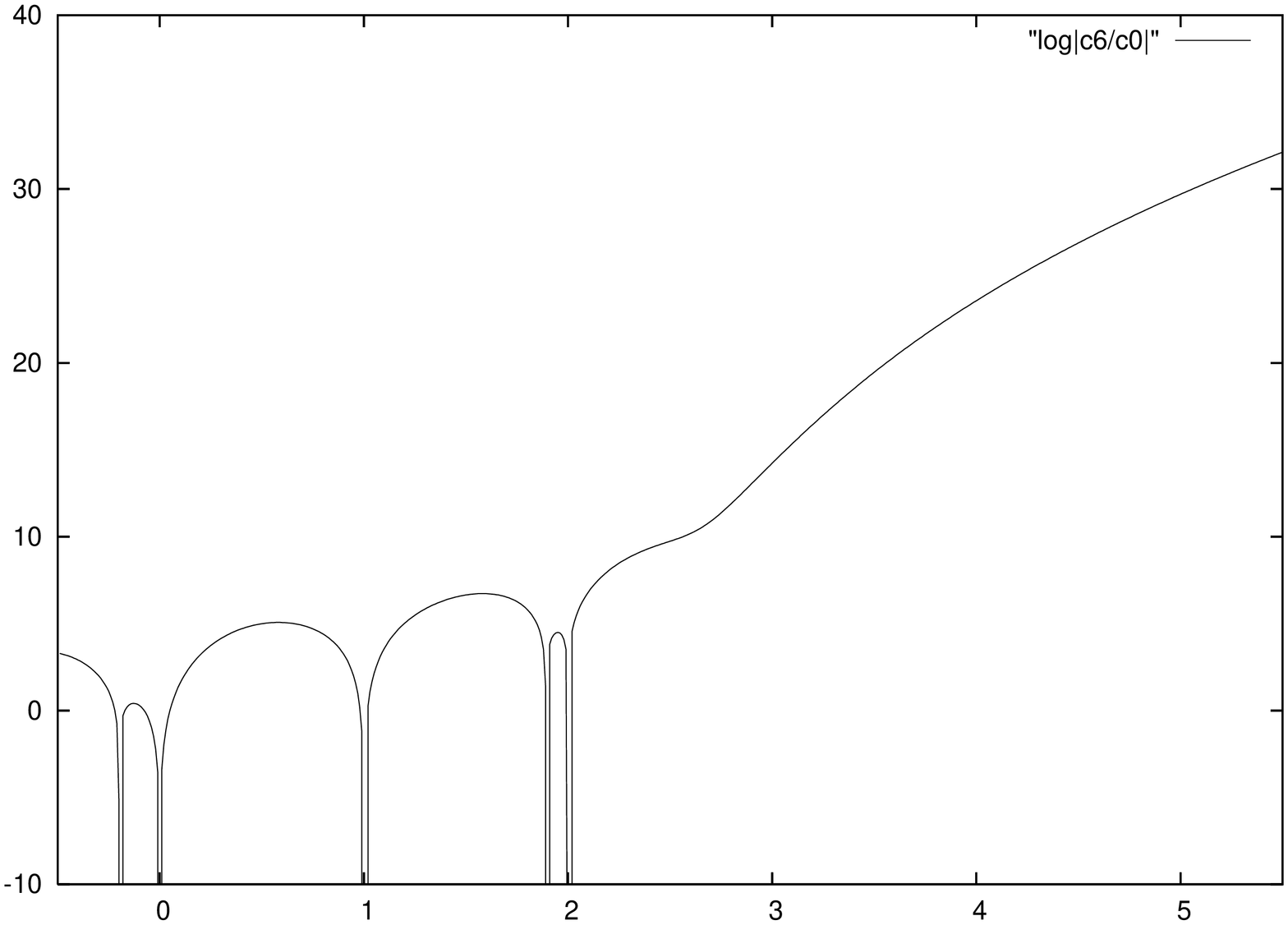,width=2.75in,angle=270}
    \psfig{figure=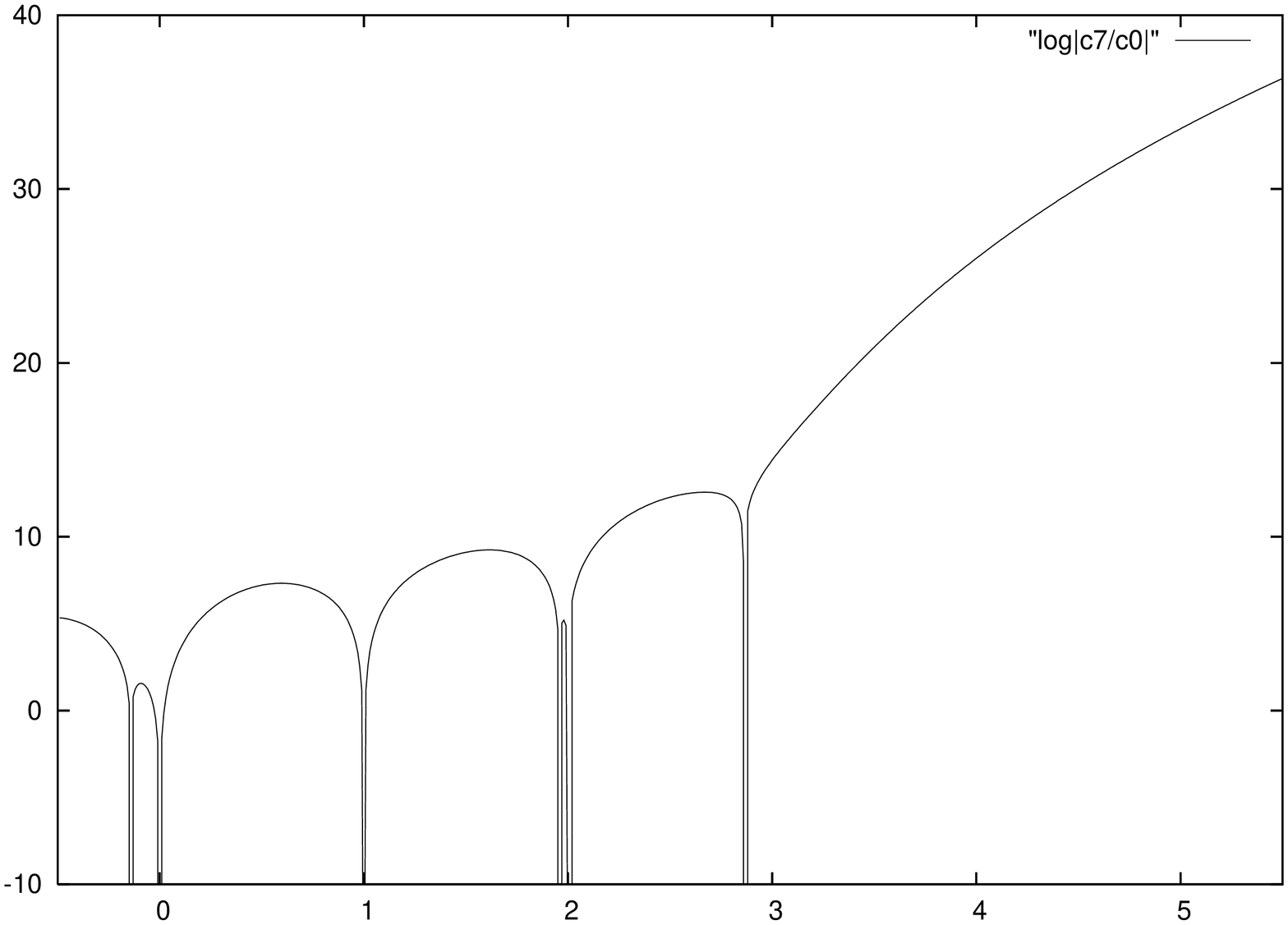,width=2.75in,angle=270}
}

\caption
            { The first figure
             depicts the graph of $\log(|c_0(k)|)$, while the next seven depict $\log(|c_r(k)/c_0(k)|)$, for $r=1\ldots,7$.
             The asymptotic behaviour of $\log(c_0(k))$ as $k\to \infty$ is implied by~\cite{CGo}~\cite{KeS} 
             and is, to leading order, $-k^2 \log(k)$. The cusps occur at zeros of $c_r(k)$, some of which are
             accounted for by the fact that, for non-negative $k \in {\mathbb Z}$,
             $P_k(x)$ is a polynomial of degree $k^2$ so that $c_r(k)=0$ if $r>k^2$.
             }
\label{fig:graphs c_r b}
\end{figure}

We present in Figures~\ref{fig:graphs c_r} and~\ref{fig:graphs c_r b}
some graphs of the coefficients $c_r(k)$ with $-1/2 < k < 11/2 $, for $r=0,1\ldots,7$.

\section{Remarks about other families: orthogonal and symplectic}
\label{section:other families}

In this paper we have explained two approaches for obtaining the coefficients
$c_r(k)$ of $P_k(x)$. The first involves explicitly determining the
residue on the r.h.s. of~(\ref{eq:P_k(x)}).
Theorems~\ref{thm:lot}--~\ref{thm:b_k}, and the procedure given in
Section~\ref{section: lower terms zeta} describe this in detail.
The second approach involves using the combinatorial sum~(\ref{eq: combinatorial sum}),
using small shifts, and high precision.

The same methods can be taken for other families of $L$-functions,
for instance in determining the lower order terms in the moments of
$L(1/2,\chi_d)$, quadratic Dirichlet $L$-functions,
or of $L_E(1/2,\chi_d)$, the $L$-functions associated to the quadratic twists
of a given elliptic curve, to name just two examples, in both cases
evaluated at the critical point. The former
is an example of a unitary symplectic family, while the latter is an
example of an orthogonal family~\cite{KaS}. See~\cite{CFKRS} where we discuss these
examples in detail. As with the Rieman zeta function, conjectures are given
for the full asymptotics of their moments, expressed in terms of multi-dimensional residues
and also as combinatorial sums. In that paper, we used Method~2
of Section~\ref{section:method2} for the analogous combinatorial sums
to numerically compute lower terms for the moments, and verify the full asymptotics.

For the elliptic curve family, the next to leading term in the asymptotics of the moments
has been worked out explicitly, and a test has been divised to
verify the first two terms in the asymptotics of that particular
family with an application to estimating the number of elliptic
curves of rank greater than zero~\cite{CPRW}. See also~\cite{BMSW} for a 
survey of results related to the latter question.

% ------------------------------------------------------------------------

\subsection*{Acknowledgment}

JBC, DWF, and MOR were supported by the NSF Focused Research
Group grant DMS 0244660. JPK was supported by an
EPSRC Senior Research Fellowship. MOR was also funded by NSERC.
NCS was supported by fellowships from the Royal Society and EPSRC.
JPK, MOR, and NCS wish to thank AIM for providing further support and an environment
for collaborative work.

% ------------------------------------------------------------------------
\bibliographystyle{amsplain}

% ------------------------------------------------------------------------

\end{document}